# Seiberg-Witten equations on three-manifolds with Euclidean ends

Yi-Jen Lee
Department of Mathematics, Princeton University
Princeton NJ 08544, U.S.A.
ylee@math.princeton.edu

This version: February 2002

# Contents

# 1 Introduction

This paper comprises of the foundational work for some versions of Seiberg-Witten theory on 3-manifolds with Euclidean ends. A *manifold with a Euclidean end, or an MEE for short*, is a smooth, orientable 3-manifold formed by connect-summing a compact closed manifold with $\mathbb{R}^3$, whose metric is Euclidean outside a compact region (cf. Definition 2.2.2 below). We consider Seiberg-Witten equations ((2.2) below) on such manifolds, with a family of perturbations parametrized by $t \in \mathbb{R}^+ \cup \{0\}$. See (2.4) below for the form of the perturbation 2-form $\omega$. Roughly speaking, it has $-t\theta/2$ as the dominant term, where $\theta$ is a harmonic 2-form asymptotic to $dx_1 \wedge dx_2$, $(x_1, x_2, x_3)$ being the coordinates of $\mathbb{R}^3$. In contrast to the well-known theory on *compact* manifolds, which is essentially independent of the choice of metric or perturbation, the theory on *non-compact* manifolds typically depends crucially on the asymptotic conditions. Indeed, in our theory the cases of $t = 0$ and $t > 0$ behave quite differently.

The following theorem summarizes our main conclusions.

**Theorem** *Let $M$ be a 3-manifold with a Euclidean end with a fixed $\mathrm{Spin}^c$ structure.*

*(a) When $t = 0$ and $b_1(M) > 1$, for a generic pair of metric $g$ and closed 2-form $w$ on $M$, the moduli space $\mathcal{M}_{g,w}$ of Seiberg-Witten solutions is a compact 0-dimensional manifold, and for two such pairs $(g_1, w_1)$, $(g_2, w_2)$, the moduli spaces $\mathcal{M}_{g_1,w_1}$, $\mathcal{M}_{g_2,w_2}$ are cobordant.*

*(b) When $t > 0$, the moduli space $\mathcal{M}_{g,w,t}$ of Seiberg-Witten solutions for a generic pair $(g, w)$ is a disjoint union of finite-dimensional (possibly non-compact) manifolds. The dimensions of the different components depend on the vortex number $n \in \mathbb{Z}^+ \cup \{0\}$. Letting $\mathcal{M}^n_{g,w,t} \subset \mathcal{M}_{g,w,t}$ denote the subset of vortex number $n$, there is a smooth, proper map from $\mathcal{M}^n_{g,w,t}$ to $\mathrm{Sym}^n \mathbb{R}^2 \times \mathrm{Sym}^n \mathbb{R}^2$. Furthermore, every point of vortex number $n$ in the end of the moduli space has a neighborhood with a product structure $\mathcal{N} \times \mathcal{R}$, where $\mathcal{N}$ is a neighborhood in $\mathcal{M}^p_{g,w,t}$, $p < n$, and $\mathcal{R}$ is a neighborhood in $\mathrm{Sym}^{n-p} \mathbb{R}^2$. Finally, for two arbitrary $t_1, t_2 > 0$, $\mathcal{M}_{g_1,w_1,t_1}$ and $\mathcal{M}_{g_2,w_2,t_2}$ are cobordant for generic pairs $(g_1, w_1)$, $(g_2, w_2)$.*

## 1.1 Motivation

The main motivation of this work is two-fold:

(A) THREE-DIMENSIONAL VERSIONS OF TAUBES'S WORK ON $SW = Gr$.

Taubes's work on the equivalence of 4-dimensional Seiberg-Witten invariants and the Gromov invariant of symplectic manifolds [42, 43, 44, 46] is arguably the most influential result of Seiberg-Witten theory. Very roughly, the idea is to consider Seiberg-Witten equations perturbed by two-forms of the form $t\omega$ (up to an inessential term), where $\omega$ is the symplectic form on the underlying symplectic 4-manifold, and $t$ is a large positive real parameter. As $t \to \infty$, the zero locus of the Higgs field of a Seiberg-Witten solution approaches a set of pseudo-holomorphic curves in the symplectic manifold, which the Gromov invariant counts. Conversely, given a suitable set of pseudo-holomorphic curves in the symplectic manifold, one may construct a Seiberg-Witten solution for large $t$ by "grafting" vortices along the pseudo-holomorphic curves. To extend this equivalence to general 4-manifolds with $b_2^+ > 0$, one may take $\omega$ to be a self-dual harmonic 2-form. (Symplectic forms, together with suitable metrics, are examples of such 2-forms.) Generically, such a 2-form vanishes along a set of circles in the manifold. The generalized Gromov invariant, which should be equivalent to the Seiberg-Witten invariant, should then count pseudo-holomorphic curves ending at these circles. In a series of articles [47, 48, 49], Taubes obtained some partial results implementing this generalization; the full equivalence, however, has not been established. In fact, even the "generalized Gromov invariant" mentioned above is yet to be defined. One of the main difficulties is to understand the behavior of pseudo-holomorphic curves or Seiberg-Witten solutions near the circles of degeneracy.

The 3-dimensional analogue of this generalization appears more accessible. Consider Seiberg-Witten equations on *closed* 3-manifolds $X$ perturbed by a large harmonic 2-form $t\omega$. Generically, $\omega$ may be chosen such that its dual 1-form $*\omega$ is Morse, with a equal number of index-1 and index-2 critical points. Taubes's argument in this case concludes that the zero loci of the Seiberg-Witten solutions converge to suitable sets of (finite length) gradient flows of $*\omega$ ending at the critical points. Due to this simpler description, we were able to define in [12, 13] an invariant $I_3$ that counts flows of a Morse 1-form, and conjectured based on Taubes' philosophy that:

**1.1.1 Conjecture** [12, 13] *Let $X$ be a closed 3-manifold with $b_1 > 0$ and orientation $o$. Then the 3-dimensional Seiberg-Witten invariant $\mathbf{Sw}_{X,o} = \pm I_3$.*

We refer the reader to section 4 of [12, 13] for the precise definition of $I_3$, a brief explanation of the ideas behind, and the notation used here. In the $b_1 = 1$ case, $\mathbf{Sw}_{X,o}$ is the Seiberg-Witten invariant in the "Taubes chamber" (i.e. large perturbation).

The counting invariant $I_3$ is shown (in a more general setting) in [12, 13] to be equivalent to a version of Reidemeister torsion $T(X, o)$ due to Turaev:

**1.1.2 Theorem** [13] *Under the same assumptions, $I_3 = \pm T(X, o)$.*

The present work may be viewed as the first step towards a proof of Conjecture 1.1.1 via Taubes's ideas, and this together with Theorem 1.1.2 would provide a geometric proof of the equivalence:

**1.1.3 Theorem** (Meng-Taubes) [28, 52] $\mathbf{Sw}_{X,o} = \pm T(X, o).$

In [28], the equivalence is established by checking that both invariants satisfy the same set of axioms, which depends most importantly on the surgery formulae of Seiberg-Witten invariants [50].

On the other hand, Donaldson proposed a scheme of proving (an averaged version of) Conjecture 1.1.1 via a topological field theory formulation of the 3-dimensional Seiberg-Witten invariant, which was recently implemented by T. Mark [27].

Among the three approaches to the Seiberg-Witten–torsion correspondence, the geometric picture following Taubes' approach mentioned above is conceptually most direct, though technically most challenging. However, another advantage of this approach is that the geometric picture has (relatively) simple extensions to Floer theory. Motivated by this picture, Hutchings and Thaddeus [14] define a "periodic Floer homology" for mapping tori of surface automorphisms, which is supposed to correspond to the Seiberg-Witten-Floer homology for such 3-manifolds. In this case, one may suppose that the harmonic 1-form $*\omega$ has no critical points, and the sets of flows counted by $I_3$ above consist only of periodic orbits. The chain groups of this "periodic Floer homology" are precisely free modules spanned by these sets of periodic orbits, and the boundary maps are defined by counting pseudo-holomorphic curves in $X \times \mathbb{R}$, with symplectic form $\omega + *\omega \wedge dt$. On the other hand, with a less straightforward twist of the same geometric picture, the author has a scheme of establishing the equivalence of the Seiberg-Witten-Floer homology (hopefully for general $X$) with the recently discovered Ozsvath-Szabo Floer homology [33], which is more computable. The details of of this have to be described elsewhere [24]; at present it seems accessible.

Here is roughly how the present work relates to the proof of Conjecture 1.1.1. A basic building block in Taubes' work is the "local model" for grafting, which is obtained by classifying Seiberg-Witten solutions on $\mathbb{R}^4$ (with suitable perturbations) and understanding their behavior. In the symplectic 4-manifold case, Taubes found them to be basically (up to gauge transformations) "pull-backs" of vortex solutions over $\mathbb{R}^2$. We obtain the analogous result in §3.2: The Seiberg-Witten solutions on $\mathbb{R}^3$ with perturbation $*dx_3$ are basically "pull-backs" of vortex solutions along the projection $\mathbb{R}^3 \to \mathbb{R}^2 : (x_1, x_2, x_3) \mapsto (x_1, x_2)$. Since over a closed 3-manifold $X$, $*\omega$, wherever it is nonzero, locally approximates $*dx_3$, this means that in the 3-dimensional case, the correspondence should be obtained by grafting vortices along flow lines of the dual vector field of the Morse 1-form $*\omega$. (Note that $x_3$ corresponds to the direction of the dual vector to $*\omega$.) More importantly, since $*\omega$ in general has critical points (i.e. where it vanishes), we also need a local model for grafting Seiberg-Witten solutions near these points. Such a local model is obtained in this paper, by understanding the Seiberg-Witten solutions on $\mathbb{R}^3$ with perturbation $*tdf$, where $f$ is an "admissible" Morse function with a pair of canceling critical points. (See Definition 2.2.6 for the definition of admissibility.)

In fact, to produce such a local model is the main reason why in this paper we choose to consider perturbations of the specified asymptotic behavior. We remark that such a local model is missing in the more difficult 4-dimensional situation, and it is one of the main obstacles for establishing the generalized $SW = Gr$ correspondence proposed by Taubes.

(B) NEW GAUGE-THEORETIC INVARIANTS OF 3-MANIFOLDS.

Instead of just $\mathbb{R}^3$, we consider the more general 3-manifolds with Euclidean ends (MEE). The motivation is a traditional one: to obtain invariants of 3-manifolds out of moduli spaces of solutions to a PDE. Since an MEE $M$ is a connected sum $M = X \# \mathbb{R}^3$, where $X$ is closed and $\mathbb{R}^3$ is endowed with the Euclidean metric, the gauge-theoretic invariant obtained may be regarded as an invariant of the closed 3-manifold $X$. In this paper we concentrate on the structure-theoretic aspects of the moduli spaces and leave the properties of such invariants for future investigation; so here we shall only very briefly indicate why this might be worth pursuing.

First, the better-known (and simpler) Seiberg-Witten invariant for *compact* 3-manifolds has been well-studied (see e.g. [25, 26, 28, 1]); however it turn out to be equivalent to previously-known invariants. For homology 3-spheres, it has been shown (with a modification due to Kronheimer to make it metric-independent) to coincide with the Casson invariant. For manifolds with $b_1 > 0$, we have already mentioned that it is equivalent to a version of Reidemeister torsion. Our version of Seiberg-Witten theory for MEE's however looks very different and might lead to something new. When $t$ is large, the philosophy described in part (A) above leads one to expect that the associated invariants should be computable via Morse-theoretic methods.

Second, even if as a *topological* invariant, the gauge-theoretic invariant obtained is not new, it might still have less straightforward applications to the geometry or topology of 3-manifolds or 4-manifolds. There are plenty of such applications in other versions of gauge theories. For example, via Weitzenböck-type formulae, the 4-dimensional Seiberg-Witten equations give various curvature estimates, and hence constraints on the existence of Einstein metric. (See e.g. [22].) The Seiberg-Witten theory for (closed) 3-manifolds or 4-manifolds has various applications to minimal-genus problems ([18]). Fintushel-Stern's construction [4] of exotic homotopy K3's was based on the Meng-Taubes theorem on the equivalence of 3-dimensional Seiberg-Witten invariants and Reidemeister torsion. In the non-compact situation, Floer proposed studying Yang-Mills-Higgs equations over asymptotically flat 3-manifolds, and was able to recover a theorem of Schoen and Yau on the geometry of such 3-manifolds [5].

## 1.2 Outline and summary of results

Below is a section-by-section outline of this paper which also serves as a brief summary of the results obtained. Due to their technical nature, here we can not be very precise in the statement of these results; the reader may find the precise statements in subsequent sections. Some non-technical explanation of the key points will be presented in §1.3.

- Section 2 contains the setup: basics of Seiberg-Witten theory, the definition of MEE and admissible 2-forms or functions, and some basic tools of analysis on such manifolds.

- In §3.1, we define "admissible configurations", which are basically $L^2_{2,loc}$-configurations with some weak assumptions on their asymptotic behavior. We shall always work with such admissible configurations. Some basic properties of these objects are discussed in §3.1.

- In §3.2, we completely solve the Seiberg-Witten equations on $\mathbb{R}^3$ with Euclidean metric and the standard perturbation $-\frac{t}{2}dx_1 \wedge dx_2$, where $t > 0$. The solutions consist of "pull-backs" of solutions to the vortex equations on $\mathbb{C}$. This leads to the description of the moduli space of Seiberg-Witten solutions on $\mathbb{R}^3$ as $\coprod_{k>0} \mathrm{Sym}^k(\mathbb{R}^2)$. We remark that this is another instance of the well-known relationship between Seiberg-Witten equations and the vortex equations: this has been used in the computation of Seiberg-Witten invariants for Kähler surfaces, and plays a crucial role in Taubes's proof of $SW = Gr$ explained above. In the 3-dimensional situation, this relationship is used to compute the Seiberg-Witten-Floer homology of $\Sigma \times S^1$ (with respect to $\mathrm{Spin}^c$-structures pulled back from $\Sigma$, a closed surface), which is the first nontrivial Seiberg-Witten-Floer homology computed. Generalizing this relationship, Mrowka-Ozsvath-Yu [32] computed the Seiberg-Witten-Floer homology for Seifert-fibered 3-manifolds.

- In §3.3, we prove that admissible Seiberg-Witten solutions on an MEE approximate solutions on $\mathbb{R}^3$ polynomially outside a compact region. This enables use to define the notion of the "vortex number" of a solution as the vortex number of the limit. The estimate is crucial for the Fredholm theory and gluing.

- In section 4, we define the configuration space $\mathcal{C}$ and construct its quotient space $\mathcal{Q}$ under the gauge group action. We show that $\mathcal{Q}$ has a Banach manifold structure. This requires more care than the usual gauge theories as the base manifold is not compact. We find the appropriate Banach spaces that work as the domain of the elliptic operators in our theory, and restrict our attention to subsets of the total configuration space containing the moduli spaces, that are Banach manifolds modeled on these Banach spaces. The cases of nontrivial vortex numbers cause additional complications by the fact that the curvature $F_A$ is not $L^2$. This difficulty is overcome by subtracting off fixed configurations and we show eventually that in these cases the quotient space has the structure of a fibered space, with Sobolev-space fibers.

- In section 5 we formulate the Fredholm theory of the relevant deformation operator $\mathcal{D}_c$. We find the appropriate Banach spaces as the domain and range of $\mathcal{D}_c$. The case when $t = 0$ or when $t > 0$ and the vortex number is zero is dealt with in §5.2, while the case of non-zero vortex numbers is discussed in the rest of the section. These two cases are completely different (cf. the discussion in §1.3). In the first case, the index of $\mathcal{D}_c$ is zero, while in the second case the index grows arbitrarily large with the vortex number.

- In §6.1 we follow the standard procedure to show that the moduli spaces are smooth finite-dimensional manifolds and enjoy some invariance properties by cobordism arguments. When the moduli spaces are non-compact, which is often the case in the situation considered in this paper, these standard invariance properties are useful (e.g. for the purpose of defining gauge-theoretic invariants) only when we also have an explicit description of the ends of the moduli spaces. In the rest of section 6, we give a recursive description of these ends: namely, the ends of higher vortex numbers are described in terms of moduli of lower vortex number numbers. For example, the end of the vortex number 1 piece of the moduli space consists of finite copies of $\mathbb{R}^3 \backslash B(R)$, where $B(R)$ denotes a 3-ball of large radius. This enables one to define 3-manifold invariants using these moduli spaces, in spite of their non-compactness. (Cf. [23] section 7).

## 1.3 Methods and historical background

Here we shall attempt to offer a brief explanation how our theory differs from the more familiar methods of gauge theory. The "philosophical" remarks here should become clear in later sections.

The more-familiar type of gauge theory on *non-compact* manifolds is that on cylindrical manifolds $\Sigma \times \mathbb{R}$, where $\Sigma$ is compact, or more generally on manifolds with cylindrical ends. This situation arises naturally in Floer theory, and in product or surgery formulae for gauge-theoretic invariants. Though analysis of this type can sometimes be rather complicated, the basic framework has by now become standard since the work of Floer, Taubes, and Mrowka. A standard reference for this type of techniques is [30]. Among other things, in this situation the gauge-theory equation is regarded as a formal $L^2$-gradient flow of a functional on the configuration space over $\Sigma$, and one attempts to establish an exponential decay towards the critical manifolds via estimates involving the functional and the gradient. This allows one to set up a Fredholm context of the theory via exponentially weighted Sobolev norms.

In the Euclidean situation considered in this paper, the story is quite different. There has been historically less study of gauge theories in this situation; the only references known to the author are papers by Taubes and Floer on Yang-Mills-Higgs theory on asymptotically Euclidean 3-manifolds in the 80s. Particularly relevant to our work are [39, 6, 7]. Aside from motivation from physics, a major reason for studying Yang-Mills-Higgs theory on *Euclidean* manifolds was the triviality of the theory on *compact* manifolds (which was also partly what motivated us). More recently, Kronheimer and Mrowka studied Seiberg-Witten theory on 4-manifolds with conical ends, where the conical ends are "symplectizations" of contact 3-manifolds [20].

Very roughly, since the usual Sobolev spaces typically do not work when the underlying manifold is non-compact, the key to such theories is to find suitable normed spaces such that the relevant differential operators in the theory have the desired invertibility or Fredholm properties. Furthermore, the norms have to be coarse enough so that the normed-spaces may include the moduli spaces, but they also have to be fine enough so that the nonlinear part of the gauge-theory equations is suitably small with respect them. How these norms should be chosen depends crucially on the asymptotic behaviors of the configurations in the theory. We find that in our $t = 0$ case, the configuration is asymptotically trivial, which is also the case with the Yang-Mills-Higgs theory studied by Taubes and Floer. Thus a suitable modification of the norms in [6, 7] works in our situation. Basically, the relevant operators in these case look like $-\nabla^2$ or $\bar{\partial}$ asymptotically, and the domain and range spaces are essentially completions of $C_0^\infty$ with respect to the norm $\|\nabla \cdot\|_2$ and $L^2$. When $t > 0$, in the case of zero vortex number the Higgs field is asymptotically a non-zero constant; in this case the basic model for the relevant operators is $-\nabla^2 + C$, $C > 0$, or a similar Dirac-type operator. This is the case when the usual Sobolev spaces work, and it is also the situation that occurs in [20]. The case when both $t$ and the vortex number are positive is the most complicated. We know that the configurations are asymptotically pull-backs of vortex solutions (§3.3), and the basic model for the relevant operator is $-\nabla^2 + V$, where $V$ is a function which is almost constant except for a few "tunnels" in the $x_3$-direction. The more familiar exponentially-weighted Sobolev norms actually work for the Fredholm theory here; however they are either too fine to describe the moduli spaces, or too coarse for the nonlinear aspects. The main reason is, even though we take the metric to be Euclidean except for a compact set, the associated admissible 2-form (which is used as our perturbation) in general approximates $t * dx_3$ only *polynomially* at infinity. As a consequence, the Seiberg-Witten solutions approxi-mates pull-backs of vortex solutions only *polynomially* as well. We therefore need polynomially weighted norms. The Fredholm theory is reduced to the Fredholm theory of a standard $\mathbb{R}^3$ case by a typical excision argument. On $\mathbb{R}^3$, we observe following [43] that the deformation operator may be decomposed into two parts: one roughly looks like $\partial/\partial x_3$, and the other, $iN'$, is self-adjoint and depends only on $x_1, x_2$. The relevant sections over $\mathbb{R}^3$ may be regarded as functions of $x_3$ taking values in the space of sections over $\mathbb{R}^2$ ($\mathbb{R}^3 \to \mathbb{R}^2 : (x_1, x_2, x_3) \mapsto (x_1, x_2)$). This can be decomposed via the decomposition of the above space of values into Coker $N'$ and its $L^2$-orthogonal complement. (Coker $N'$ is nontrivial precisely when the vortex number is non-zero.) The deformation operator preserves this decomposition. Over the part taking value in the orthogonal complement, it is invertible between polynomially-weighted norms of the same weights. Over the part taking value in Coker $N'$, the deformation operator looks like $\partial/\partial x_3$ and is Fredholm between polynomially weighted norms of *different* weights. Our choice of norms is a combination of the two. The parameter for the weight, $\epsilon$ below, is in this case taken to be a number between 1 and 3/2: it has to be larger than 1 because the weight for the range has to be one-less than that for the domain, while still has to be positive; it also has to be smaller than 3/2, because the associated space has to contain the moduli space, which can be guaranteed by our decay estimate (Proposition 3.3.3) only when $\epsilon > 3/2$. Long after this work was completed, we discovered while working on another project that Floer has used similar norms for pseudo-holomorphic disks ending at degenerate Lagrangian intersection points in [8].

On the other hand, the Fredholm theory in the $t = 0$ case and the vortex number zero case when $t > 0$ works by a modification of the techniques in [39], which is also what was done in [20].

We also briefly comment on the proof of the decay estimates (Proposition 3.3.3). Again it is very different from the cylindrical situation; it is however simpler and uses the nice pointwise estimates via maximum-principle-type arguments that are particular to Seiberg-Witten theory. This is also the case with [20]. However, though the proof of the analogous decay estimate in [20] has the same starting point as ours, what they did was to first estimate all invariant quantities (under gauge group action) via maximum-principle-type arguments, and then choose a gauge with nice asymptotic behavior. In comparison, we only use maximum-principle-type arguments to estimate the magnitude of $\beta$ and its derivatives. To get estimates for other quantities, we use the $\mathbb{R}$-action on $\mathbb{R}^3$ which is intrinsic to our situation. Namely, we use the Seiberg-Witten equations in a "temporal" gauge to relate the derivative of these quantities in the $x_3$-direction to the magnitude of $\beta$ and its derivatives; then integrate over $x_3$. (Cf. §3.3. step 3).

Finally, some history about this work itself. The author began working on this project since the summer of 1995; at that time there was no literature on Seiberg-Witten theory on 3-manifolds, and we had to start from scratch. The case $t > 0$, vortex number zero, which is the simplest case in this paper and is in several places similar to the analysis in [20], was done long before the preprint of [20] was available. In the long process of revising this paper, we have tried to eliminate any part that is now in the existing literature. Except for abridgement, expository and stylistic improvements, there has been basically no change in mathematical content since the version of 1998. That this work has taken so long to see publication is due to a combination of many unfortunate factors.

To shorten this paper, we have in several places omitted standard arguments; the reader may consult the earlier version of this paper [23] for more details. We have also omitted the section that defines gauge-theoretic invariants from the moduli spaces obtained here, and the proof of Proposition 2.2.7. They were respectively, section 7 and Appendix 1 in [23].

### 1.4 Notation and conventions.

The reader is advised to first browse through this subsection for a guide to the conventions, then return later for reference of notation.

Throughout this paper we let $C$, $C'$, or $C_i$, $i \in \mathbb{Z}^+$ denote positive constants of order 1 varying with the context. Similarly, $\varepsilon, \varepsilon', \varepsilon_i$ will denote some small positive constants.

*In contrast, the plain epsilon $\epsilon$ parameterizes the weight on the norms, $\epsilon \in [0, 3/2)$.*

$\langle \cdot, \cdot \rangle$ or $\langle \cdot, \cdot \rangle_2$ usually denotes the $L^2$-product of two functions, while $(\eta, \chi)$, or $\eta \cdot \chi$ denotes the pointwise inner product of the functions $\eta$ and $\chi$. $|\cdot|$ denotes the Euclidean norm of a vector in $\mathbb{R}^n$.

We denote the trivial $\mathbb{R}$-bundle or $\mathbb{C}$-bundle on a manifold by $\mathbb{R}$ or $\mathbb{C}$.

$\Gamma(M, E)$ denotes the space of sections of the bundle $E$ over the manifold $M$. $L^2(M, E)$ denotes the space of $L^2$ sections of $E$. Similar notation is used for for the completions with respect to other norms. When $M$ is the 3-manifold with a Euclidean end under discussion in this paper, we often omit $M$ in the notation, and denote the spaces as $\Gamma(E)$ etc. $\Omega^k := \Omega^k(M)$ is the space of $k$-forms on $M$.

Unless otherwise specified, any norm in this paper is a norm of functions or sections on the 3-manifold $M$. We use notation such as $\|\cdot\|_{L^2_k(\mathbb{C})}$ to denote norms of functions on a different space (in this case, $\mathbb{C}$).

Regarding the different notions of adjoint operators used in this paper: we shall use $D^*, D^\dagger, D^{\mathsf{t}}$ to denote, respectively, the formal $L^2$-adjoint, formal $L^2_\epsilon$-adjoint, and the adjoint in the sense of [16] of the differential operator $D$.

We shall often omit subscripts or superscripts when there is no danger of confusion.

The objects in our theory (e.g. configuration spaces, moduli spaces) are typically indexed by $t$, which is the real parameter in the perturbation (2.4); $n$, which is the vortex number, and $l$, which indicates the differentiability.

Other notation or conventions will be specified along the way.

**Acknowledgments.**

The author would like to thank Prof. Taubes for suggesting this problem, for his indispensable help and careful reading of the manuscripts of this work. She also thanks P. Feehan, P. Kronheimer, P. Ozsvath, M. Shubin for helpful discussions/comments.

## 2 Preliminaries

This section contains the set-up of this work.

### 2.1 Review of Seiberg-Witten theory on 3-manifolds

Here we quickly review elements of 3-dimensional Seiberg-Witten theory, and fix some notation to be used. Some references for the theory on *compact* 3-manifolds or 4-manifolds are [28, 25, 29].

Seiberg-Witten theory is concerned with an oriented 3-manifold $M$ endowed with a fixed Spin$^c$-structure $s$. (All orientable 3-manifolds are Spin$^c$ [17]). A "spinor bundle" $S$ associated to a Spin$^c$-structure is a rank 2 hermitian bundle. An "associated line bundle" is $L := \det S$.

A pair $(A, \psi)$ is usually called a *configuration*, where $A$ is a unitary connection on $L$ (or equivalently a Spin$^c$ connection on $S$—*we shall often confuse the two*), and $\psi$ is a section of $S$.

Given a harmonic 2-form $\theta$ (which will be our perturbation later on), its action on the spinor bundle by Clifford multiplication has eigenvalues $\pm i|\theta|$; away from the zeroes of $\theta$, this gives rise to a splitting of the spinor bundle into a pair of complex line bundles

$$S = E \oplus E', \text{ where } E' := E \otimes K^{-1}. \tag{2.1}$$

Here $K^{-1}$ is the sub-bundle of $TM$ (away from the zeroes of $\theta$) whose fibers consist of tangent vectors annihilated by $*\theta$. $\theta$ endows $K^{-1}$ with a complex structure $J$: For two tangent vectors $w$, $v$ in a fiber of $K^{-1}$, $\langle w, Jv \rangle := \frac{\theta}{|\theta|}(v, w)$. We shall often identify $K^{-1}$ with the $\overline{\text{Hom}}_{\mathbb{C}}(K^{-1}; \mathbb{C})$ component of $T^*M \otimes \mathbb{C}$.

Wherever the splitting (2.1) makes sense, we denote by $A^E$ the connection on $E$ induced by $A$. Given $\theta$ and a metric on $M$, $A$ can be equivalently specified by $A^E$. In general, with respect to the decomposition (2.1), the Spin$^c$ connection will have off-diagonal components depending on the metric and $\theta$.

On 3-manifolds the Seiberg-Witten equations read

$$\begin{cases} \not{\partial}_A \psi = 0 \\ \rho(F_A) = i\sigma(\psi, \psi) + i\rho(\omega), \end{cases} \tag{2.2}$$

where $\partial_A$ is the Dirac operator and $\sigma$ as usual denotes the map $\mathbb{C} \times \bar{\mathbb{C}} \to \mathfrak{su}_2$

$$(v, w) \mapsto i\Big(\frac{wv^\dagger + vw^\dagger}{2}\Big)_0, \qquad (2.3)$$

where the subscript 0 denotes the traceless part. And for $F \in \bigwedge^2 T^*M$, $\rho(F)$ stands for the $\mathfrak{su}_2$ representation of $F$ via its action on the spinor bundle by Clifford multiplication. In 3-dimensions, elements in $\bigwedge^1 T^*(M)$ act by Clifford multiplication on the same rank two spinor bundle. We will denote the representation by the same notation $\rho$, taking $\rho(dx_1 \wedge dx_2) = \rho(dx_3)$.

In these equations, $\omega$ is a closed 2-form usually called "perturbation". In this paper it takes the form:

$$\omega = -\frac{t}{2}\theta + w, \qquad (2.4)$$

where $\theta$ is an "admissible form" defined below, and $w$ is small in the sense of section 6. *We take $w = 0$ in sections 2-5 of this paper unless otherwise specified.* In this paper, a "Seiberg-Witten solution" always refer to a solution of (2.2) with perturbation (2.4) on an MEE, which will be defined next.

Solutions of the Seiberg-Witten equations can be regarded as minima of the "energy functional":

$$\mathcal{E}(A, \psi) = \int_M |\partial_A\psi|^2 + \frac{1}{2}\int_M |\rho(F_A) - i\sigma(\psi, \psi) - i\rho(\omega)|^2. \qquad (2.5)$$

## 2.2 Analysis on manifolds with Euclidean ends

Part (A) below defines MEE's and admissible metrics and 1-forms on them. Part (B) introduces some useful norms on MEE's and their basic properties.

(A) MEE'S AND ADMISSIBLE PAIRS ON THEM

**2.2.1 Notation** $B(R)$ denotes an open 3-ball of radius $R$; $D(R)$ denotes an open disc of radius $R$.

**2.2.2 Definition** A *3-manifold with a Euclidean end*, $(M, g)$, (or a "MEE" for short) is a complete, orientable manifold with a metric $g$ in $L^2_{l,loc}(\text{Sym}^2 T^*M)$ for some $l \geq 4$ such that there is a $\Re \in \mathbb{R}^+$ and an injective smooth map

$$\Omega : \mathbb{R}^3 \backslash B(\Re) \to M,$$

so that: (1) $M_\Re := M\backslash\Omega(\mathbb{R}^3\backslash B(\Re))$ is compact; (2) on $\Omega^{-1}(M\backslash M_\Re)$, $\Omega^*g - g_0 = 0$, where $g_0$ is the Euclidean metric on $\mathbb{R}^3$.

A metric satisfying the above is called an *(l-)admissible metric*.

Note that in the above $L^2_{l,loc}$ is defined with respect to a *fixed* smooth metric; however it is independent of the choice.

**2.2.3 Notation** We will often use the same notation to denote a function or a section on $\mathbb{R}^3 \backslash B(\Re)$ with its corresponding function or section on $M \backslash M_\Re$ induced by $\Omega$. For $x \in M \backslash M_\Re$, $|x|$ makes sense via $\Omega$ in a similar way. For any $R > \Re$, we define the open set $M_R := \{x \in M : |x| < R\}$.

**2.2.4 Notation** Let $\chi_1$ be a smooth cutoff function on $\mathbb{R}^3$ which is 1 for $|x| \leq 1$, and 0 for $|x| \geq 2$. Generalizing, we use $\chi_R$ to denote a cutoff function on $M$ such that $\chi_R(x) = \chi_1(\Omega^{-1}(x)/R)$ on $M \backslash M_\Re$, and 1 on $M_\Re$.

**2.2.5 Notation** Let $(x_1, x_2, x_3)$ be the Cartesian coordinates of $x \in \mathbb{R}^3$. Later we will combine $x_1, x_2$ as $z := x_1 + ix_2$, and write $x = (z, x_3)$. On $M \backslash M_\Re \sim \mathbb{R}^3 \backslash B(\Re)$, we will often decompose a configuration as: $(A, \psi) = ((A_z, A_3), (\alpha, \beta))$, where $\alpha$, $\beta$ are components in the splitting (2.1) of $S$.

**2.2.6 Definition** For any integer $l \geq 4$, and a fixed $l$-admissible metric $g$, we call a 2-form $\theta$ *(l-)admissible* (with respect to $g$) if it is $C^{l-3}$, harmonic (with respect to $g$), and over $M \backslash M_\Re$, $q := \Omega^*(*\theta) - dx_3$ satisfies:

$$\sum_{k=0}^{l-3} \left(|x|^{(3+k)} |\nabla^k q|\right) \leq C. \tag{2.6}$$

When $\theta$ is coexact with $*\theta = df$, we call the harmonic function $f$ *(l-)admissible* if over $M \backslash M_\Re$,

$$\sum_{k=0}^{l-2} \left(|x|^{(2+k)} |\nabla^k (\Omega^* f - x_3)|\right) \leq C'. \tag{2.7}$$

An *(l-)admissible pair* is a pair of an $l$-admissible metric and a corresponding admissible form.

We shall choose the constant $\Re$ in Definition 2.2.2 so that $C < 1$ in (2.6), and $C' < 1$ in (2.7).

In fact admissible pairs exist in plenty as the following proposition shows. Versions of this proposition have appeared in the literature (e.g. [10].) A proof can be found in Appendix 1 of [23].

**2.2.7 Proposition** *For an arbitrary admissible metric on a 3-manifold with a Euclidean end, there exists a unique corresponding admissible function. Furthermore, in either the $L_l^2$ or $C^\infty$ category, the admissible function corresponding to a generic admissible metric is Morse.*

We shall from now on until section 6 fix an admissible pair $(g, \theta)$ on an MEE $M$. Since $*\theta$ is always exact on $M \backslash M_\Re \subset \mathbb{R}^3$, there is a function $f$ over $M \backslash M_\Re$ satisfying (2.7). Such a function will suffice for our later purposes, and we shall also call such a function an admissible function, though it is not necessarily globally defined on $M$.

(B) SOME USEFUL BANACH SPACES ON MEE'S.

**2.2.8 Definition** Let $V$ be a Euclidean or hermitian bundle over an MEE $M$, and let $A$ be a metric-preserving connection on $V$. Define $L^p(V)$, $V^p_{k/A}(V)$, and $L^p_{k/A}(V)$ to be the completions of $C_0^\infty(V)$ (the space of compactly supported smooth sections), with respect to the following norms respectively:

$$\|\xi\|_p^p := \int_M (\xi,\xi)^{p/2},$$

$$\|\xi\|_{p;k/A} := \sum_{i=1}^k \|\underbrace{\nabla_A \cdots \nabla_A}_{i} \xi\|_p,$$

where $\nabla_A$ is a metric-preserving connection induced from $A$ and the Levi-Civita connection, and

$$\|\xi\|_{p,k/A} := \|\xi\|_p + \|\xi\|_{p;k/A}.$$

**2.2.9 Notation** Later on we will often omit the subscript $A$. For example, $\|\cdot\| = \|\cdot\|_{\cdot/A}$ if $A$ is imposed with certain asymptotic conditions such that the Banach spaces defined are independent of the different choices of $A$ (which is often the case). Furthermore, in these cases we have the inequality

$$C\|\xi\|_{\cdot/A} \leq \|\xi\|_{\cdot/A'} \leq C'\|\xi\|_{\cdot/A} \tag{2.8}$$

for different connections $A$, $A'$ satisfying the asymptotic conditions and some constants $C$, $C'$ depending on $A$, $A'$. (cf. Lemma 4.2.2, Lemma 4.3.3). Therefore, if an inequality involving $\|\cdot\|_A$ holds, a similar formula (which differ only in some constants) holds for $\|\cdot\|_{A'}$, for any connection $A'$ satisfying the same conditions. Another case to drop the subscript $A$ is when $V$ is trivial (with a fixed trivialization) or a bundle derived from $T^*M$; in this case $A$ is assumed to be obvious choice, namely, the trivial connection or the connection induced from the Levi-Civita connection.

$V_k^2$ will be often denoted by $V_k$.

**2.2.10 Definition** Let $\lambda_1'$ be a cutoff function on $\mathbb{R}$ so that $\lambda_1'(s) = 1$ for $|s| < 1$ and $\lambda_1'(s) = 0$ for $|s| > 2$. Let $\lambda_R'(s) := \lambda_1'(s/R)$.

Let $\varsigma$ be a real function on $M$ defined by:

$$\varsigma := \lambda_R' \circ f + (1 - \lambda_R' \circ f)|f|/R \tag{2.9}$$

In the above we take $R > \Re$ to be large enough such that $\|\nabla\varsigma\|_\infty \leq C/R \ll 1$, where $C$ is a positive constant.

**2.2.11 Definition** (weighted version of 2.2.8) Let $\xi, A$ be as in Definition 2.2.8 and let $\epsilon \in \mathbb{R}$, $\epsilon \geq 0$. We define the following norms:

$$\|\xi\|_{p:\epsilon} := \|\varsigma^\epsilon \xi\|_p.$$
$$\|\xi\|_{p;k:\epsilon/A} := \|\varsigma^\epsilon \nabla_A \xi\|_p + \cdots + \|\varsigma^\epsilon \nabla_A^k \xi\|_p.$$
$$\|\xi\|_{p,k:\epsilon/A} := \|\varsigma^\epsilon \xi\|_p + \|\varsigma^\epsilon \nabla_A \xi\|_p + \cdots + \|\varsigma^\epsilon \nabla_A^k \xi\|_p.$$

Define the weighted Banach spaces $L^p_\epsilon(V)$, $V^p_{k:\epsilon/A}(V)$, $L^p_{k:\epsilon/A}(V)$ as the completions with respect to the above three norms similarly to Definition 2.2.8.

Note that when $\epsilon = 0$, the norms and their corresponding Banach spaces reduce to the unweighted case in Definition 2.2.8. When $p = 2$, $L^p_\epsilon$ is equipped with a Hilbert space structure:

$$\langle \xi, v \rangle_{2:\epsilon} := \langle \varsigma^\epsilon \xi, \varsigma^\epsilon v \rangle_2$$

The Hölder inequality in the weighted norms takes the form

$$\|fg\|_{b:\epsilon} \leq \|f\|_{p:\epsilon_1} \|g\|_{q:\epsilon_2}, \tag{2.10}$$

for $b^{-1} = p^{-1} + q^{-1}$ and $\epsilon_1 + \epsilon_2 = \epsilon$. In particular, if $f \in L^p_\epsilon$ and $g \in L^q_\epsilon$, then the right hand side is bounded by $\|f\|_{p:\epsilon} \|g\|_{q:\epsilon}$.

We have the following versions of Sobolev inequalities on an MEE.

**2.2.12 Lemma** ([6], Lemmas 13 & 14, with errors corrected) *Let $V, A$ be as in Definition 2.2.8. Let $p \in (1,3)$ and $\hat{V}^p_{1/A}(E)$ be the space of sections $\xi$ of $V$ satisfying $\|\nabla_A \xi\|_p < \infty$. Then we have embedding $\hat{V}^p_{1/A}(V) \hookrightarrow L^q_{loc}(V)$ for $q \leq \bar{p} \equiv \frac{3p}{3-p}$. Moreover, there exists a constant $C_s$ depending on $p$ and $V$ but not on $A$, and a continuous map $\mu : \hat{V}^p_{1/A}(V) \to [0, \infty)$ that factors through the map $\xi \to |\xi|$, such that*

1. $\mu^{-1}(0) = V^p_{1/A}(E);$
2. $\||\xi| - \mu(\xi)\|_{\bar{p}} \leq C_s \|\nabla \xi\|_p.$ \hfill (2.11)

*If moreover $\|\nabla_A \xi\|_{p,1} < \infty$, then*

$$\lim_{R \to \infty} \sup_{|x| \geq R} (|\xi \circ \Omega| - \mu(\xi)) = 0. \tag{2.12}$$

*Also if $1 < p < 3 < q < \infty$, we have for any $\xi \in C_0^\infty(E)$,*

$$\|\xi\|_\infty \leq C \|\xi\|_{q,1/A}; \tag{2.13}$$
$$\|\xi\|_\infty \leq C'(\|\nabla_A \xi\|_p + \|\nabla_A \xi\|_q). \tag{2.14}$$

For example, when $\epsilon = 0$ and $E$ is a trivial real line bundle, $\xi$ may be identified with a scalar function $f$. In this case (2.11), (2.12) simply says that for all $\xi \in L^p_{1,loc}$ with $\nabla f \in L^p$, there exists a constant $c = \mu(f) \geq 0$, such that $f - c \in L^{\bar{p}}$ and such that $f \to c$ uniformly at infinity.

# 3 Properties of Seiberg-Witten solutions

In §3.1 we define admissible configurations and discuss some basic properties of admissible Seiberg-Witten solutions. In §3.2 we show that all admissible Seiberg-Witten solutions on $\mathbb{R}^3$ with the standard perturbation arise as pull-backs of vortex solutions up to gauge transformations. §3.3 contains a crucial decay estimate of Seiberg-Witten solutions.

## 3.1 Seiberg-Witten solutions on 3-manifolds with Euclidean ends

We mentioned that the Seiberg-Witten theories on MEEs corresponding to the $t = 0$ case and the $t > 0$ case are very different. The following observation is a first manifestation of this fact. When integration by parts is applicable, the Weitzenböck formula implies that in the $t = 0$ case, a finite energy solution has $L^2$-integrable $\nabla_A \psi$, $F_A$, and $|\psi|^2$. This is no longer true in the $t > 0$ case (cf. especially the $\mathbb{R}^3$ example in §2.2). Instead, we shall introduce configurations of finite "vortex numbers" for which $F_A$, $\nabla_A \psi$ are only $L^2$-integrable over generic hypersurfaces in $M$.

Let $A_0$ be a reference connection on $L$ such that $F_{A_0} \in C_0^\infty$, and let $(A, \psi)$ be a configuration on a MEE $M$. Then $A - A_0$ is a 1-form on $M$.

**3.1.1 Definition** In the above notation, a configuration $(A, \psi)$ is *admissible* if it satisfies:

1. $(A - A_0, \psi) \in L^2_{2/A_0, loc}$.

2. There exists a real number $R > 0$, such that $|\psi(x)|$ is bounded on $M \backslash M_R$.

3. If $t = 0$, then $(A - A_0, \psi) \in V_{2/A_0}$.
   If $t > 0$, then: (i) on $M \backslash M_\Re$, $\beta := (\frac{\rho(\theta)}{|\theta|} - i)\psi \in V_{1/A}(S)$; (Note that on $M \backslash M_\Re$, the splitting of $S$ (2.1) makes sense.); (ii) there exists a positive constant $R > \Re$ such that for any real number $C$, $|C| \geq R$, the integral of $|F_A|^2 + |\nabla_A \psi|^2$ over the plane $P_C$ in $M \backslash M_\Re \sim \mathbb{R}^3 \backslash B(\Re)$ given by $x_3 = C$ is finite.

Note that the above definition does not depend on the choice of $A_0$. In the $t = 0$ case, conditions 1 and 3 imply condition 2 by lemma 2.2.12. In the $t > 0$ case, condition 3 (ii) is an alternative way of saying that the configuration has finite "vortex number" (cf. Definition 3.3.8).

*In this paper, we shall always assume that the configurations are admissible.*

By the Sobolev embedding theorem, we have:

**3.1.2 Lemma** *If $(A, \psi)$ is an admissible configuration, then $\psi \in L^\infty$ for any $t$.*

**3.1.3 Lemma** *Suppose that $M$ is an MEE with a $(l+6)$-admissible metric; $l \geq 3$, and $w \in L_l^2$ (hence $\omega \in L_{l,loc}^2$ in (2.4)). Then the following holds for any admissible Seiberg-Witten solution*

1. $(A - A_0, \psi) \in C^2$.

2. *For any small $\varepsilon > 0$, there exists an $\varepsilon$-dependent constant $R > 0$, such that $|\psi|^2(x) \leq z + \varepsilon$ in the region where $|x| > R$, where $z := \|\omega\|_\infty$ in the $t > 0$ case, and $z := 0$ when $t = 0$.*

Note that conversely, statements 1 and 2 in this lemma obviously imply condition 1 and 2 in Definition 3.1.1.

*Proof.* The first claim follows from part (1) of the next lemma by the Sobolev embedding theorem. The second claim follows from Proposition 3.3.3 below in the $t > 0$ case, and in the $t = 0$ case it follows from part (2) of the next lemma and lemma 2.2.12. □

**3.1.4 Lemma** *Let $(A, \psi)$ be as in the previous lemma. Then (1) $(A, \psi) \in L_{l+1/A_0, loc}^2$; (2) in the $t = 0$ case, $(A, \psi) \in V_{l+1/A_0}$.*

*Proof.* Both (1) and (2) follow from an elliptic bootstrapping argument similar to that in [29] section 5.3, using the $L^\infty$-bound on $\psi$. □

Thus we have the following finer uniform $L^\infty$-bound on $\psi$ when $(A, \psi)$ is a Seiberg-Witten solution.

**3.1.5 Lemma** *Let $(A, \psi)$ be as in lemma 3.1.3 and $z' := \|sup(-s, 0)\|_\infty$, where $s$ is the scalar curvature of $M$. Then $|\psi(x)|^2 \leq z := \|\omega\|_\infty + z'$ for all $x \in M$. Via (2.2), this gives a $L^\infty$-bound for $F_A$.*

This follows easily from lemma 3.1.3 and a standard argument via a Weitzenböck formula (cf. [19]).

## 3.2 Solutions on $\mathbb{R}^3$

The theory on $\mathbb{R}^3$ is the simplest example and will be the building block for later sections. In this case the Seiberg-Witten equations can be completely solved.

First notice that when $t = 0$ and $w = 0$, (2.5), the Weitzenböck formula and Definition 3.1.1 tell us that the only solution (up to gauge equivalence) is the trivial one: $F_A \equiv 0$, $\psi \equiv 0$. So we quickly move on to the $t > 0$ case.

In the fundamental representation of $su(2)$, the Lie algebra is spanned by the basis [1]

$$\gamma_1 = \begin{pmatrix} 0 & i \\ i & 0 \end{pmatrix}, \ \gamma_2 = \begin{pmatrix} 0 & 1 \\ -1 & 0 \end{pmatrix}, \ \gamma_3 = \begin{pmatrix} i & 0 \\ 0 & -i \end{pmatrix}. \tag{3.1}$$

On $\mathbb{R}^3$, we choose the Euclidean coordinates and metric, and choose the perturbing harmonic 2-form to be $-t\theta/2$, $\theta := dx_1 \wedge dx_2 = *dx_3$. Since $|\theta|(x) > 0$ for all $x \in \mathbb{R}^3$, $\theta$ splits the spinor bundle into two trivial complex line bundles (cf. §2.1). Denote according to the splitting $\psi = (\alpha, \beta)$. The Euclidean coordinates and metric induce a canonical trivialization of $S$ (compatible with the above splitting), with respect to which the Spin-connection is trivial (as a matrix-valued 1-form). A general Spin$^c$-connection with respect to the same trivialization is

$$A = \begin{pmatrix} A^E & 0 \\ 0 & A^E \end{pmatrix}. \tag{3.2}$$

Working in the temporal gauge and letting $z := x_1 + ix_2$, (2.2) reduces to:

$$\begin{cases} 2F^E_{12} = -i/2(t - |\alpha|^2 + |\beta|^2), \\ 2\partial_3 A^E_{\bar{z}} = -\alpha\bar{\beta}, \end{cases} \tag{3.3}$$

where $A^E_{\bar{z}} = A^E_1 + iA^E_2$, and

$$\begin{cases} \partial_3 \alpha + 2\partial_{A^E}\beta = 0, \\ -\partial_3 \beta + 2\bar{\partial}_{A^E}\alpha = 0. \end{cases} \tag{3.4}$$

For $t > 0$, (3.3), (3.4) can be reduced to the $t = 1$ case by rescaling:

$$\delta_t(x) := t^{-1/2}x; \qquad \psi = t^{1/2}(\alpha_1, \beta_1). \tag{3.5}$$

We will hence concentrate on the $t = 1$ case for the rest of this subsection.

There is an obvious family of solutions: we simply take $\beta \equiv 0$; the Seiberg-Witten equations then require $\alpha$, $A^E$ to be independent of $x_3$, and their dependence on $z$ is described by the equations

$$\begin{cases} \bar{\partial}_{A^E}\alpha = 0, \\ 2F^E_{12} = -\frac{i}{2}(t - |\alpha|^2). \end{cases} \tag{3.6}$$

When $t = 1$ this is exactly the vortex equations on $\mathbb{C}$, whose solutions are described in e.g. [40] and [15]. We refer the reader to the Appendix for a list of some of their important properties which are frequently used in this paper. The solutions for the case $t > 0$, $t \neq 1$ will be called *"t-rescaled vortex solutions"*. The prefix *"t-rescaled"* will be often omitted when it is clear from the context.

---
[1] We adopt the convention in quantum mechanics, where $\gamma_2$ differs by a sign with that in many mathematical literature.

Note that since the vortex solutions $(A, \alpha)$ have nonzero $L^2(\mathbb{C})$-norms for $F_A$ and $\nabla_{A^E}\alpha$, the family of Seiberg-Witten solutions described above have infinite $L^2(\mathbb{R}^3)$-norms for the curvature and $\nabla_A \psi$ unless the vortex number is zero.

**3.2.1 Proposition** *The above solutions are the only admissible solutions of (3.3), (3.4) for $t > 0$ in temporal gauge.*

*Proof.* It suffices to show that if $(A, \psi)$ be an admissible solution, then $\beta \equiv 0$, because (3.4) then implies that $\partial_3 \alpha = \bar{\partial}_{A^E}\alpha \equiv 0$, and hence the conclusion of the Proposition.

The Seiberg-Witten equations imply

$$\left(-2\nabla_{A^E}^2 + \frac{1}{2}(1 + |\alpha|^2 + |\beta|^2)\right)\beta = 0. \tag{3.7}$$

Let $\chi_R$ be the cutoff function in Definition 2.2.4. Taking $L^2$-product of (3.7) with $\chi_R \beta$, we have

$$\begin{aligned}
\left\langle \chi_R \beta, (-2\nabla_{A^E}^2 + \frac{1}{2}(1 + |\alpha|^2 + |\beta|^2))\beta \right\rangle &= 0 \\
&= \left\langle 2\nabla_{A^E}(\beta \chi_R), \nabla_{A^E}\beta \right\rangle + \left\langle \chi_R \beta, \frac{1}{2}(1 + |\alpha|^2 + |\beta|^2)\beta \right\rangle \\
&= \left\langle 2\beta \nabla \chi_R, \nabla_{A^E}\beta \right\rangle + \left\langle 2\chi_R \nabla_{A^E}\beta, \nabla_{A^E}\beta \right\rangle \\
&\quad + \left\langle \chi_R \beta, \frac{1}{2}(1 + |\alpha|^2 + |\beta|^2)\beta \right\rangle.
\end{aligned} \tag{3.8}$$

Now

$$|\langle \beta \nabla \chi_R, \nabla_{A^E}\beta \rangle| \leq \|\beta\|_6 \|\nabla \chi_R\|_3 \|\nabla_{A^E}\beta\|_{2,R} \to 0 \text{ as } R \to \infty,$$

since by scaling $\|\nabla \chi_R\|_3 = \|\nabla \chi_1\|_3$. (Cf. [15] lemma VI.3.3.) Here $\|\cdot\|_{2,R}$ denotes $L^2$-norm over the space outside the sphere of radius $R$, as $\nabla \chi_R$ is supported on the annulus of $R \leq |x| \leq 2R$. If $\beta \in V_{1/A^E}$, by lemma 2.2.12 $\|\beta\|_6 \leq \|\nabla_{A^E}\beta\|_2 < \infty$. As both $\langle \chi_R \nabla_{A^E}\beta, \nabla_{A^E}\beta \rangle$ and $\langle \chi_R \beta, \frac{1}{2}(1 + |\alpha|^2 + |\beta|^2)\beta \rangle$ are positive, taking the limit $R \to \infty$ in (3.8), we obtain $\beta \equiv 0$. $\square$

Thus there is a 1-1 map, $j$, from the space of vortex solutions on $\mathbb{C}$ to the space of admissible Seiberg-Witten solutions on $\mathbb{R}^3$,

$$j((A, \alpha)) := \delta_t^*(2p^* A, (t^{1/2} p^* \alpha, 0)), \tag{3.9}$$

where $p$ is the projection from $\mathbb{R}^3$ to $\mathbb{C}$: $x \mapsto z$ (cf. Notation 2.2.5). We call therefore these Seiberg-Witten solutions the *"pull-backs of vortex solutions"*.

## 3.3 Asymptotics of the Solutions

*Throughout this subsection we restrict our attention to $M\backslash M_\Re \sim \mathbb{R}^3\backslash B(\Re)$.* The goal of this subsection is to establish a pointwise asymptotic estimate (Proposition 3.3.3 for $t > 0$; Proposition 3.3.11 for $t = 0$) which tells us that a Seiberg-Witten solution approximates a "reference configuration" (§4.1). To obtain better estimates, we shall use an non-Euclidean coordinate system.

**3.3.1 An atlas on $M\backslash M_{2\Re}$.** We will describe an atlas on $M\backslash M_{2\Re} \sim \mathbb{R}^3\backslash B^3(2\Re)$, induced by the gradient flow of $f$. Let $U_\pm := \{(z, x_3) \in \mathbb{R}^3 : \pm x_3 > \Re \text{ or } |z| > \Re\} \subset M\backslash M_\Re$. $M\backslash M_{2\Re} \subset U_+ \cup U_-$; we will specify coordinate systems on $U_+, U_-$ respectively.

Since $df$ is never zero on this region, $f$ is a good coordinate function. The level surfaces of $f$ form a 2-dimensional foliation of $M\backslash M_\Re$, and the Clifford action of $df$ induces a complex structure on these surfaces. The gradient flow of $f$ induces diffeomorphisms between level surfaces when $f > \Re$ or $f < -\Re$; furthermore, the direct limits of the systems of level surfaces $f \to \pm\infty$ exist and are diffeomorphic to $\mathbb{C}$, which we shall denote as $P_{\pm\infty}$. (This follows directly from the asymptotic condition on $f$ (2.7), which implies, for example, that if the distance $\delta$ between two points on the $f = \Lambda > \Re$ level surface is small, then the distance between their images under the gradient flow on the $P_{+\infty}$ is between $\delta e^{C\Lambda^{-3}}$ and $\delta e^{-C\Lambda^{-3}}$.) Let $\partial_\pm : U_\pm \to P_{\pm\infty}$ be the differentiable maps induced by forward and backward gradient flows respectively, and let $z_\pm$ be the complex coordinates on $P_{\pm\infty} = \mathbb{C}$. Using the same notation for its pull-back via $\partial_\pm$, $z_\pm$ together with $f$ form a coordinate system on $U_\pm$. Of course, in this coordinate system the metric $g_{ij} \neq \delta_{ij}$.

Note that on $M\backslash M_\Re \supset U_+ \cup U_-$, the complex structure on $K^{-1}$ (whose fibers are tangent spaces to the level surfaces) is given by the Clifford action by $df$ described in §2.1. We let $\partial$ denote the corresponding complex differentiation on the level surfaces, which is *not* equal to $\frac{\partial}{\partial z_\pm}$. However in the $f \to \pm\infty$ limits, they coincide because of the asymptotic conditions of the metric $g$ and the function $f$.

**3.3.2 $E$ and $T^*M$ as pull-backs over $U_\pm$.** A trivialization of $E\big|_{U_\pm}$ induces an isomorphism $E\big|_{U_\pm} \to \partial_\pm^* \mathbb{C}\big|_{\mathbb{C}}$ (the latter denotes the pullback of the Higgs bundle over $\mathbb{C} = P_{\pm\infty}$).

On the other hand, over $M\backslash M_\Re$ we will often write:

$$T^*M = K^{-1} \oplus \mathbb{R},$$

where $\mathbb{R}$ denotes the trivial $\mathbb{R}$ bundle spanned by $df$. On $U_\pm$, the connections induced from the Levi-Civita connection identify both components on the right hand side as pull-backs: $K^{-1} \simeq \partial_\pm^* T\mathbb{C}$ (Regarding $K^{-1}$ as a sub-bundle of the tangent bundle); similarly $\mathbb{R} \simeq \partial_{\pm}^* N$, where $N$ is the normal bundle of $\mathbb{C} \times \{0\} \subset \mathbb{C} \times \mathbb{R}$.

**3.3.3 The decay estimate for $t > 0$.** Let $t > 0$ from now on until 3.3.11. Let $\lambda_1 : \mathbb{R} \to [0, 1]$ be a smooth cutoff function, which is 1 on $[1, \infty)$ and vanishes on $(-\infty, -1]$. Let $\lambda_d(x) := \lambda_1(x/d)$.

Given two $t$-rescaled vortex solutions $a, b \in \Gamma(T^*\mathbb{C} \oplus \mathbb{C})$ and a real number $0 < d \leq \Re/2$, with respect to a trivialization of $E$ on $M \backslash M_{2\Re}$ we define the following pair of connection and section on $E\big|_{M \backslash M_{2\Re}}$:

$$v(a, b) := \lambda_d(f)a(z_+) + (1 - \lambda_d(f))b(z_-) \in \Gamma(K^{-1} \oplus E). \qquad (3.10)$$

This makes sense according to 3.3.2, because $\text{support}(\lambda_d \circ f) \cap M \backslash M_{2\Re} \subset U_+$ and $\text{support}((1-\lambda_d) \circ f) \cap M \backslash M_{2\Re} \subset U_-$ on which we adopt respectively the coordinate systems on $U_+, U_-$ described above.

**Proposition** *Let $M$ be an MEE with a $k$-admissible metric, $k \geq l + 5$, $k, l \in \mathbb{Z}^+$, and let $t > 0$, $w = 0$. Let $(A, \psi) = (A, (\alpha, \beta))$ be an admissible Seiberg-Witten solution on $M$.*

*Then for a large enough $R \geq \Re$ and $d$, $R/2 > d \gg 0$ (which may be chosen independently of $(A, \psi)$), there exist $t$-rescaled vortex solutions $v_+, v_- \in C^\infty(T^*\mathbb{C} \oplus \mathbb{C})$ such that $(A_0^E, \alpha_0) := v(v_+, v_-)$, regarded as a pair of connection and section on $E$ via a trivialization of $E\big|_{M \backslash M_R}$ with respect to which $A^E(\nabla f) = 0$, satisfies*

$$\sum_{i=0}^{l} |\nabla_{A^{E'}}^i \beta| + |x|^{-1}\Big(\sum_{i=0}^{l} |\nabla^i(A^E - A_0^E)| + \sum_{i=0}^{l} |\nabla_{A^E}^i(\alpha - \alpha_0)|\Big) \leq C|x|^{-4}, \quad (3.11)$$

*on $M \backslash M_R$, where $C$ is a positive constant dependent on $R, d$, but not on the configuration $(A, \psi)$.*

It will be clear from the proof that though $v_+, v_-$ depend on the choice of trivialization above, the connection and section $(A_0^E, \alpha_0)$ does not. Also, obviously we could have chosen the number $\Re$ in Definition 2.2.2 to be as large as the $R$ in this Proposition. To simplify notation, *we shall let $\Re = R$ in later sections.*

*Proof.* The proof will occupy the rest of this section and is divided into 3 steps. The first step starts with Definition 3.1.1 condition 3 and introduces pointwise estimates for $\beta$ and its derivatives. At Step 2 we deduce the existence of the limits of $\big(A(z_\pm, f), \psi(z_\pm, f)\big)$ as $f \to \pm\infty$. At Step 3 we deduce the estimates for quantities involving $A^E - A_0^E$ and $\alpha - \alpha_0$.

STEP 1. ESTIMATING $\beta$. We will repeatedly make use of the following version of the maximum principle, whose proof is very similar to the proof of Proposition 3.2.1 earlier.

**3.3.4 Lemma** ([15] Proposition VI.3.2, weak maximum principle.) *Let $M$ be as in the previous Proposition and let $V := M \backslash M_R \sim \mathbb{R}^3 \backslash B(R) \subset \mathbb{R}^3$, where $R \geq \Re$. Suppose that*

1. *$u$ is a $L_1^2$ function on $V$ such that $\nabla u \in L^2(V)$, and $u < 0$ on $\partial V$.*

2. *For all compactly supported function $\xi$, $0 \leq \xi \in L_1^2(V)$, the following holds:*

$$\int_V \left( (\nabla \xi) \cdot (\nabla u) + c\xi u \right) \leq 0,$$

*where $c$ is a positive function.*

*Then $u \leq 0$ in $V$.*

The original statement in [15] is in fact stronger: in the above we replaced the original $L^p$ ($p = 6$ in our case) condition for $u$ by $V_1$, which implies the former by Lemma 2.2.12, and is easier to verify in our case.

**3.3.5 Notation** We let $o^n$ or $o_i^n$, $i \in \mathbb{Z}^+$ denote a positive polynomially decaying function $o^n \leq C|x|^{-n}$ for some positive constant $C$, which varies with the context. Because $f = x_3 + o^3$, $z_\pm = z + o^2$, we see that $o^n$ exhibits the same decaying behavior in the new coordinates: $o^n \leq C'(|z_\pm|^2 + |f|^2)^{-n/2}$.

To start the ball rolling, project the equation $\partial_A^2 \psi = 0$ to the $\beta$ component by inner producting with $\beta$ as in [42]:

$$\frac{1}{2} d^* d |\beta|^2 + |\nabla_A \beta|^2 + \frac{1}{2}|\beta|^2 (|\beta|^2 + t|\nabla f| + |\alpha|^2) \leq o_1^4 |\beta|^2 + o_2^4 |\nabla_{A^E} \alpha||\beta|,$$

where the functions $o_1^4$ and $o_2^4$ has the designated decay because of the asymptotic condition of $g$, $\theta$ (cf. Definitions 2.2.2, 2.2.6). Using the triangle inequality, one gets

$$\frac{1}{2} d^* d |\beta|^2 + |\nabla_A \beta|^2 + \frac{1}{2}|\beta|^2(|\beta|^2 + t|\nabla f| + |\alpha|^2) \leq \varepsilon |\beta|^2 + o_3^8 |\nabla_{A^E} \alpha|^2. \quad (3.12)$$

Where $|x| = r > R$ is large enough; $\varepsilon$ is a small positive number depending on $R$. To get rid of the undesired $|\nabla_{A^E} \alpha|^2$ term, perform a similar estimate for $|\alpha|^2$:

$$\frac{1}{2} d^* d |\alpha|^2 + |\nabla_{A^E} \alpha|^2 + \frac{1}{2}|\alpha|^2(|\alpha|^2 - t|\nabla f| + |\beta|^2) \leq \varepsilon' |\alpha|^2 + o_4^8 |\nabla_A \beta|^2. \quad (3.13)$$

We may find a constant $C$ such that by adding $Cr^{-8}$ times (3.13) to (3.12), and using the $L^\infty$-bound for $\psi$, we have with

$$\nu_0 := |\beta|^2 + Cr^{-8}|\alpha|^2, \quad (3.14)$$

$$\frac{1}{2}d^*d\nu_0 + \frac{1}{2}(|\alpha|^2 + |\beta|^2 + t|\nabla f|)\nu_0 \leq \varepsilon_1 \nu_0 + o^8.$$

Take a comparison function
$$h := C'(t) r^{-8}, \qquad (3.15)$$

and let $u := \nu_0 - h$. We choose the constant $C'(t)$ in (3.15) such that $u < 0$ where $r = R$ and $\frac{1}{2}d^*du + C_1 u < 0$ for another constant $C_1$. To apply Lemma 3.3.4 for $u$, we need only verify $\nabla u \in L^2$. This holds because $h$ is evidently in $V_1$, and $\nu_0 \in V_1$ by Lemma 3.3.6 below. Lemma 3.3.4 thus implies that $\nu_0 \leq o^8$ and therefore
$$|\beta| \leq o_2^4. \qquad (3.16)$$

Next we show by induction that
$$|\nabla^i_{A_{E'}}\beta| \leq o^4; \quad |\nabla^i_{A_E}\alpha| \in L^\infty \qquad \text{for } i \leq l+1. \qquad (3.17)$$

Note first that (3.17) holds for $i = 0$ by previous arguments. Assume that (3.17) holds for all $i < k$. We will show that (3.17) holds for $i = k$ as well.

Let $\nabla^k_{A_{E'}}$ acts on the $\beta$ component of the equation $\partial^2_A \psi = 0$. Then inner product with $\nabla^k_{A_{E'}}\beta$. Let $\pi_\beta$ denote the projection to the $\beta$ component. Note that the commutator $[\nabla^k_{A_{E'}}, \pi_\beta \nabla^*_A \nabla_A]$ yields terms involving $F_A, F_{A^K}, \nabla \omega$ and their derivatives. Using the Seiberg-Witten equations, $F_A$ and its derivatives may be substituted by terms involving $\alpha$, $\beta$, and their derivatives; the derivatives of $F_{A^K}, \nabla \omega$ are $o^4$ by the decay conditions in the definition of admissible pairs. Summing up, we thus have

$$\frac{1}{2}d^*d|\nabla^k_{A_{E'}}\beta|^2 + |\nabla^{k+1}_{A_{E'}}\beta|^2 + \frac{1}{2}|\nabla^k_{A_{E'}}\beta|^2(|\beta|^2 + t|\nabla f| + |\alpha|^2)$$
$$\leq (C_3 t + o_5^4)|\nabla^k_{A_{E'}}\beta|^2 + |\nabla^k_{A_{E'}}\beta|(o_6^4|\nabla^k_{A_E}\alpha| + o_7^4|\nabla^{k+1}_{A_E}\alpha| + o^4)$$

by the induction hypothesis. Using the triangle inequality again, the above inequality may be simplified as

$$\frac{1}{2}d^*d|\nabla^k_{A_{E'}}\beta|^2 + |\nabla^{k+1}_{A_{E'}}\beta|^2 + \frac{1}{2}|\nabla^k_{A_{E'}}\beta|^2(|\beta|^2 + t|\nabla f| + |\alpha|^2)$$
$$\leq (C_4 + \varepsilon_2)|\nabla^k_{A_{E'}}\beta|^2 + o_8^8|\nabla^k_{A_E}\alpha|^2 + o_9^8$$
$$+ o_{10}^8|\nabla^{k+1}_{A_{E'}}\alpha|^2 \qquad (3.18)$$

Again to get rid of the unwanted $|\nabla^{k+1}_{A_E}\alpha|^2$ term, we perform a similar estimate for $|\nabla^k_{A_E}\alpha|$:

$$\frac{1}{2}d^*d|\nabla^k_{A_E}\alpha|^2 + |\nabla^{k+1}_{A_E}\alpha|^2 + \frac{1}{2}|\nabla^k_{A_E}\alpha|^2(|\beta|^2 - t|\nabla f| + |\alpha|^2)$$
$$\leq (C_5 t + \varepsilon_3)|\nabla^k_{A_E}\alpha|^2 + o_{12}^8|\nabla^k_{A_E}\beta|^2 + o_{13}^8$$
$$+ o_{14}^8|\nabla^{k+1}_{A_{E'}}\beta|^2. \qquad (3.19)$$

Again we may choose constants $C$, $C'$ and $C_6$ such that by the induction hypothesis, letting

$$\nu_k := |\nabla^k_{A^{E'}}\beta|^2 + Cr^{-8}|\nabla^k_{A^E}\alpha|^2 + C_6|\nabla^{k-1}_{A^{E'}}\beta|^2 + C'r^{-8}|\nabla^{k-1}_{A^E}\alpha|^2, \qquad (3.20)$$

we obtain by summing up suitable multiples of (3.18), (3.19) and $r^{-8}$ times of their $k-1$ versions:

$$\frac{1}{2}d^*d(\nu_k) + \frac{t}{2}\nu_k(|\beta|^2 + |\nabla f| + |\alpha|^2) \leq \varepsilon_4 \nu_k + o_{20}^8.$$

Thus by the comparison principle, using a comparison function $h'$ of the same form as (3.15), we may apply Lemma 3.3.4 to conclude that $\nu_k \leq o^8$, *provided* that $\nu_k \in V_1$. This is verified in Lemma 3.3.6. Thus we obtain

$$|\nabla^k_{A^E}\beta| \leq o^4; \qquad (3.21)$$
$$|\nabla^k_A \alpha| \leq C_7 \qquad (3.22)$$

for large enough $r$, verifying (3.17) for $i = k$, and by induction (3.17) is proved up to Lemma 3.3.6.

**3.3.6 Lemma** *Let $\nu_0$ and $\nu_k$ be as in (3.14) and (3.20) respectively. Then for $0 \leq k \leq l+1$, $\nu_k \in V_1(M \backslash M_R)$.*

*Proof.* We first show that $\nu_0 \in V_1$. By Definition 3.1.1, condition 3, we need only to show that

$$|\nabla_{A^E}\alpha|/r^8 \in L^2(M \backslash M_R). \qquad (3.23)$$

To show this, let $\chi_\Lambda$ be the family of cutoff functions given in Notation 2.2.4 (parameterized by $\Lambda$, $\Lambda > R$), and let $\lambda_\Lambda := (1 - \chi_R)\chi_\Lambda$. Then for $n \in \mathbb{Z}^+$:

$$0 = \left\langle \frac{\partial_A(\lambda_\Lambda \psi)}{r^n}, \frac{\partial_A \psi}{r^n} \right\rangle$$
$$= \left\langle \lambda_\Lambda \frac{\nabla_A \psi}{r^n}, \frac{\nabla_A \psi}{r^n} \right\rangle + \left\langle \lambda_\Lambda \frac{\psi}{r^n}, s\frac{\psi}{r^n} \right\rangle + \frac{1}{2}\sum_{i,j}\left\langle \lambda_\Lambda \psi, \nabla_j\left(\frac{1}{r^{2n}}\right)[\gamma_j, \gamma_i]\nabla_{Ai}\psi \right\rangle$$
$$+ \left\langle \lambda_\Lambda \frac{\psi}{r^n}, \rho(F_A)\frac{\psi}{r^n} \right\rangle + \left\langle \nabla \lambda_\Lambda \frac{\psi}{r^n}, \frac{\nabla_A \psi}{r^n} \right\rangle. \qquad (3.24)$$

Now using the $L^\infty$-bounds of $s$ (scalar curvature) and $F_A$, and applying the triangle inequality, we have

$$\left\| \lambda_\Lambda^{1/2} \frac{\nabla_A \psi}{r^n} \right\|_2^2 \leq (C + C_2/\Lambda)\left\| \lambda_\Lambda^{1/2}\frac{\psi}{r^n} \right\|_2^2 + \varepsilon\left\| \lambda_\Lambda^{1/2}\frac{\nabla_A \psi}{r^n} \right\|_2^2 + C_3(\psi),$$

where $C_3(\psi)$ is a positive constant dependent on $\psi, \varepsilon$, but not on $\Lambda$, that comes from integrals over compactly supported integrands involving $\psi$. Now taking $\Lambda \to$ $\infty$, this means that $(1-\chi_R)^{1/2}\frac{\nabla_A\psi}{r^n} \in L^2$ when $n \geq 2$. Combined with the fact that $\psi \in L^2_{1,loc}$, this implies that $\frac{\nabla_A\psi}{r^n} \in L^2(M\backslash M_R)$ for $n \geq 2$, which implies (3.23).

In fact, arguing in a similar fashion with the equality ($j \in \mathbb{Z}^+$)

$$0 = \Big\langle \frac{\nabla_A^{j-1}\partial_A(\lambda_\Lambda\psi)}{r^n}, \frac{\nabla_A^{j-1}\partial_A\psi}{r^n}\Big\rangle,$$

we see by induction that

$$\frac{\nabla_A^j\psi}{r^n} \in L^2(M\backslash M_R), \quad \text{for } 1 \leq j \leq l+3. \tag{3.25}$$

Next, to show that $\nu_k \in V_1(M\backslash M_R)$ when $k > 0$, it suffices to show that

$$\text{(i) } |\nabla_{A^{E'}}^k\beta|^2 \in V_1(M\backslash M_R) \text{ and} \tag{3.26}$$

$$\text{(ii) } r^{-8}|\nabla_{A^E}^k\alpha|^2 \in V_1(M\backslash M_R) \text{ for all } k \in \mathbb{Z}^+, k \leq l+1. \tag{3.27}$$

To show (3.26), we multiply (3.18) by $\lambda_\Lambda|\nabla_{A^{E'}}^k\beta|^2$ and integrate over $M$, then take $\Lambda \to \infty$. We obtain

$$\|(1-\chi_R)^{1/2}d|\nabla_{A^E}^k\beta|^2\|_2^2$$
$$\leq C\Big(\|(1-\chi_R)^{1/4}\nabla_{A^{E'}}^k\beta\|_4^4 + \|(1-\chi_R)^{1/2}r^{-4}\nabla_{A^{E'}}^k\beta\|_2^2\Big)$$
$$+ C'\|(1-\chi_R)^{1/4}\nabla_{A^{E'}}^k\beta\|_4^2\Big(\|(1-\chi_R)^{1/4}r^{-4}\nabla_{A^E}^k\alpha\|_4^2$$
$$+ \|(1-\chi_R)^{1/4}r^{-4}\nabla_{A^E}^{k+1}\alpha\|_4^2\Big) + C_1(\beta).$$

By Lemma 2.2.12, and (3.25), this is finite if

$$\nabla_{A^{E'}}^j\beta \in L^2(M\backslash M_R) \quad \text{for all } j \in \mathbb{Z}^+, 1 \leq j \leq l+2. \tag{3.28}$$

To verify this, we multiply (3.18) by $\lambda_\Lambda$ and integrate over $M$. Then take $\Lambda \to \infty$. This gives

$$\|(1-\chi_R)^{1/2}\nabla_{A^{E'}}^{k+1}\beta\|_2^2 \leq C\Big(\|(1-\chi_R)^{1/2}\nabla_{A^{E'}}^k\beta\|_2^2 + C'$$
$$+ \|(1-\chi_R)^{1/2}r^{-4}\nabla_{A^E}^k\alpha\|_2^2 + \|(1-\chi_R)^{1/2}r^{-4}\nabla_{A^E}^{k+1}\alpha\|_2^2\Big), \tag{3.29}$$

By (3.25), the right hand side is finite except perhaps for the first term. Note that the $j=1$ case of (3.28) holds by Definition 3.1.1. By induction, (3.29) implies that (3.28) holds in general. Thus (3.26) is proved.

To show (3.27), we multiply (3.19) by $\lambda_\Lambda r^{-16}|\nabla_{A^E}^k \alpha|^2$ and integrate over $M$. Then take $\Lambda \to \infty$. This gives

$$\|(1-\chi_R)^{1/2} r^{-8} d|\nabla_{A^E}^k \alpha|^2\|_2^2$$
$$\leq C_1 \|(1-\chi_R)^{1/6} r^{-4} \nabla_{A^E}^k \alpha\|_6^3 \|(1-\chi_R)^{1/2} r^{-5} \nabla_{A^E}^{k+1}\alpha\|_2$$
$$+ C_2 \|(1-\chi_R)^{1/2} r^{-12} \nabla_A^k \psi\|_2^2 + C_3(\psi)$$
$$+ C' \|(1-\chi_R)^{1/4} r^{-4} \nabla_A^k \psi\|_4^2 \Big( \sum_{i=[k/3]}^{k} \|(1-\chi_R)^{1/4} r^{-4} \nabla_A^i \psi\|_4^4$$
$$+ \|(1-\chi_R)^{1/4} r^{-4} \nabla_A^{k+1} \psi\|_4^2 \Big), \tag{3.30}$$

which is finite by Lemma 2.2.12 and (3.25). This implies (3.27). Lemma 3.3.6 is proved. □

STEP 2. APPEARANCE OF VORTEX SOLUTIONS AND VORTEX NUMBERS.

**3.3.7 Lemma** *With respect to the trivialization of $E$ specified in Proposition 3.3.3, $(A^E(z_\pm, f), \alpha(z_\pm, f))$ converge pointwise to $(A_\pm(z_\pm), \alpha_\pm(z_\pm))$ uniformly (in $z_\pm$) as $f \to \pm\infty$. Moreover, $(A_\pm(z_\pm), \alpha_\pm(z_\pm))$ are both t-rescaled vortex solutions.*

*Proof.* With respect to this trivialization, $A^E(\nabla f) = 0$, so $A^E = (A^E)^{(1,0)} + (A^E)^{(0,1)}$ is fully encoded in the anti-holomorphic 1-form $(A^E)^{(0,1)} \in K^{-1}$. We shall often drop the superscript $(0,1)$ when it is clear from the context. The Seiberg-Witten equations implies

$$|\nabla f|\partial_f \alpha = -2\partial_{A^E}\beta + o^4; \tag{3.31}$$
$$2|\nabla f|\partial_f (A^E)^{(0,1)} = -\bar\alpha \beta. \tag{3.32}$$

by the asymptotic condition of $\theta$ and the $L^\infty$-boundedness of $\psi$. Note that had we adopted the Cartesian coordinates, the right hand side of last equation would have involved an additional $o^3$-term, leading to less desirable estimates.

In Step 1 we established the $o^4$ decay of $|\beta|$ and $|\partial_{A^{E'}}\beta|$; furthermore $\Big||\nabla f| - 1\Big| \leq o^3$, therefore both $\partial_f \alpha$, $\partial_f A^E$ decay as $o^4$. By integrating along the gradient flow lines, we see that the equations imply that both $A^E(z_\pm, f)$ and $\alpha(z_\pm, f)$ converge uniformly to $A_\pm(z_\pm)$ and $\alpha_\pm(z_\pm)$ respectively as $f \to \pm\infty$. On the other hand, the rest of the Seiberg-Witten equations

$$2\bar\partial_{A^E}\alpha = |\nabla f|\partial_f \beta + o^4;$$
$$2F^E_{*df} + i\frac{1}{2}(t - |\alpha|^2) = -i\frac{1}{2}|\beta|^2 + o^3$$

($F^E_{*df}$ denoting the $*df$ component of $F^E$) tell us that

$$2\bar{\partial}_{A_\pm}\alpha_\pm = 0;$$
$$2(F_{A_\pm})_{12} + i\frac{1}{2}(t - |\alpha_\pm|^2) = 0. \tag{3.33}$$

Namely, the limits solve the $t$-rescaled vortex equations (3.6) on $\mathbb{C}$. □

**3.3.8 Digression on Vortex Numbers.** The convergence in the last lemma allows us to introduce the notion of "vortex number" of a Seiberg-Witten solution:

**Definition** The *vortex number* of a Seiberg-Witten solution $(A, \psi)$ is the vortex number (cf. Appendix pt. 1) of the vortex solutions $(A_+, \alpha_+)$ or $(A_-, \alpha_-)$ from the previous lemma.

It is justified below that the vortex numbers of $(A_+, \alpha_+)$ and $(A_-, \alpha_-)$ are the same, so we have an unambiguous definition of vortex numbers for Seiberg-Witten solutions.

The vortex numbers have an alternative description as Chern numbers. First observe the simple fact:

**Lemma** Let $(A_+, \alpha_+)$ be as the above. Then (i) the vortex number of $(A_+, \alpha_+)$ is finite, and (ii) $\sup_{|z_+|=d}|F_{A^E}(z_+, f)| \to 0$ uniformly in $f$ as $d \to \infty$.

The $U_-$ versions of the above statements are true by the same reason.

*Proof.* By differentiating (3.32) with $\partial$ and by the estimates (3.16), (3.21), (3.22) and the $L^\infty$-bound on $\psi$, we have

$$|\partial_f|F^E_{*df}|| \leq |\partial_f F^E_{*df}| \leq o^4. \tag{3.34}$$

Claim (ii) of the lemma follows from integrating the above over gradient flow lines of $f$, and using the exponential decay (with $|z_+|$) of curvature of the vortex solution $(A_+, \alpha_+)$ (Appendix pt. 3) when claim (i) is true, plus (3.32).

To prove claim (i), by Appendix pt. 1 we only have to show that

$$\int_\mathbb{C}(|F_{A_+}|^2 + |\nabla_{A_+}\alpha_+|^2)dz_+d\bar{z}_+$$
$$\leq \iiint_{\{(z_+,f):x_3(z_+,f)\geq\Lambda>>\Re\}}\partial_f(|F^E_{*df}|^2 + |\nabla_{A^E}\alpha|^2)dz_+d\bar{z}_+df$$
$$+(1+C\Lambda^{-3})\int_\mathbb{C}(|F^E|^2 + |\nabla_{A^E}\alpha_+|^2)dzd\bar{z}$$

is finite. But this is true because on the right hand side, the integrand of the first term is $o^4$ by (3.34) and a similar estimate for $\partial_f|\nabla_{A^E}\alpha|$, and the second term is finite by Definition 3.1.1 condition 3. □

Let $H$ be a surface in $U_+$ such that $\partial_+\big|_H$ is surjective and $f(H)$ is bounded. Then the above lemma implies: (i) $\int_H \frac{iF_{A^E}}{2\pi}$ is the relative Chern number of $E\big|_H$, and thus an integer; (ii) by the Bianchi identity this integral is the same for any two such surfaces; (iii) therefore all such integrals equal the $f \to \infty$ vortex number of the Seiberg-Witten solution. Obviously, the analogous statements for the $U_-$ region and the $f \to -\infty$ limit hold also.

Now If one deforms the surface $H$ out of $U_+$, the above Chern number changes only when the surface sweeps across points in $M_\Re$ where $\theta$ vanishes (at which the splitting of the spin bundle (2.1) fails). However, *the $f \to \pm\infty$ vortex numbers are the same*, because $F_A$ is well-defined on the whole $M$, and by the Bianchi identity, last lemma, (2.1), and the asymptotic condition on $\theta$,

$$\lim_{L \to \infty} \frac{i}{2\pi} \int_{\mathbf{P}_{\pm L}} F_A = \lim_{L \to \infty} \left( \frac{2i}{2\pi} \int_{\mathbf{P}_{\pm L}} F_{A^E} + \frac{i}{2\pi} \int_{\mathbf{P}_{\pm L}} F_{A^K} \right)$$
$$= 2 \times (\text{the } f \to \pm\infty \text{ vortex number}),$$

where $\mathbf{P}_L$ is the level surface $f = L$ in $M$. End of the digression.

<u>Step 3</u>. Estimating $a$ and $\eta$. Let $v_+ := (A_+, \alpha_+)$ and $v_- := (A_-, \alpha_-)$ and define $(A_0^E, \alpha_0) := v(v_+, v_-)$ as in the statement of the proposition. We shall now derive the estimates for $a := A^E - A_0^E$ and $\eta := \alpha - \alpha_0$ in the proposition. In fact, below we will only focus on estimates for $a$, since those for $\eta$ are entirely parallel.

Note that $\partial_f A_0^E$ is supported in the region where $|z| \geq \Re$, which is in $U_+ \cap U_-$. In this region $z_-$ can be expressed as a function of $z_+$, and vice versa. By (3.32)

$$\begin{aligned} |\partial_f A_0^E(z_+, f)| &\leq |\nabla \lambda_d(f)| \, |A_+ - A_- \circ z_-|(z_+) \\ &\leq |\nabla \lambda_d(f)| \int_{\mathbb{R}} |\partial_f A^E(z_+, f)| df \leq C|z|^{-3} |\nabla \lambda_d(f)|; \end{aligned}$$

also from (3.32):

$$|\partial_f a| \leq |-\bar{\alpha}\beta| + |\partial_f A_0^E| \leq o^4 + C|z|^{-3} |\nabla \lambda_d(f)|. \tag{3.35}$$

From the construction of $(A_0^E, \alpha_0)$, we know $a$ and $\eta$ approach zero as $|f| \to \infty$; we may integrate over the gradient flow lines using the fact that $\partial_f |a| \leq |\partial_f a|$ to obtain $|a| \leq o^3$.

To obtain pointwise estimates for the derivatives of $a$, we first differentiate (3.32), and argue similarly as above using the estimates for $|\alpha|$, $|\beta|$, and their derivatives obtained in step 1.

This completes the proof of the proposition. □

The following corollary has important applications in sections 4 and 5.

**3.3.9 Corollary** *In the notation in Proposition 3.3.3:*
*(1) $(A^E - A_0^E, (\alpha - \alpha_0, \beta))$ is bounded in the weighted Sobolev norm $L^2_{l;\epsilon}$ for $0 \leq \epsilon < 3/2$.*
*(2) Let $\Theta_c$ be the operator defined by (A.1) (the deformation operator of the vortex equations at c). Let $\pi_{v_\pm}$ be the projection from $L^2(T^*\mathbb{C} \oplus \mathbb{C})$ to $\operatorname{Ker} \Theta_{v_\pm} \subset L^2(T^*\mathbb{C} \oplus \mathbb{C})$. Then regarding $(A^E - A_0^E, \alpha - \alpha_0)$ as a family of functions in $z_\pm$ parametrized by $f$, $f \in \mathbb{R}, |f| > R$, $\pi_{v_\pm}(A^E - A_0^E, \alpha - \alpha_0)$ gives a $\mathbb{C}^n$-valued function over $f$ since $\operatorname{Ker} \Theta_{v_\pm} \simeq \mathbb{C}^n$ (cf. Appendix, point 6). Proposition 3.3.3 implies that $|\pi_{v_\pm}(A^E - A_0^E, \alpha - \alpha_0)| \leq C|f|^{-3}$.*
*(3) $\sum_{i=0}^{l} |\nabla^i v_+ - \nabla^i v_-| \leq C_\pm |z_\pm|^{-3}$, where $C_\pm$ are positive constants independent of $(A, \psi)$.*

*Proof.* The proof for (1) is immediate from the proposition. For (2), using the exponential decay property of elements in $\operatorname{Ker} \Theta_c$ (Appendix, point 6) we may bound $|\pi_{v_\pm}(A^E - A_0^E, \alpha - \alpha_0)|$ by:

$$C\left|\int |x|^{-3} e^{-\gamma|z_\pm|} dz_\pm d\bar{z}_\pm\right| \leq \frac{C_1}{|f|^3} \left|\int e^{-\gamma|z_\pm|} dz_\pm d\bar{z}_\pm\right| \leq C_\pm |f|^{-3}.$$

(3) is a by-product of Step 3 of the proof for the proposition. □

We remark that in section 4, we work in another gauge $d^*(A - A_0) + i \operatorname{Im} \psi_0 \cdot (\psi - \psi_0) = 0$. ($\psi_0 = (\alpha_0, 0)$; $A_0$ is the connection on $L$ associated to $A_0^E$.) By Proposition 3.3.3, that gauge is asymptotic to the present gauge in the sense that in the present gauge, $|d^*(A - A_0) + i \operatorname{Im} \psi_0 \cdot (\psi - \psi_0)| \leq o^3$.

**3.3.10 Definition** The "centers" of vortex solutions on $\mathbb{C}$ mean the zero loci of their Higgs fields (which consist of points in this case, cf. Appendix). Similarly, the "centers" of perturbed Seiberg-Witten solutions ($t > 0$) are the zero loci of the spinor fields (which consist of paths that are asymptotic to lines in the $x_3$-direction where $x_3 \to \pm\infty$).

**3.3.11 The $t = 0$ version.** Note that in the proof of Proposition 3.3.3 above, we did not actually use the condition $t > 0$ except at Step 2. When $t = 0$, the proof of Lemma 3.3.7 still works, but now (3.33) implies that $\alpha_\pm \equiv 0$; $A_\pm(z_\pm) = id\xi_\pm(z_\pm)$ are flat connections. We then may define $\alpha_0 \equiv 0$, and

$$A_0^E := id\xi_0, \text{ where } \xi_0 := \lambda_d(f)\xi_+(z_+) + (1 - \lambda_d(f))\xi_-(z_-),$$

interpreted as with (3.10).

**Proposition** *Let $t = 0$. Under the same assumptions as Proposition 3.3.3, an admissible Seiberg-Witten solution $(A, \psi)$ in this case also satisfies (3.11), with $(A_0^E, \alpha_0)$ reinterpreted as above.*

# 4 The configuration space and the quotient space

In §4.1, we define the relevant norms for the purpose of proving sliceness of the configuration spaces under the gauge group actions. We construct the configuration spaces and establish the Banach manifold structure of the quotient space in the $t = 0$ case and the $t > 0$ case in §4.2, 4.3 respectively.

## 4.1 Analytical preliminaries

Suppose that $M$ is an MEE with a $k$-admissible pair $(g, \theta)$, $k \geq l + 5$.

**4.1.1 The Reference Configurations.** It turns out that the configuration spaces we consider will be contained in some coordinate patches based on certain *reference configurations*. Before we define them, let us fix some notation. Let

$$\mathcal{A}_l^n \subset L_l^2(\mathbb{C}, T^*\mathbb{C}) \times L_{l,loc}^2(\mathbb{C}, \mathbb{C})$$

be a neighborhood (in the $L_l^2$ sense) of $\mathcal{M}_{vortex}^n$–the embedding of the moduli space of vortex solutions with vortex number $n$ described in the Appendix. We shall take the neighborhood to be small enough for our purposes. In particular, $\Theta_c$ (cf. (A.1)) is invertible when $c \in \mathcal{A}_l^n$. An element in $\mathcal{A}_l^n$ is called an *approximate vortex solution*.

**Definition** (1) In the $t = 0$ case, $c_0 = (A_0, \psi_0)$ is a *reference configuration* if there is a trivialization of $E\big|_{M \setminus M_\Re}$ with respect to which $(A_0^E, \psi_0) = (0, 0)$ on $M \setminus M_\Re$.

(2) In the $t > 0$ case, a *reference configuration* $c_{a_1 a_2} = (A_0, \psi_0)$ associated to two approximate vortex solutions $a_1, a_2 \in \mathcal{A}_l^n$ is a configuration such that over $M \setminus M_\Re$, $(A_0^E, \psi_0) = (A_0^E, (\alpha_0, 0))$ where $(A_0^E, \alpha_0) = v(a_1, a_2)$ defined in (3.10).

We next introduce some useful norms and discuss their properties.

**4.1.2 Definition** Let $M$ be an MEE, and let $V$ be a Euclidean/hermitian bundle over $M$ constructed from $S$ and/or $T^*M$. ($V$ may be trivial.) Let $c = (A, \psi)$ be a reference configuration. Let $\nabla_A$ denote the covariant derivative of $V$ derived from the Levi-Civita connection on $T^*M$ and/or $A$. (When $V$ is trivial, $\nabla$ is the ordinary derivative.) Let $\xi \in C_0^\infty(V)$. Define

$$\|\xi\|_c^2 := \|\nabla_A \xi\|_{2:\epsilon}^2 + \|\xi|\psi|\|_{2:\epsilon}^2,$$

$\epsilon$ being zero for the $t = 0$ case and $\epsilon \in [0, 3/2)$ in the $t > 0$ case. Let

$$\|\xi\|_{Z_{l,c}} := \begin{cases} \|\xi\|_{2,l/A} + \|\xi\|_{6/5} & \text{in the } t = 0 \text{ case;} \\ \|\xi\|_{2,l:\epsilon/A} & \text{in the } t > 0 \text{ case.} \end{cases}$$

If $V$ is a trivial $\mathbb{C}$-bundle or $\mathbb{R}$-bundle and $l \in \mathbb{Z}^+, l \geq 2$, define

$$\|\xi\|_{X_{l,c}} := \begin{cases} \|\xi\|_c + \|\nabla\nabla\xi\|_{2,l-2} + \|-\nabla^2\xi + \xi|\psi|^2\|_{6/5} & \text{when } t = 0; \\ \|\nabla\xi\|_c + \|\xi\psi\|_c & \text{when } t > 0, l = 2; \\ \|\nabla\xi\|_c + \|\xi\psi\|_c + \|\nabla\nabla^2\xi\|_{2,l-3:\epsilon} & \text{when } t > 0, l \geq 2. \end{cases}$$

For another configuration $c'$, with $(a, \eta) := c' - c$ and $l \geq 1$, define the following norm in the tangent space to the configuration space

$$\|(a,\eta)\|_{Y_{l,c}} := \begin{cases} \|a\|_c + \|\nabla a\|_{2,l-1} + \|\eta\|_c + \|\nabla_A \eta\|_{2,l-1/A} \\ \quad + \|d^*a + i\operatorname{Im}\psi \cdot \eta\|_{6/5} & \text{in the } t = 0 \text{ case}; \\ \|a\|_c + \|\nabla a\|_{2,l-1:\epsilon} + \|\eta\|_c + \|\nabla_A\eta\|_{2,l-1/A:\epsilon} & \text{in the } t > 0 \text{ case}. \end{cases}$$

Denote by $H_c(V)$, $Z_{l,c}(V)$, $X_{l,c}$, $Y_{l,c}$, the corresponding completions.

**4.1.3 Remark** (1) These norms are obviously defined for more general configuration $c$. Later we will encounter versions of these norms for a general $c$ in the configuration space $\mathcal{C}$ (cf. §4.2, 4.3).

(2) The norms are chosen so that some invertibility results (Proposition 4.1.8, Lemma 4.2.5 below) needed for the existence of the quotient spaces hold. The reason why the norms needed in the $t > 0$ and $t = 0$ cases are so different was briefly discussed in §1.3. (3) When $V$ is not specified, it is implied to be a trivial bundle. It is clear from the definition that in this case $Z_{l,c}$ does not depend on $c$. Furthermore, we shall see from lemmas 4.1.4, 4.2.2, 4.3.3 that these spaces in fact do not depend on the choice of the reference configuration $c$. We will then drop the subscript $c$.

These spaces have many equivalent definitions according to the following lemma:

**4.1.4 Lemma** *Let $V$, $c$, $\xi$ be as in Definition 4.1.2. Then there exist (c-dependent) constants $\mu_c$, $\mu'$, $\nu_c$ such that:*

$$\mu_c \|\xi\|_{2,1:\epsilon/A} \leq \|\xi\|_c \leq \mu' \|\xi\|_{2,1:\epsilon/A}; \quad \text{in the } t > 0 \text{ case}, \tag{4.1}$$

$$\|\xi\|_{2;1/A} \leq \|\xi\|_c \leq \nu_c \|\xi\|_{2;1/A} \quad \text{in the } t = 0 \text{ case}. \tag{4.2}$$

*Consequently, in the $t > 0$ case, $H_c(V)$, $Y_{l,c}$, $Z_{l,c}(V)$ and $X_{l,c}$ are commensurate with $L^2_{l:\epsilon/A}(V)$, $L^2_{l:\epsilon/A}(iT^*M \oplus S)$, $L^2_{l:\epsilon/A}(V)$ and $L^2_{l:\epsilon}$ respectively. In the $t = 0$ case, by Lemma 2.2.12, $H_c$ is commensurate with $V_{1/A}$; $\|\xi\|_{X_{l,c}}$ is bounded above and below by multiples of $\|\xi\|_{2;l/A} + \|\nabla^2\xi\|_{6/5}$; $\|(a,\eta)\|_{Y_{l,c}}$ is bounded above and below by multiples of $\|(a,\eta)\|_{2;l/A} + \|d^*a\|_{6/5}$.*

*Proof.* The inequality $\|\xi\|_c \leq \mu' \|\xi\|_{2,1:\epsilon/A}$ in the $t > 0$ case follows from the $L^\infty$-bound on $\psi$ by the definition of reference configurations, and the inequality $\|\xi\|_c \geq$

$\|\xi\|_{2;1/A}$ in the $t = 0$ case follows immediately from the definition. We still need to show:

$$\|\xi\|_c \geq \mu_c \|\xi\|_{2,1:\epsilon/A} \qquad \text{for the } t > 0 \text{ case, and} \tag{4.3}$$

$$\|\xi\|_c \leq \nu_c \|\xi\|_{2;1/A} \qquad \text{in the } t = 0 \text{ case.} \tag{4.4}$$

(4.4) is proved by lemma 2.2.12, as

$$\|\xi\|_c \leq \|\nabla_A \xi\|_2 + \|\xi|\psi|\|_2 \leq \|\nabla_A \xi\|_2 + \|\psi\|_3 \|\nabla_A \xi\|_2,$$

and $\|\psi\|_3$ is finite because by the definition of the reference configuration, $\psi \in C_0^\infty$.

We use a partition of unity argument for (4.3). Say $c = c_{a_1 a_2}$ for two approximate vortex solutions $a_1, a_2$. Let $\xi = \xi_1 + \xi_2$, where $\xi_1 = \chi_R \xi$; $\xi_2 := (1 - \chi_R)\xi$. Here we choose $R = R(c) > \Re$ to be large enough so that when $|f| \leq d$ ($d$ is as in Proposition 3.3.3), the zero locus of $\psi$ lies in the region where $|x| < R/2$, and that $R$ is much larger than the distance of the centers of $a_1$, $a_2$ from origin.

Since $c$ is a reference configuration, on the support of $\xi_2$ (which is in $\mathbb{R}^3$) $c := ((A^E, \alpha), \beta) := (\lambda_d(f) a_1(z_+) + (1 - \lambda_d(f)) a_2(z_-), 0)$ with respect to the trivialization of $E$ specified earlier. In this case, we separate the variables $x = (z_\pm, f)$ and use the following similar result in 2-dimensions.

**4.1.5 Sublemma** Let $(a, \phi)$ be an approximate vortex solution on $\mathbb{C}$; let $q$ be a $C_0^\infty$ function on $\mathbb{C}$. Then

$$\|\nabla_a q\|_2^2 + \|q\|_2^2 \leq C(\|\nabla_a q\|_2^2 + \|q|\phi|\|_2^2)$$

for some positive $C$ depending only on the vortex number of $(a, \phi)$.

This is the key inequality for proving the invertibility of the operator $\Theta$ in (A.1). When $(a, \phi)$ is an honest vortex solution this is well-known (an analogue is Lemma 4.6 in [43]); according to our definition it is also true for approximate vortex solutions. Roughly speaking, the proof follows from the fact that $|\phi|(z)$ is almost constant where $|z| > R, R \gg 0$ [40, 15], and that the zero locus of $\phi$ is compactly supported. (Use a partition of unity argument similar to (4.7) below.) □

Accordingly,

$$\|\xi_2\|_c^2 = \|\nabla_A \xi_2\|_{2:\epsilon}^2 + \|\xi_2|\psi|\|_{2:\epsilon}^2$$

$$\geq (1 - CR^{-3}) \int df \varsigma^{2\epsilon} (\|\nabla_A \xi_2\|_{2,z_\pm}^2 + \|\xi_2|\psi|\|_{2,z_\pm}^2) \tag{4.5}$$

$$\geq C' \int df \varsigma^{2\epsilon} (\|\nabla_A \xi_2\|_{2,z_\pm}^2 + \|\xi_2\|_{2,z_\pm}^2) \tag{4.6}$$

$$\geq C'' (\|\nabla_A \xi_2\|_{2:\epsilon}^2 + \|\xi_2\|_{2:\epsilon}^2),$$

where $\|\cdot\|_{2,z_\pm}$ means the $L^2$ norm is taken over the 2-dimensional space parameterized by $z_\pm$ (the sign $\pm$ is determined by the sign of $f$). Note that to go from (4.5) to (4.6), we used the above sublemma where $|f| > d$, and we need the condition on $R$ to obtain a lower bound for $|\psi|$ where $|f| \leq d$ on the support of $\xi_2$.

On the other hand, $\xi_1$ also satisfies

$$\|\xi_1\|_c^2 \geq \nu^2(\|\nabla_A\xi_1\|_{2:\epsilon}^2 + \|\xi_1\|_{2:\epsilon}^2)$$

for some positive constant $\nu$ because it is compactly supported: $\|\xi_1\|_2 = \|\chi_{2R}\xi_1\|_2 \leq C\|\chi_{2R}\|_3\|\nabla_A\xi_1\|_2$. Thus

$$\begin{aligned}
\|\nabla_A\xi\|_{2:\epsilon}^2 + \||\psi|\xi\|_{2:\epsilon}^2 &= \|\nabla_A(\xi_1 + \xi_2)\|_{2:\epsilon}^2 + \||\psi|(\xi_1 + \xi_2)\|_{2:\epsilon}^2 \\
&= \|\nabla_A\xi_1\|_{2:\epsilon}^2 + \|\nabla_A\xi_2\|_{2:\epsilon}^2 + \||\psi|\xi_1\|_{2:\epsilon}^2 + \||\psi|\xi_2\|_{2:\epsilon}^2 \\
&\quad + 2\operatorname{Re}\langle|\psi|\xi_1, |\psi|\xi_2\rangle_{2:\epsilon} - 2\|(\nabla\chi_R)\xi\|_{2:\epsilon}^2 \\
&\quad + 2\operatorname{Re}\langle(\nabla\chi_R)\xi, (1-2\chi_R)\nabla\xi\rangle_{2:\epsilon} + 2\langle\chi_R\nabla_A\xi, (1-\chi_R)\nabla_A\xi\rangle_{2:\epsilon} \\
&\geq \frac{C_1}{2}(\|\nabla_A\xi\|_{2:\epsilon}^2 + \|\xi\|_{2:\epsilon}^2) - C_2\|\varsigma^\epsilon\nabla\chi_R\|_3^2\|\nabla\xi\|_2^2 \\
&\quad - C_3\|\varsigma^\epsilon(\nabla\chi_R)(1-2\chi_R)\|_3\|\nabla_A\xi\|_2^2, \\
&\geq C_4(\|\nabla_A\xi\|_{2:\epsilon}^2 + \|\xi\|_{2:\epsilon}^2) - C_5\|\nabla_A\xi\|_{2:\epsilon}^2 \quad (4.7)
\end{aligned}$$

where $C_2, C_3$ are positive constants independent of $c$. Rearranging we get the desired inequality. $\square$

**4.1.6 Remark** In the $t > 0$ case, Lemma 4.1.4 in fact holds for any $c \in \mathcal{C}_l$ (cf. §4.3) by combining the above proof with Lemma 2.2.12.

**4.1.7 Lemma** *In the $t = 0$ case, $\langle \zeta, \eta \rangle \leq C\|\zeta\|_{Z_{l,c}}\|\eta\|_c$ for some $c$-independent constant $C$.*

*Proof.* $\langle \zeta, \eta \rangle \leq \|\eta\|_6\|\zeta\|_{6/5} \leq C\|\eta\|_c\|\zeta\|_{Z_{l,c}}$. $\square$

The norms in Definition 4.1.2 are useful due to the following Lemma.

**4.1.8 Lemma** *Let $V$ in Definition 4.1.2 be a trivial $\mathbb{R}$-bundle or $\mathbb{C}$-bundle. Then for either the $t = 0$ or the $t > 0$ cases, and for any $l \geq 1$, the operator*

$$L_c \equiv -\nabla^\dagger\nabla + |\psi|^2$$

*is an isomorphism between $X_{l+1,c}$ and $Z_{l-1}$, where $\nabla^\dagger$ denotes the formal $L_\epsilon^2$-adjoint of $\nabla$.*

*Proof.* The continuity of $L_c$ is obvious for the $t > 0$ case by lemma 4.1.4. $L_c$ is continuous in the $t = 0$ case by the following estimate:

$$\| -\nabla^2 u + u|\psi|^2 \|_{Z_{l-1}}$$
$$\leq \| -\nabla^2 u + u|\psi|^2 \|_{6/5} + \|\psi\|_\infty \|u\psi\|_2 + C\|u\|_{2;l+1} + \|u\|_{2;l} \|\psi^2\|_{2;l}$$
$$\leq C' \|u\|_{X_{l+1,c}}$$

for constants $C$, $C'$ depending on $c$.

To construct an inverse for $L_c$, for each $\xi \in Z_{l-1}$, we use a standard variational argument (cf. e.g. [6] p.55) using the functional $S_\xi(u) = \frac{\|u\|_c^2}{2} + \langle u, \xi \rangle_{2:\epsilon}$ on $H_c$ to produce a $u_\xi \in L^2_{3,loc}$ minimizing $S_\xi$; $L_c u_\xi = \xi$ almost everywhere. This $u_\xi$ is the candidate for $L_c^{-1}\xi$, but we still need to obtain some bounds on $u_\xi$ to show that $u_\xi \in X_{l+1,c}$.

For this purpose, note that since by the definition of $u_\xi$, $\|u_\xi\|_c^2 = \langle u_\xi, \xi \rangle$,

$$\|u_\xi\|_c \leq C\|\xi\|_{Z_{l-1}} \tag{4.8}$$

by lemma 4.1.7 (the case $t = 0$) and lemma 4.1.4 (the case $t > 0$). To prove that $\|u_\xi\|_{X_{l+1,c}}$ is bounded by $\|\xi\|_{Z_{l-1}}$, in the $t = 0$ case it follows from the following two additional estimates: Firstly,

$$\|L_c u_\xi\|_{6/5} = \|\xi\|_{6/5} \leq \|\xi\|_{Z_{l-1}};$$

secondly, we show by induction that

$$\|\nabla^i u_\xi\|_2 \leq C\|\xi\|_{Z_{l-1}} \tag{4.9}$$

for all $i \in \mathbb{Z}^+$, $i \leq l+1$. For $i = 1$, this is true by (4.8); for $i \geq 2$, by the induction assumption and via integration by parts

$$\|\nabla^{i-2} \nabla \nabla u_\xi\|_2 = \|\nabla^{i-2} \nabla^\dagger \nabla u_\xi\|_2$$
$$= \|\nabla^{i-2}(L_c u_\xi - u_\xi |\psi|^2)\|_2$$
$$\leq \|\xi\|_{2,l-1} + C\|u_\xi\|_{2;i-1} \||\psi|^2\|_{3,i-2}$$
$$\leq C'\|\xi\|_{Z_{l-1}}.$$

Therefore $L_c$ has a bounded inverse.

In the $t > 0$ case, (4.8) and an estimate similar to (4.9) do the job. □

### 4.2 Configurations and gauge group action: the case $t = 0$

In this subsection we concentrate on the case $t = 0$, though some lemmas will be proved for both cases $t = 0$ and $t > 0$.

(A) THE CONFIGURATION SPACE.

We start with some general discussions. Let $\mathcal{W}_0$ be the vector bundle

$$\mathcal{W}_0 := iT^*M \oplus S, \tag{4.10}$$

where $S$ is the spinor bundle. Most generally, the (total) configuration space $\tilde{\mathcal{C}}$ should consists of all admissible configurations in $L^2_{2,loc}(\mathcal{W}_0)$, with the (total) gauge group $\tilde{\mathcal{G}} = L^2_{3,loc}(\mathrm{Aut}(L))$ acting on it by

$$g(c) = (A - (dg)g^{-1}, g^{1/2}\psi), \tag{4.11}$$

where $c = (A, \psi) \in \tilde{\mathcal{C}}$, $A$ is a connection of the auxiliary line bundle $L$ which we may identify with an element of $L^2_{2,loc}(T^*M)$ by subtracting off a fixed connection; $\psi$ is a section of $S$. Given the topology of the local Sobolev spaces on $\tilde{\mathcal{C}}$ and $\tilde{\mathcal{G}}$, we can endow the quotient space $\tilde{\mathcal{Q}} := \tilde{\mathcal{C}}/\tilde{\mathcal{G}}$ with the quotient topology (Lemma 4.2.4 below shows that it is Hausdorff) and that in turn induces a topology on the moduli space by regarding the moduli space as a subset of the quotient space. In these topologies, Lemma 3.3.7 defines two continuous maps $\partial_\pm$ from the ($t > 0$) Seiberg-Witten moduli space to the moduli space of vortices on $\mathbb{C}$ given by the $f \to \pm\infty$ limits of the solution. However, as the local Sobolev spaces that $\tilde{\mathcal{C}}, \tilde{\mathcal{G}}$ model on are quite intractable, this is not sufficient to ensure Banach space structures on the quotient space and the moduli space, and we need to restrict our attention to a better behaved subset. On the other hand, as we are mainly concerned with the moduli space of the Seiberg-Witten equations in $\tilde{\mathcal{Q}}$, it suffices to look at neighborhoods of the space of solutions: $\mathcal{C}_l \subset \tilde{\mathcal{C}}$.

We start with a fixed configuration $c_0$ close to the solution space:

**4.2.1 Definition** Set $t = 0$. Let $c_0 = (A_0, \psi_0)$ be a reference configuration. Define the *configuration space* $\mathcal{C}_l$, $l \in \mathbb{Z}^+$ in this case as:

$$\mathcal{C}_l := \left\{ c = (A, \psi) \in L^2_{l,loc}(\mathcal{W}_0) : \|(a, \eta)\|_{Y_{l,c_0}} < \infty \right\}, \tag{4.12}$$

where $(a, \eta) := (A, \psi) - (A_0, \psi_0)$. Define the topology on $\mathcal{C}_l$ by the $Y_{l,c_0}$-norm.

Note that $\mathcal{C}_l$ defined this way does not depend on the choice of the reference configuration: suppose $c_1$ and $c_2$ are two reference configurations, then $\left| \|c - c_1\|_{Y_{l,c_1}} - \|c - c_1\|_{Y_{l,c_2}} \right| + \|c_1 - c_2\|_{Y_{l,c_2}}$ depends on finite integrals over a compact space; also $Y_{l,c_1}, Y_{l,c_2}$ are commensurate by lemmas 2.2.12 and 4.1.4. In fact,

**4.2.2 Lemma** *For any $k \in \mathbb{Z}^+$, the norms $V_{k/A}$ are commensurate for different reference configurations $c = (A, \psi)$ and any $k \in \mathbb{Z}^+$. The spaces $L^2_{k/A}$ are commensurate for different $c \in \mathcal{C}_l$ if $k \leq l$, and inequalities of the type of (2.8) hold.*

Henceforth we drop the subscript $A$ (cf. Remark 2.2.9).

(B) THE GAUGE GROUP ACTION AND THE QUOTIENT MANIFOLD.

Let $\mathcal{G}_{l+1} \subset \tilde{\mathcal{G}}$ be the stabilizer of $\mathcal{C}_l \subset \tilde{\mathcal{C}}$. It follows immediately from (4.11) that the action of $\mathcal{G}_{l+1}$ on $\mathcal{C}_l$ is free except at configurations of the form $(A, 0)$, which we call *reducible configurations*.

Let $\mathcal{C}_l^*(M) \in \mathcal{C}_l$ denote the set of irreducible configurations; we aim to describe the local structure of the quotient $\mathcal{Q}_l^* := \mathcal{C}_l^*(M)/\mathcal{G}_{l+1} \subset \tilde{\mathcal{C}}/\tilde{\mathcal{G}}$.

We begin by enumerating some simple facts.

The linearization of the gauge group action at the configuration $c$ is $d_c$,

$$d_c \xi := (-\nabla \xi, \frac{1}{2} \xi \psi); \tag{4.13}$$

where $\xi$ is an imaginary function and $\xi$ acts on $\psi$ by complex multiplication. We take the pointwise (real) inner product of $(a, \eta)$ and $(b, \chi)$ to be

$$(a, \eta) \cdot (b, \chi) = a \cdot b + 2 \operatorname{Re} \eta \cdot \chi. \tag{4.14}$$

$d_c$ then has the formal $L^2$-adjoint

$$d_c^*(a, \eta) = -d^* a - i \operatorname{Im} \psi \cdot \eta.$$

As usual, we construct a slice for the action of $\mathcal{G}_{l+1}$ by fixing a gauge. For technical reasons we choose a less conventional gauge $d_c^*(a, \eta) = 0$.

**4.2.3 Theorem** *Let $t = 0$. Then the quotient $\mathcal{Q}_l^*$ above is a Banach manifold covered by local coordinate patches modeled on*

$$Q_{l,c} := \left\{(a, \eta) \in \Gamma(\mathcal{W}_0) : \|a\|_{2;l} + \|\eta\|_{2;l/A} < \infty, \ d_c^*(a, \eta) = 0\right\}, \tag{4.15}$$

*where $c = (A, \psi)$ is an irreducible reference configuration.*

*Proof.* The proof of theorem 4.2.3 is an adaptation of the standard procedures.

First we note the following lemma which shows that $\tilde{\mathcal{C}}/\tilde{\mathcal{G}}$ is Hausdorff. Since $\mathcal{Q}_l^* \subset \tilde{\mathcal{C}}/\tilde{\mathcal{G}}$ and has finer topology than the induced topology from $\tilde{\mathcal{C}}/\tilde{\mathcal{G}}$, this implies that $\mathcal{Q}_l^*$ is Hausdorff as well.

**4.2.4 Lemma** (both cases $t = 0$ and $t > 0$) *The map*

$$\tilde{\mathcal{G}} \times \tilde{\mathcal{C}} \to \tilde{\mathcal{C}} : \ (g, c) \to g(c) \tag{4.16}$$

*with $g(c)$ defined in (4.11) is continuous in the topology of the local Sobolev norms. Moreover, the action of $\tilde{\mathcal{G}}$ on $\tilde{\mathcal{C}}^*$ has a closed graph. That is, let $(c_n) \subset \tilde{\mathcal{C}}^*$ and $(g_n) \subset \tilde{\mathcal{G}}$ be sequences so that $c_n \to c$ and $g_n(c_n) \to c'$ in $\tilde{\mathcal{C}}^*$. Then there is a subsequence $g_n \to g$ in $\tilde{\mathcal{G}}$ with $g \in \tilde{\mathcal{G}}$ satisfying $g(c) = c'$.*

This lemma is proved by a standard argument via Sobolev embedding (cf. e.g. [9] p.50, Appendix A).

Next, to prove sliceness we need an invertibility lemma.

**4.2.5 Lemma** (both cases $t = 0$ and $t > 0$) *Let $c$ be a reference configuration. Then the linear maps $d_c$, $d_c^*$:*

$$X_{l+1,c} \xrightarrow{d_c} Y_{l,c} \xrightarrow{d_c^*} Z_{l-1} \qquad (4.17)$$

*are continuous, and $\operatorname{Ker} d_c^* \subset Y_{l,c}$ is commensurate with the Banach space $Q_{l,c}$, the tangent space of $\mathcal{Q}_l$. Moreover, $d_c^* d_c$ is an isomorphism.*

*Proof.* For the $t > 0$ case the continuity is obvious from lemma 4.1.4. The continuity for the $t = 0$ case follows from the following routine estimates

$$\begin{aligned} &\|d^*a + i\operatorname{Im}\psi \cdot \eta\|_{Z_{l-1}} \\ &\leq \|\nabla a\|_{2,l-1} + \||\psi|\eta\|_{2,l-1} + \|d^*a + i\operatorname{Im}\psi \cdot \eta\|_{6/5} \\ &\leq C\|(a,\eta)\|_{Y_{l,c}}; \end{aligned} \qquad (4.18)$$

$$\begin{aligned} &\|(-2\nabla\xi, \xi\psi)\|_{Y_l} \\ &\leq 2\|\nabla\xi\|_{2,l} + C'\|\xi\|_{2;l+1}\|\psi\|_{2;l+1/A} + 3(\|\xi\|_{2;l} + \|\xi\|_c)\|\psi\|_\infty + \|L_c\xi\|_{6/5} \\ &\leq C\|\xi\|_{X_{l+1,c}}. \end{aligned} \qquad (4.19)$$

For the $t = 0$ case and the $\epsilon = 0$ case of the $t > 0$ case, the isomorphism between $X_{l+1,c}$ and $Z_{l-1}$ is the direct consequence of Lemma 4.1.8, as $d_c^* d_c = -L_c$. When $\epsilon > 0$, the isomorphism still holds because $d_c^* d_c$ differs from $-L_c$ by a small perturbation ($\sim 2\epsilon\varsigma^{-1}\partial_3\varsigma$, cf. Definition 2.2.10). □

**4.2.6 Remark** The claims in this lemma in fact hold for any $c \in \mathcal{C}_l$ when $t > 0$ (cf. §4.3).

We now define the Banach Lie group $\mathcal{G}_{X_{l+1}} \subset \mathcal{G}_{l+1}$ for the $t = 0$ case. Let $X_{l+1} := X_{l+1,c}$ for some reference configuration $c$. (Note that the definition of $X_{l+1}$ does not depend on the choice of $c$, cf. lemma 4.1.4.)

First we note that the exponential map $\exp : i\Omega^0 \to \operatorname{Aut}(L)$ is a bijection on a neighborhood $U_X$ of $0 \in X_{l+1}$.

**4.2.7 Lemma** *The exponential map above extends to the whole $X_{l+1}$. Its image, denoted $\mathcal{G}_{X_{l+1}}$, is a Banach Lie group with Banach Lie algebra $X_{l+1}$, and acts smoothly on $Y_{l,c}$ by $g[y] = g(c+y) - c$.*

*Proof.* The infinitesimal action of $X_{l+1}$ on $Y_{l,c}$ is

$$\phi: \ (\xi, (a, \eta)) \mapsto (-\nabla \xi, \frac{\xi}{2}\psi) + (0, \frac{\xi}{2}\eta). \tag{4.20}$$

The $Y_{l,c}$-norm of the first term is bounded by the $X_{l+1}$-norm of $\xi$ by lemma 4.2.5; the $Y_{l,c}$-norm of the second term can be bounded by

$$\mu \|\nabla_A(\xi\eta)\|_{2,l-1} + \|\xi \operatorname{Re} \psi \cdot \eta\|_{6/5}$$
$$\leq \mu \|d\xi\|_{3,l-1} \|\nabla_A \eta\|_{2;l} + C(\|d\xi\|_{2,l-1} + C'\|d\xi\|_{4,l-1})\|\nabla_A \eta\|_{2,l-1}$$
$$+ \|\xi\|_6 \|\psi\|_2 \|\eta\|_6$$
$$\leq C_2 \|\eta\|_{Y_{l,c}} \|\xi\|_{X_{l+1}},$$

so the infinitesimal action is continuous.

The rest of proof is standard. (Cf. e.g. [23].) □

The implicit function theorem plus Lemma 4.2.5 imply via a standard argument (cf. e.g. [2] chapter 4) that

**4.2.8 Lemma** *When $t = 0$, $\mathcal{C}_l^*/\mathcal{G}_{X_{l+1}}$ is a Banach manifold covered by local coordinate patches modeled on $Q_{l,c}$ given in (4.15).*

This is not quite what we want yet, since our $\mathcal{Q}_l^* := \mathcal{C}_l^*(M)/\mathcal{G}_{l+1}$ is the quotient space by $\mathcal{G}_{l+1}$ (the stabilizer of $\mathcal{C}_l$) instead. However, we observe:

**4.2.9 Lemma** *$\mathcal{G}_{X_{l+1}}$ is the component of $\mathcal{G}_{l+1}$ containing 1.*

*Proof.* As $\mathcal{G}_{l+1}/\mathcal{G}_{X_{l+1}}$ acts freely on $\mathcal{C}_l^*/\mathcal{G}_{X_{l+1}}$, if $\mathcal{G}_{l+1}/\mathcal{G}_{X_{l+1}}$ is not discrete, by the local model for $\mathcal{C}_l^*/\mathcal{G}_{X_{l+1}}$ given in the previous lemma, there exists a nonzero $\xi \in L^2_{3,loc}$ which satisfies $d_c^* d_c \xi = 0$, which, by the unique continuation theorem, is impossible if $c$ is irreducible. □

The above two lemmas thus complete the proof of Theorem 4.2.3. □

**4.2.10 Remark** From basic algebraic topology (cf. [37] §8.1), $\mathcal{G}_{l+1}$ has $H^1(M; \mathbb{Z})$ components.

### 4.3 Configurations and gauge group action: the case $t > 0$

Throughout this subsection, let $t > 0$. For the case of positive vortex-numbers we need an extra fibration construction. As we have seen in §3.2 for the case of $\mathbb{R}^3$, the Seiberg-Witten solutions with $t > 0$ and positive vortex numbers have $L^2$-unbounded curvature, thus to apply the analysis in the previous subsection we must subtract off some fixed configuration.

(A) THE CONFIGURATION SPACE $\mathcal{C}_l$.

Recall from §4.1 the definition of $\mathcal{A}_l^n$ as a neighborhood of the moduli space of the vortex equations. Let $G_{l+1}$ be the gauge group for the vortex equations; we take it to be a Banach Lie group modeled on $L_{l+1}^2(\mathbb{C}, i\mathbb{R})$. Denote $\Xi_l^n = \mathcal{A}_l^n/G_{l+1}$ and $\Xi := \coprod_{n=0}^{\infty} \Xi_l^n$.

**4.3.1 Definition** Given $a_1, a_2 \in \mathcal{A}_l^n$, let $c_{a_1 a_2}$ be a reference configuration associated with $a_1, a_2$.

For $t > 0$, let $\mathcal{C}_{l, a_1 a_2}$ be the following subspace in $\tilde{\mathcal{C}}$.

$$\mathcal{C}_{l, a_1 a_2} := \{c_{a_1 a_2} + e : e \in Y_{l, c_{a_1 a_2}}\}.$$

The configuration space $\mathcal{C}_l$ in this case is defined by $\mathcal{C}_l^n := \bigcup_{a_1, a_2 \in \mathcal{A}_l^n} \mathcal{C}_{l, a_1 a_2}$; $\mathcal{C}_l := \bigcup_{n \geq 0, n \in \mathbb{Z}} \mathcal{C}_l^n$.

Note again that the space $\mathcal{C}_{l, a_1 a_2}$ depends only on $(a_1, a_2) \in \mathcal{A}_l^n \times \mathcal{A}_l^n$, not the choice of the reference connection $c_{a_1 a_2}$.

**4.3.2 Notation** The above constructions depend on the parameters $t$ and $\epsilon$. We shall add the subscripts $t, \epsilon$ when we want to emphasize the dependence.

Generalizing Definition 4.1.1, we may define $H$, $X$, $Y$, $Z$ norms associated with any configuration in $\mathcal{C}_l$.

**4.3.3 Lemma** *The norms $Y_{l, c_1}$, $Y_{l, c_2}$ are commensurate for any $c_1 = (A_1, \psi_1) \in \mathcal{C}_{l, t_1}^n$, $c_2 = (A_2, \psi_2) \in \mathcal{C}_{l, t_2}^n$. Similarly, any pairs $X_{l+1, c_1}$, $X_{l+1, c_2}$ are also commensurate. In fact, $L_{k/A_1:\epsilon}^2$, $L_{k/A_2:\epsilon}^2$ are commensurate if $k \leq l$, or if $k = l+1$ and $k \geq 2$. Again an inequality of the form of (2.8) holds.*

*Proof.* We shall only present the proof for the $Y$-norms since the proof for the other norms are similar.

Let $e := c_2 - c_1$. We need to show that for all $\xi \in Y_{l, c_1}$, there is a constant $C(e)$ depending on $e$, such that

$$\left| \|\xi\|_{Y_{l, c_2}} - \|\xi\|_{Y_{l, c_1}} \right| \leq C(e) \|\xi\|_{Y_{l, c_1}}. \tag{4.21}$$

(4.21) holds in general if it holds for the following special cases: (i) when $c_1, c_2 \in C_{l, t_1}^n$ have the same left and right limits, (ii) when $c_1 \in \mathcal{C}_{l, t_1}^n$ and $c_2 \in \mathcal{C}_{l, t_2}^n$ are both reference configurations.

For case (i), note first that in this case $e \in Y_{l, c_1}$, and (4.21) follows from routine estimates using Lemma 2.2.12. For case (ii), it follows from the fact that $e \in L^{\infty}$ in this case. $\square$

Thus we shall drop $c$ and $A$ from the subscript.

**4.3.4 Remark** We fix the admissible pair $(g, \theta)$ in this section. However, by the same arguments we can see that the norms for different $(g, \theta)$ are commensurate. In fact, by the asymptotic property of $\theta$, the configuration space does not depend on the admissible pair either.

The previous lemma and the definition of $\mathcal{C}_l^n$ implies that $\mathcal{C}_l^n$ is the trivial product $Y_l \times (\mathcal{A}_l^n \times \mathcal{A}_l^n)$, and taking the right and left-hand limits gives the following fibration:
$$\mathcal{C}_{l,a_1a_2} \longrightarrow \mathcal{C}_l^n \xrightarrow{\partial_+ \times \partial_-} \mathcal{A}_l^n \times \mathcal{A}_l^n, \tag{4.22}$$
where $\partial_+, \partial_-$ are defined respectively by taking the $f \to +\infty$, $f \to -\infty$ limits of the configuration.

We take the topology of $\mathcal{C}_l^n$ to be the product topology. This makes $\mathcal{C}_l^n$ (and hence $\mathcal{C}_l$) a Banach manifold.

Note that Proposition 3.3.3 and Corollary 3.3.9 guarantee that all Seiberg-Witten solutions lie in $\mathcal{C}_l$ (modulo gauge transformations). As we will be looking at the moduli space of the Seiberg-Witten equations eventually, $\mathcal{C}_l$ is sufficient for our purpose.

(B) THE GAUGE GROUP ACTION AND THE QUOTIENT MANIFOLD.

In the $t > 0$ case, we define the subgroup, $\mathcal{G}_{l+1}$, of the total gauge group $\tilde{\mathcal{G}}$ in the following steps:

1. Embed $G_{l+1} \times G_{l+1} \to \tilde{\mathcal{G}}$ in the following way: Let $g_1$, $g_2$ be any two elements in $G_{l+1}$. Write $g_2 = e^{-iv}g_1 \in G_{l+1}$, where $v \in L^2_{l+1}(\mathbb{C}, \mathbb{R})$. Define an element $h_{g_1 g_2}$ in the total gauge group $\tilde{\mathcal{G}}$ out of $g_1, g_2$ as follows: let
$$h_{g_1 g_2} = e^{i\lambda_d(f)v} g_2 \tag{4.23}$$
on $M \backslash M_\Re$. ($\lambda_d$ is as in 3.3.3). It is not hard to see that we may extend $h_{g_1 g_2}$ over $M_\Re$ (continuously with respect to $g_1, g_1$) to define an element in $\tilde{\mathcal{G}}$.

2. Similarly to the $t = 0$ case, let $\mathcal{G}_{l+1,a_1a_2}$ be the stabilizer of the fiber $\mathcal{C}_{l,a_1a_2}$, and let $\mathcal{G}_{X_{l+1}}$ be the exponential of $X_{l+1}$. An analogue of lemma 4.2.7 shows that $\mathcal{G}_{X_{l+1}}$ is a Banach Lie group. Here in the $t > 0$ case, we can see directly that $\mathcal{G}_{X_{l+1}}$ is the component of $\mathcal{G}_{l+1,a_1a_2}$ containing unity: here it requires $\|(-2\nabla\xi, \xi\psi_y)\|_{Y_l} < \infty$ for $\xi \in T_e\mathcal{G}$ to stabilize $Y_l$, which is exactly the requirement for $\xi$ to be in $X_{l+1}$. Note that this implies that $\mathcal{G}_{l+1,a_1a_2}$ is independent of $a_1, a_2$.

   $\mathcal{G}_{X_{l+1}}$ acts on $\mathcal{C}_{a_1,a_2}$ smoothly. The relevant estimate here is
$$\|\xi\eta\|_{Y_l} \leq \|\xi\eta\|_{2,l;\epsilon}$$
$$\leq \|\xi\|_{3,l;\epsilon}\|\eta\|_{2,l} + (C\|\nabla\xi\|_2 + C'\|\nabla\xi\|_4)\|\nabla_A^l \eta\|_{2;\epsilon}$$
$$\leq C\|\xi\|_{X_{l+1}}\|\eta\|_{Y_l}. \tag{4.24}$$

3. Define the gauge group $\mathcal{G}_{l+1}$ to be the subgroup of $\tilde{\mathcal{G}}$ generated by $\mathcal{G}_{l+1,a_1a_2}$ together with the image of the embedding of $G_{l+1} \times G_{l+1}$ constructed in point 1 above.

$$\mathcal{G}_{l+1} \simeq G_{l+1} \times G_{l+1} \times \mathcal{G}_{l+1,a_1a_2} \tag{4.25}$$

Note that the definition of $\mathcal{G}_{l+1}$ does not depend on the choice of $\lambda_d$ or how we extend $h_{g_1g_2}$ over $M_\Re$ in point 1.

From the above construction, we observe that $\mathcal{G}_{l+1}$ and $\tilde{\mathcal{G}}$ have the same orbits in $\mathcal{C}_l \subset \tilde{\mathcal{C}}$. (Note that $\mathcal{C}_l$ is not $\mathcal{G}_{l+1}$ invariant.) Letting $\mathfrak{Q}_l^n \subset \tilde{\mathfrak{Q}}$ be the image of $\mathcal{C}_l^n$ under the quotient by $\mathcal{G}_{l+1}$, we have the following fibration

$$\mathfrak{Q}_{l,a_1a_2} := \mathcal{C}_{l,a_1a_2}/\mathcal{G}_{l+1,a_1a_2} \longrightarrow \mathfrak{Q}_l^n \overset{\partial_+ \times \partial_-}{\longrightarrow} \Xi_l^n \times \Xi_l^n, \tag{4.26}$$

where we have used the same notation $\partial_+ \times \partial_-$ to denote the map induced from the fibration map $\mathcal{C}_l^n \to \mathcal{A}_l^n \times \mathcal{A}_l^n$.

Let $\mathfrak{Q}_l := \bigcup_n \mathfrak{Q}_l^n$.

In the $t > 0$ case, lemma 4.1.8 generalizes to state that $d_c^*$ is an isomorphism between $Y_{l,c}$ and $Z_{l-1}$ for any $c \in \mathcal{C}_l$. Using the implicit function theorem as in §4.2 we thus obtain:

**4.3.5 Theorem** *Let $t > 0$. In this case, for each $y = (a_1, a_2) \in \Xi_l^n \times \Xi_l^n$, $\mathfrak{Q}_{l,y}$ is a Hilbert manifold modeled on*

$$Q_{l,c} := \Big\{f : f \in Y_l,\ d_c^* f = 0\Big\}, \tag{4.27}$$

*for any $c \in \mathcal{C}_{l,y}$.*

*The quotient manifold $\mathfrak{Q}_l^n$ is a fiber-bundle over $\Xi_l^n \times \Xi_l^n$ with fibers $\mathfrak{Q}_{l,y}$. Endowed with the product topology, $\mathfrak{Q}_l^n$ (and hence $\mathfrak{Q}_l$) is a Hilbert manifold modeled on*

$$\hat{Q}_{l,c} := \Big\{(b_1, b_2, f) : b_1, b_2 \in L_l^2(\mathbb{C}, T^*\mathbb{C} \oplus \mathbb{C}), f \in Y_l,$$
$$d_c^* f = 0, \delta_{a_1}^1 b_1 = \delta_{a_2}^1 b_2 = 0.\Big\} \tag{4.28}$$

*Here $c$ is an arbitrary configuration in $\mathcal{C}_l^n$, and $a_1, a_2$ are images of $c$ under the maps $\partial_+, \partial_-$ respectively.*

*$\delta^1$ is the gauge condition for vortex solutions given in (A.2).*

In fact, we can define a global gauge slice $(*)$:

**4.3.6 Definition** *We say that a configuration $c \in \mathcal{C}_l$ is in* the $*$-gauge *if*

$$c = c_{a_1a_2} + e, \tag{4.29}$$

*where: (i) $a_1$, $a_2$ are both in the gauge slice (**vor**). (cf. Appendix). (ii) $e \in Y_l$ satisfies the gauge constraint $d_{c_{a_1a_2}}^* e = 0$.*

# 5 Fredholm theory

Here we deal mainly with the $t > 0$ case. The $t = 0$ case is similar but simpler; it will be discussed briefly in §5.2.

## 5.1 The deformation operator $\mathcal{D}_c$

From the discussion in 3.3.2, over $U_\pm$, the trivialization of $E$ specified in 3.3.3 identifies $\mathcal{W}_0 = iT^*M \oplus S$ with the pull-back $i(\partial_\pm^* T\mathbb{C} \oplus \partial_\pm^* N) \oplus (\partial_\pm^* \mathbb{C} \oplus (\partial_\pm^* \mathbb{C} \otimes \partial_\pm^* T\mathbb{C}))$. It is often more convenient to regroup this as:

$$\mathcal{W}_0 = \partial_\pm^*(iT\mathbb{C} \oplus \mathbb{C}) \oplus \partial_\pm^*(iN \oplus (\mathbb{C} \otimes T\mathbb{C})). \tag{5.1}$$

For example, with respect to the same trivialization, a reference configuration has vanishing second component in the above over $U_\pm$.

**5.1.1 Definition** We define the space $\mathcal{V}_l$ as follows:

**(a)** In the $t = 0$ case, $\mathcal{V}_l := L_l^2(M, \mathcal{W}_0)$.

**(b)** In the $t > 0$ case, $\mathcal{V}_l \subset L_{l,loc}^2(M, \mathcal{W}_0)$,

$$\begin{aligned}\mathcal{V}_l &:= \Big\{ F(v_1, v_2) + e : e \in L_{l;\epsilon}^2(M, \mathcal{W}_0); \\ & \qquad\qquad v_1, v_2 \in L_l^2(\mathbb{C}, iN \oplus (\mathbb{C} \otimes T\mathbb{C})) \Big\} \tag{5.2} \\ &\simeq L_l^2(\mathbb{C}, iN \oplus (\mathbb{C} \otimes T\mathbb{C}))^{\oplus 2} \oplus L_{l;\epsilon}^2(M, \mathcal{W}_0), \tag{5.3}\end{aligned}$$

where

$$F(v_1, v_2)(x) := \Big(0, (1 - \chi_\Re)(x)\Big(\lambda_d(f)v_1(z_+) + (1 - \lambda_d)(f)v_2(z_-)\Big)\Big)$$

with respect to the decomposition (5.1). The above expression makes sense because it is supported on $M \backslash M_\Re$, over which we have the interpretation for $\mathcal{W}_0$ noted above. The direct sum decomposition (5.3) endows a norm on $\mathcal{V}_l$.

Consider the map $\mathcal{S} : \mathcal{C}_l \to \mathcal{V}_{l-1}$,

$$\mathcal{S}(A, \psi) := \Big( *(F_A - i\rho^{-1} \circ \sigma(\psi, \psi) - i\omega), \partial_A \psi \Big). \tag{5.4}$$

The zero set $\mathcal{S}^{-1}(0) \subset \mathcal{C}_l$ is the space of Seiberg-Witten solutions. Its image under the quotient by $\mathcal{G}_{l+1}$, $\mathcal{M}_l \subset \mathcal{Q}_l$, is *the moduli space* of Seiberg-Witten solutions.

The above definition of the moduli space seems to depend on the choice of $l$, $\epsilon$, but actually it *is* sufficiently universal considering the following fact:

**5.1.2 Lemma** *If the metric $g$ is $(k+5)$-admissible, $k \in \mathbb{Z}^+$, and the perturbation $\omega$ is given by (2.4) with $w \in L^2_{k-1}$, then in the $t > 0$ case, $\mathcal{M}_{l:\epsilon} = \mathcal{M}_{l':\epsilon'}$ for arbitrary $l, l' \in \mathbb{Z}^+$, $1 < l, l' \leq k$, and $\epsilon, \epsilon' \in [0, 3/2)$. Similarly, in the $t = 0$ case, $\mathcal{M}_l = \mathcal{M}_{l'}$.*

*Proof.* We present the proof for the $t > 0$ case. The $t = 0$ case is similar. Suppose $l \leq l'$, $\epsilon \leq \epsilon'$. Then by definition $\mathcal{M}_{l':\epsilon'} \subset \mathcal{M}_{l:\epsilon}$. The converse $\mathcal{M}_{l:\epsilon} \subset \mathcal{M}_{l':\epsilon'}$ is also true by elliptic bootstrapping (cf. Lemma 3.1.4) and Proposition 3.3.3. These together imply the lemma. $\square$

We will henceforth drop the subscripts and denote the moduli space by $\mathcal{M}$. Without loss of generality, *for the rest of this paper let $l = 2$, and drop the subscript $l$ in $\mathcal{C}_l$, $\mathcal{Q}_l$, etc.* For the case of positive vortex numbers, $\epsilon$ will be a fixed number $\epsilon \in (1, 3/2)$ (cf. Corollary 3.3.9 and Lemma 5.3.4 for the reasons of this choice). All discussion below may be easily generalized to higher $l$ cases.

To study the local structure of the moduli space, we need to consider the linearization of $\mathcal{S}$. Formally, the linearization of $\mathcal{S}$ at $c$ is a map from $\Gamma(M, \mathcal{W}_0)$ to itself given by:

$$D_c(a, \eta) = \left( * (da - 2i\rho^{-1} \circ \sigma(\eta, \psi)),\ \partial\!\!\!/_A \eta + \frac{\rho(a)}{2}\psi \right). \tag{5.5}$$

A straightforward computation shows that when $c$ is a Seiberg-Witten solution, $D_c$ fits into the elliptic complex:

$$i\Omega^0(M) \xrightarrow{d_c} i\Omega^1(M) \oplus \Gamma(S) \xrightarrow{D_c} i\Omega^1(M) \oplus \Gamma(S) \xrightarrow{d_c^*} i\Omega^0(M). \tag{5.6}$$

Let $\mathcal{W}_1$ be the bundle $i\mathbb{R} \oplus \mathcal{W}_0$ over $M$. It has a Euclidean metric by direct summing the one on $\mathcal{W}_0$ given by (4.14), and the standard one on $\mathbb{R}$.

By a standard construction, from this elliptic complex we obtain a formally $L^2$-selfadjoint operator (see lemma 5.2.3), $\mathcal{D}_c$, from $\Gamma(M, \mathcal{W}_1)$ to itself:

$$\mathcal{D}_c(\gamma, q) := (0, d_c\gamma) + (d_c^* q, D_c q). \tag{5.7}$$

Following (5.1), over $M \backslash M_\Re$, we shall also decompose $\mathcal{W}_1$ alternatively as

$$\left( K^{-1} \oplus E \right) \oplus \left( i\mathbb{C} \oplus (E \otimes K^{-1}) \right), \tag{5.8}$$

where the trivial bundle $\mathbb{C} := \mathbb{R} \oplus i\mathbb{R}$ consists of the $\mathbb{R}df$ component from $T^*M$, and the $i\mathbb{R}$ component in $\mathcal{W}_1$ complementing $\mathcal{W}_0$. Again over $U_\pm$ this can be identified with a pull-back bundle.

The goal of this section is to prove the Fredholmness of $\mathcal{D}_c$.

To make the definition of $\mathcal{D}_c$ precise, we need to specify its domain and range in different cases. Before doing so we need the following digression.

**5.1.3 A Decomposition of $\mathcal{D}_c$ on $\mathbb{R}^3$.** Let $c$ be a Seiberg-Witten solution on $\mathbb{R}^3$ with vortex number $n$. We have seen in §3.2 that $c = (v, 0)$ with respect to the decomposition (5.1).

It is straightforward to see that in this case $\mathcal{D}_c$ decomposes into the horizontal and normal parts $T'$, $N'$, each involving only the $x_3$ and $z$ variable respectively: $\mathcal{D}_c = i(T' + N')$, where

$$T' := \begin{pmatrix} \partial_3 & 0 \\ 0 & -\partial_3 \end{pmatrix}; \tag{5.9}$$

$$N' := \begin{pmatrix} 0 & -\Theta_v^* \\ \Theta_v & 0 \end{pmatrix} \tag{5.10}$$

as endomorphisms of $\Gamma(\mathcal{W}_1)$ with respect to (5.8).

It is clear from (5.10) that the kernel and cokernel of $N'$ both consist precisely of elements of the form $(k(z, x_3), 0)$ where $k(z, x_3) \in \text{Ker}\,\Theta_v \simeq \mathbb{C}^n$ for fixed $x_3$. ($n$ is the vortex number of $c$). So these can be regarded as sections of a trivial $\mathbb{C}^n$-bundle over $\mathbb{R}$ (parameterized by $x_3$), which we call $\mathcal{K}_c$. Let $L_\epsilon^2(\mathcal{K}_c)$ be the space of $L_\epsilon^2$ sections of $\mathcal{K}_c$ regarded as a subspace of $L_\epsilon^2(\mathcal{W}_1)$ as above. We let $L_\epsilon^2(\mathcal{W}_1; \mathcal{K}_c)$ denote the $L_\epsilon^2$-orthogonal complement of $L_\epsilon^2(\mathcal{K}_c)$ in $L_\epsilon^2(\mathcal{W}_1)$.

Let $\Pi_c$ denote the projection from $L_\epsilon^2(\mathcal{W}_1)$ onto $L_\epsilon^2(\mathcal{K}_c)$.

**5.1.4 Generalizations of this decomposition.** For our later partition of unity argument, we need to generalize the above definitions on $\mathbb{R}^3$ to the case where $f$, $c$ come from extending restrictions to $M \backslash M_R$ in the following sense.

Let $f$ be an $k$-admissible function on $M$ (recall that $k \geq l + 5$); we define $\tilde{f}$ to be a function which extends $f\big|_{M \backslash M_R}$ ($R \geq \Re$) over the whole $\mathbb{R}^3$, such that $\|\tilde{f} - x_3\|_{C_k} \leq \varepsilon(R)$. Note that $\tilde{f}$ is no longer harmonic, and we may choose $R$ large enough so that $\varepsilon$ is as small as we want. We can then use $\tilde{f}$ to decompose $\mathcal{W}_1$ as in (5.8) over the whole $\mathbb{R}^3$, and both $(z_+, f)$ and $(z_-, f)$ are good coordinate systems on $\mathbb{R}^3$.

Let $c$ be a configuration on $\mathbb{R}^3$ of the form $(v, 0)$ with respect to (5.8) and the trivialization on $E$ and $K^{-1}$ described above. For any fixed value of $f$, we assume that $v(z_\pm, f) \in \Gamma(\mathbb{C}, T\mathbb{C} \oplus \mathbb{C})$ is an approximate vortex solution. (The condition of being an approximate vortex solution is the same for either coordinate $z_\pm$ if $\varepsilon$ is small enough.)

In this case, $T'$, $N'$, $\mathcal{K}_c$, $\Pi_c$ have straightforward generalizations: $T', N'$ are given by the same formulas, except that $\partial_3$ is replaced by $|\nabla f|\partial_f$, and $\partial$ in (A.1) now is respect to the complex structure on $K^{-1}$. $\mathcal{K}_c$ is still a $\mathbb{C}^n$-bundle over $\mathbb{R}$; the fiber over the surface $f = C$ is $\text{Ker}\,\Theta_{v|_{f=C}}$. However, $\mathcal{D}'_c$ is only "approximately" decomposable:

$$\mathcal{D}'_c := i(T' + N') + R', \tag{5.11}$$

where $R'$ comes from the off-diagonal terms of the Spin$^c$ and Riemannian connections and involves only (matrix) multiplication; $\|R'\|_{C_{k-2}} \leq C\varepsilon$.

For general $M$ and $c$, $\Pi_c$ is partially defined as follows.

**5.1.5 Definition** Let $c \in \mathcal{C}_t(M)$, $t > 0$, and let $\lambda' := \lambda'_R$ as in Definition 2.2.10. For any $u \in L^2(\mathcal{W}_1)$ such that $u = (v, h)$ over $M \backslash M_\Re$ according to (5.8), we define

$$\Pi_c\Big((1 - \lambda'(f))(v,h)(x)\Big) := \begin{cases} (P_{\partial_+ c} v(z_+, f), 0) & \text{in the region where } f > R, \\ (P_{\partial_- c} v(z_-, f), 0) & \text{in the region where } f < -R, \\ 0 & \text{for the rest of } M, \end{cases}$$

where $P_{\partial_\pm c}$ denotes the projection from $L^2(\mathbb{C}, T\mathbb{C} \oplus \mathbb{C})$ to $\operatorname{Ker} \Theta_{\partial_\pm c}$ respectively.

In §5.3 case (i), the above definition of $\Pi_c$ agrees with that in 5.1.4 when the parameter $R$ is large enough.

**5.1.6 Domain and Range of $\mathcal{D}_c$: the case $t = 0$ and the case $t > 0, n = 0$.**
In the $t = 0$ case, $\mathcal{D}_c$ maps $V_2(M, \mathcal{W}_1)$ to $L_1^2(M, \mathcal{W}_1)$. [2]

In the $t > 0$ case, however, the domain and range of $\mathcal{D}_c$ have fibration structures due to the fibration of $\mathcal{C}$. We start with the $n = 0$ case ($n$ being the vortex number). Since we want $\operatorname{Dom}(\mathcal{D}_c)$ to include $T_c\mathcal{C}$, a natural choice of $\operatorname{Dom}(\mathcal{D}_c)$ is:

$$\Big\{q(v_1, v_2) + e(x) : v_1, v_2 \in L_2^2(\mathbb{C}, T\mathbb{C} \oplus \mathbb{C}), e \in L_{2:\epsilon}^2(M, \mathcal{W}_1)\Big\} \subset L_{loc}^2(M, \mathcal{W}_1), \tag{5.12}$$

where

$$q(v_1, v_2) := (1 - \chi_\Re)(x)\Big(\lambda_d(f)v_1(z_+) + (1 - \lambda_d)(f)v_2(z_-), 0\Big),$$

in the decomposition (5.8) according to the usual interpretation. Correspondingly, since $N'$ is off-diagonal, $\operatorname{Range}(\mathcal{D}_c)$ is given by the same formula as the right hand side of (5.2), but with $\mathcal{W}_0$ replaced by $\mathcal{W}_1$, and $N$ replaced by the trivial bundle $\mathbb{C}$.

It is obvious from the definition that in this case,

$$\operatorname{Dom}(\mathcal{D}_c) \simeq L_2^2(\mathbb{C}, T\mathbb{C} \oplus \mathbb{C})^{\oplus 2} \oplus L_{2:\epsilon}^2(M, \mathcal{W}_1);$$
$$\operatorname{Range}(\mathcal{D}_c) \simeq L_1^2(\mathbb{C}, \mathbb{C} \oplus \mathbb{C} \otimes T\mathbb{C})^{\oplus 2} \oplus L_{1:\epsilon}^2(M, \mathcal{W}_1).$$

If $a_1, a_2$ are respectively $\partial_+ c$ and $\partial_- c$, then in these decompositions,

$$\mathcal{D}_c(v_1, v_2, e) = (\Theta_{a_1} v_1, \Theta_{a_2} v_2, \mathcal{D}'_c e) + (0, 0, \mathfrak{X}_c), \tag{5.13}$$

---
[2] Recall that we have set $l = 2$.

where $\mathfrak{X}_c$ is the remainder term

$$\mathfrak{X}_c(v_1, v_2, e) := \sigma_{\mathcal{D}}(-d\chi_{\mathfrak{R}})\Big(\lambda_d v_1 + (1-\lambda_d)v_2, 0\Big) + \sigma_{\mathcal{D}}(d\lambda_d)\Big((1-\chi_{\mathfrak{R}})(v_1 - v_2), 0\Big)$$
$$+(1-\chi_{\mathfrak{R}})\Big(\lambda_d(\mathcal{D}_c - \mathcal{D}_{(a_1,0)})(v_1, 0) + (1-\lambda_d)(\mathcal{D}_c - \mathcal{D}_{(a_2,0)})(v_2, 0)\Big) \quad (5.14)$$

according to the usual interpretation. In the above, to simplify notation we identify $\lambda_d, v_1, v_2$ with their composition with $f, z_+, z_-$ respectively; also $\sigma_{\mathcal{D}}$ denotes the principal symbol of $\mathcal{D}_c$. $\mathcal{D}'_c$ is an operator of the form (5.7) from $L^2_{2:\epsilon}(\mathcal{W}_1)$ to $L^2_{1:\epsilon}(\mathcal{W}_1)$.

**5.1.7 Domain and Range of $\mathcal{D}_c$: the case $t > 0, n > 0$.** In this case, the spaces above are not good enough for the purpose of Fredholm theory, and we need to extend the domain of $\mathcal{D}_c$ a little. The trick is to replace the $e$ in the formula (5.12) with elements in a larger space, $K_{2,c}$, $L^2_{2:\epsilon} \subset K_{2,c} \subset L^2_{2:\epsilon-1}$. (From this point on we assume $\epsilon > 1$.) The moduli space corresponding to $K_{2,c}$ will include $\mathcal{M}_{2:\epsilon}$, yet be included by $\mathcal{M}_{2:\epsilon-1}$. By Lemma 5.1.2, this moduli space is isomorphic to $\mathcal{M}_{2:\epsilon} = \mathcal{M}_{2:\epsilon-1}$. Thus it does not matter to work with the $K_{2,c}$ norm. [3]

Below we define the space $K_{l,c}$ for general $l \in \mathbb{Z} \cup \{0\}$.

**Definition** Let $\epsilon \in (1, 3/2)$, and let $c := (A, \psi) \in \mathcal{C}^n$ for $n > 0$ be a configuration in the $t > 0$ case. Let $l \in \mathbb{Z}^+ \cup \{0\}$. Let $\Pi_c, \lambda'$ be as in Definition 5.1.5.
(1) For $\zeta \in C_0^\infty(M, \mathcal{W}_1)$, define the $K_{l,c}$-norm of $\zeta$ as

$$\|\zeta\|_{K_{l,c}} := \begin{cases} \|\nabla \zeta\|_{2,l-1:\epsilon} + \|(1-\Pi_c)((1-\lambda' \circ f)\zeta)\|_{2:\epsilon} + \|\zeta\|_{2:\epsilon-1} & \text{if } l > 0; \\ \|(1-\Pi_c)((1-\lambda' \circ f)\zeta)\|_{2:\epsilon} + \|\zeta\|_{2:\epsilon-1} & \text{if } l = 0. \end{cases}$$

Let $K_{l,c}$ be the completion of $C_0^\infty(M, \mathcal{W}_1)$ with respect to the above norm.
(2) Define the space $R_{l,c}$ to be the completion of $C_0^\infty(M, \mathcal{W}_1)$ with respect to the norm:

$\|\xi\|_{R_{l,c}} :=$
$$\begin{cases} \|\xi\|_{2,l:\epsilon} + \|\partial_f(\Pi_c(1-\lambda' \circ f)\xi)\|_{2:1+\epsilon} + \sum_{i=2}^l \|\partial_f^i(\Pi_c(1-\lambda' \circ f)\xi)\|_{2:2+\epsilon} & \text{if } l > 1, \\ \|\xi\|_{2,l:\epsilon} + \|\partial_f(\Pi_c(1-\lambda' \circ f)\xi)\|_{2:1+\epsilon} & \text{if } l = 1, \\ \|\xi\|_{2:\epsilon} & \text{if } l = 0. \end{cases}$$

From now on, we let $\mathcal{D}'_c$ denote the operator between $K_{2,c}$ and $L^2_{1:\epsilon}(\mathcal{W}_1)$ given by (5.7), and let $\mathcal{D}_c$ be the operator (5.13) from $L^2_2(\mathbb{C}, \mathbb{C}^2)^{\oplus 2} \oplus K_{2,c}$ to $L^2_1(\mathbb{C}, \mathbb{C}^2)^{\oplus 2} \oplus L^2_{1:\epsilon}(M, \mathcal{W}_1)$. The (formal $L^2_\epsilon$-) adjoint of $\mathcal{D}'_c$, $\mathcal{D}'^\dagger_c$, will map from $R_{2,c}$ to $K_{1,c}$. More generally, we may regard $\mathcal{D}'_c$ as an operator between $K_{l,c}$ and $L^2_{l-1:\epsilon}(\mathcal{W}_0)$, and $\mathcal{D}'^\dagger_c$ as an operator between $R_{l,c}$ and $K_{l-1,c}$. However unless otherwise specified, we take $l = 2$. Note that when we take $n = 0$ above, $K_{l,c}$, $R_{l,c}$ both reduce to $L^2_{l:\epsilon}$, and thus the domain and range of $\mathcal{D}_c$ defined in Definition 3 agree with the previous simpler definition for the $n = 0$ case.

---
[3] Alternatively, one can rework §4.3 using more complicated norms corresponding to $K_{2,c}$.

**5.1.8** Similar to the case in 5.1.6, $\mathcal{D}_c$ is more naturally interpreted as an operator from $\hat{K}_1$ to $\hat{L}^2_{1:\epsilon}(\mathcal{W}_1)$, which are defined below via an obvious modification of (5.12):

**Definition** (1) Let $\hat{K}_1$ be the space defined by (5.12), but with $L^2_{2:\epsilon}(M, \mathcal{W}_1)$ replaced by $K_{2,c}$ in the formula. In the notation there, we define the norm on $\hat{K}_1$ as:

$$\left\|q(a_1, a_2) + e\right\|_{\hat{K}_1} := \|a_1\|_{L^2_2(\mathbb{C}, T\mathbb{C} \oplus \mathbb{C})} + \|a_2\|_{L^2_2(\mathbb{C}, T\mathbb{C} \oplus \mathbb{C})} + \|e\|_{K_{2,c}}. \qquad (5.15)$$

(2) Let $\hat{K} := \hat{K}_1 \cap \Gamma(M, \mathcal{W}_0)$, often regarded as a subspace of $L^2_{2,loc}(\mathcal{W}_0)$. It inherits a norm from $\hat{K}_1$.

(3) Let $V$ be a Euclidean or hermitian bundle over $M$ with a chosen trivialization on $M\setminus M_\Re$, then again we may identify $V = \partial^*_\pm V'$ on $U_\pm$ correspondingly, where $V'$ is the trivial $\mathbb{R}^r$ or $\mathbb{C}^r$ bundle over $\mathbb{C}$ (depending on whether $V$ is Euclidean or hermitian; $r$ is the rank of $V$). $\hat{L}^p_{k:\epsilon}(V)$ will denote the space given again by (5.12), but with $a_1, a_2$ now in $L^p_k(\mathbb{C}, V')$ and $e \in L^p_k(M, V)$; and $q(a_1, a_1)$ is given by the same formula, but *without* the second component (i.e. 0). The norm on this space is defined in the usual way.

**5.1.9 Notation** In particular, when $V$ is a trivial $\mathbb{R}$-bundle, we often denote $\hat{L}^2_{3:\epsilon}(V)$ by $\hat{X}$ in light of Lemma 4.1.4. We shall sometimes also abuse notation and use $\hat{L}^2_{k:\epsilon}(\mathcal{W}_1)$ to denote the analogue of (5.12), with $e$ there replaced by an element in $L^2_{k:\epsilon}$. Similarly, we shall write $\hat{L}^2_{k:\epsilon}(\mathcal{W}_0)$ for $\hat{L}^2_{k:\epsilon}(\mathcal{W}_1) \cap \Gamma(\mathcal{W}_0)$.

We shall often confuse elements in $\Gamma(\mathcal{W}_0)$ with their images in $\Gamma(\mathcal{W}_1)$ under the natural embedding $\Gamma(\mathcal{W}_0) \hookrightarrow \Gamma(\mathcal{W}_1)$.

Now we are ready to state the main result in this section. Though we shall only need the result for $l = 2$, we state the Fredholmness result for general $l$ because we shall use an induction argument in $l$.

**5.1.10 Theorem** *Let $t > 0$ and $l \in \mathbb{Z}^+$. For $k \geq max(l, 2)$, let $c \in \mathcal{C}^n_{l:\epsilon}$. Then:*
*(0) $\mathcal{D}'_c$ and $\mathcal{D}'^\dagger_c$ are bounded (uniformly in c) with respect to the domains and ranges specified in 5.1.6, 5.1.7.*
*(1) $\mathcal{D}'_c$ is Fredholm between $K_{l,c}$ and $L^2_{l-1:\epsilon}(\mathcal{W}_1)$ for $\epsilon \in (1, 3/2)$ and $l \in \mathbb{Z}^+$. Its index does not depend on $l$ or $\epsilon$.*
*(2) The formal $L^2_\epsilon$-adjoint of $\mathcal{D}'_c$, $\mathcal{D}'^\dagger_c$, is also Fredholm of index $-\mathrm{Ind}(\mathcal{D}'_c)$ between $R_{l,c}$ and $K_{l-1,c}$.*

(3) If $n = 0$, we have in addition that for $\epsilon = 0$ (unweighted case), $\mathcal{D}'_c$ is Fredholm between $L^2_l(\mathcal{W}_1)$ and $L^2_{l-1}(\mathcal{W}_1)$ and has index 0; in fact, $\mathcal{D}'_c$ is a self-adjoint operator between $L^2_1$ and $L^2$.

The version for the case $t = 0$ may be found in the end of §5.2.

The proof of Theorem 5.1.10 is divided into two parts. The first part is contained in §5.2, and deals with the simpler case $\epsilon = 0, n = 0$. The approach adopted in this part is different from the $\epsilon > 0$ cases and may be readily adapted to the $t = 0$ case. The second part occupies §5.3–5.4, and deals with the $\epsilon > 0$ cases via excision and a separation of variables argument which appeared in [43].

The extension to $\mathcal{D}_c$ is a consequence of Theorem 5.1.10, and is presented in §5.5.

## 5.2 Fredholmness via admissibility

In this subsection, we prove Theorem 5.1.10 in *the unweighted case (i.e. $t > 0$, $\epsilon = 0$—this requires the vortex number $n = 0$ also)*. The same method yields the analogous result for the case $t = 0$ (Theorem 5.2.5 below).

First, we note that the boundedness of $\mathcal{D}'_c$ and $\mathcal{D}'^{\dagger}_c$ in this case follows from routine estimates.

Next, we observe that the $l > 1$ case of the Fredholmness assertion may be reduced to the $l = 1$ case by elliptic regularity, since a straightforward computation shows that

$$\|\nabla^{k+1}\xi\|^2_2 \leq C(\|\nabla^k \mathcal{D}'_c \xi\|^2_2 + \|\xi\|^2_{2,k}) \tag{5.16}$$

for $k < l$, $c \in \mathcal{C}_{l:0}$ and $\xi \in L^2_{k+1}(\mathcal{W}_1)$. We shall therefore take $l = 1$ *for the rest of this subsection*.

To show Fredholmness in this situation, we follow the index theory on Euclidean spaces developed in [39].

Similar to [39] and [20], we define "admissible operators" as follows.

**5.2.1 Definition** If $V$ is a Euclidean/hermitian bundle over an MEE $M$, and $D$ is a first order, elliptic operator from $\Gamma(V)$ to $\Gamma(V)$. Let $D^*$ be the formal $L^2$-adjoint of $D$. We say that $D$ is *admissible* if $D$ satisfies:

(1) Both quadratic forms $\langle D(\cdot), D(\cdot) \rangle_2$, and $\langle D^*(\cdot), D^*(\cdot) \rangle_2$ may be written in the form

$$\langle \cdot, \cdot \rangle_Q + \langle \cdot, \mathcal{R}(\cdot) \rangle_2, \tag{5.17}$$

where the first term is an inner product on $\Gamma(V)$ defined as

$$\langle \cdot, \cdot \rangle_Q := \langle \nabla_A(\cdot), \nabla_A(\cdot) \rangle_2 + \langle q(\cdot), q(\cdot) \rangle_2,$$

where $A$ is a metric-preserving connection on $V$, and $q$, $\mathcal{R}$ are some endomorphisms of $L^2(V)$.

(2) $\mathcal{R}$ in (5.17) satisfies: for any $\varepsilon > 0$, there exists $R \in \mathbb{R}^+, R \geq \mathfrak{R}$ such that if $\chi_R$ is a characteristic function of $M \backslash M_R$, then

$$|\langle \zeta, (1 - \chi_R)\mathcal{R}(\eta)\rangle_2| \leq \varepsilon \langle \zeta, \zeta \rangle_Q^{1/2} \langle \eta, \eta \rangle_Q^{1/2}. \tag{5.18}$$

Denote the $Q, q, \mathcal{R}$ in (5.17) corresponding to $\langle D(\cdot), D(\cdot)\rangle_2$ by $Q_D, q_D, \mathcal{R}_D$. Likewise denote those corresponding to $\langle D^*(\cdot), D^*(\cdot)\rangle_2$ by $Q_{D^*}, q_{D^*}, \mathcal{R}_{D^*}$.

Define the $Q_D$-norm for $\zeta \in \Gamma(V)$ by

$$\|\zeta\|_{Q_D} := \langle \zeta, \zeta \rangle_{Q_D}^{1/2},$$

and let $Q_D$ also denote the completion of $C_0^\infty(V)$ with respect to this norm. Define $Q_{D^*}$ likewise. Then

**5.2.2 Proposition** *If $D$ is an admissible operator, then $D$ is Fredholm between $Q_D$ and $L^2(V)$.*

For a proof see [19] Proposition 3.6 and [39] Proposition 7.2.

To apply the above Proposition to our situation, we need to check (5.17), (5.18) for both $\mathcal{D}'_c$ and $\mathcal{D}'^*_c$; the following lemma however shows that they need only be checked for $\mathcal{D}'_c$.

**5.2.3 Lemma** *$\mathcal{D}_c$ is formally $L^2$-self-adjoint.*

The proof follows from direct computation.

We now check the admissibility of $\mathcal{D}'_c$. Let $\zeta := (f', a', \eta'), \xi := (f, a, \eta) \in i\Omega^0(M) \oplus i\Omega^1(M) \oplus \Gamma(S)$; note that

$$\langle \mathcal{D}'_c \zeta, \mathcal{D}'_c \xi \rangle = \langle df', df \rangle + \langle da', da \rangle + \langle d^*a', d^*a \rangle$$
$$+ 1/2 \int_M (a' \cdot a + f'f)|\psi|^2 + \int_M \left(2 \operatorname{Tr}[\sigma(\psi, \eta')\sigma(\psi, \eta)] + \operatorname{Im}(\psi \cdot \eta) \operatorname{Im}(\psi \cdot \eta')\right)$$
$$+ 2\langle \partial_A \eta', \partial_A \eta \rangle + \text{cross-term1} + \text{cross-term2}, \tag{5.19}$$

where the cross terms are given by

$$\text{cross-term1}$$
$$= \langle -\partial M, -2i\sigma(\psi, \eta) + i \operatorname{Im} \psi \cdot \eta \rangle + \langle M\psi, \partial_A \eta \rangle$$
$$= -\operatorname{Re} \int (\operatorname{Tr}[(\eta^\dagger \otimes \psi)\partial M]) + \langle M\psi, \partial_A \eta \rangle$$
$$= -\operatorname{Re} \int \operatorname{Tr}(M\partial(\eta^\dagger \otimes \psi)) + \operatorname{Tr}(M(\psi^\dagger \otimes \partial_A \eta))$$
$$= -\operatorname{Re} \int \operatorname{Tr}(M[\eta^\dagger \otimes \partial_A \psi]) + 2\operatorname{Re}\langle a^* \cdot \nabla_A \psi, \eta \rangle + i \operatorname{Im}\langle M\psi, \partial_A \eta \rangle$$
$$= -\operatorname{Re}\langle M\eta, \partial_A \psi \rangle + 2\operatorname{Re}\langle a^* \cdot \nabla_A \psi, \eta \rangle + i \operatorname{Im}\langle M\psi, \partial_A \eta \rangle, \tag{5.20}$$

and similarly

$$\text{cross-term2} = -\operatorname{Re}\langle \partial\!\!\!/_A\psi, N\eta'\rangle + 2\operatorname{Re}\langle a\eta', \nabla_A\psi\rangle$$
$$+ i\operatorname{Im}\langle \partial\!\!\!/_A\eta', N\psi\rangle,$$

where $N := \rho(a) + fI$; $M := \rho(a') + f'I$. By the Weitzenböck formula we have among the diagonal terms

$$\langle da', da\rangle + \langle d^*a', d^*a\rangle + 2\langle \partial\!\!\!/_A\eta', \partial\!\!\!/_A\eta\rangle$$
$$= \langle \nabla a', \nabla a\rangle + \langle \nabla_A\eta', \nabla_A\eta\rangle + 2\langle \eta', \kappa\eta\rangle$$
$$- \langle \eta', \rho(F_A)\eta\rangle, \tag{5.21}$$

where $\kappa$ is a linear combination of the components of the curvature of the manifold $M$.

We now see that $\langle \mathcal{D}'_c\zeta, \mathcal{D}'_c\xi\rangle$ is of the required form (5.17), with

$$q(f, a, \eta) := \left(-2i\sigma(\psi, \eta) + i\operatorname{Im}\eta \cdot \psi, \frac{\rho(a) + fI}{2}\psi\right);$$
$$\langle \zeta, \mathcal{R}\xi\rangle := \langle \eta', (2\kappa - \rho(F_A))\eta\rangle + \text{cross-terms}.$$

Furthermore, the $Q_{\mathcal{D}}$-norm here is commensurate with the $L_1^2$-norm as some computation shows

$$\|\zeta\|_{Q_{\mathcal{D}}}^2 = \|\nabla_A\zeta\|_2^2 + \frac{1}{2}\||\psi|\zeta\|_2^2.$$

It follows from lemma 4.1.4 that

$$\mu'\|\zeta\|_{2,1}^2 \geq \|\zeta\|_{Q_{\mathcal{D}}}^2 \geq \mu\|\zeta\|_{2,1}^2.$$

It remains to show that $\mathcal{R}$ satisfies

$$|\langle (1 - \chi_R)\zeta, \mathcal{R}\xi\rangle| \leq \varepsilon\|\zeta\|_{2,1}\|\xi\|_{2,1}$$

for some $R$ depending on $\varepsilon$, for any small $\varepsilon > 0$. As $(1 - \chi_R)\mathcal{R}$ is hermitian it suffices to show that for all $\zeta \in L_{1:\epsilon}^2$,

$$|\langle (1 - \chi_R)\zeta, \mathcal{R}\zeta\rangle| \leq \varepsilon\|\zeta\|_{2,1}^2.$$

To show this, note that from (5.21), (5.20)

$$\langle \zeta, (1-\chi_R)\mathcal{R}\zeta\rangle = \langle \eta, (1-\chi_R)(2\kappa - \rho(F_A))\eta\rangle$$
$$- \operatorname{Re}(\langle N\eta, (1-\chi_R)\partial\!\!\!/_A\psi\rangle + \int 2(1-\chi_R)\operatorname{Re}(a^* \cdot \langle \nabla_A\psi, \eta\rangle)$$
$$+ \text{similar terms from cross-term2}.$$

As $\kappa$ is compactly supported (since $M$ is an MEE), we have

$$|\langle \eta, \kappa(1-\chi_R)\eta\rangle| \leq \|\eta\|_3 \|\kappa(1-\chi_R)\|_2 \|\nabla_A \eta\|_2 \leq \varepsilon \|\eta\|_{2,1}^2$$

for any small $\varepsilon > 0$, if we choose $R$ large enough. For the term $\langle \eta, (1-\chi_R)\rho(F_A), \eta\rangle$, take $R > \Re$ so that we can work on $\mathbb{R}^3$, and note that

$$|\langle \eta, \rho(F_A)(1-\chi_R)\eta\rangle| \leq \|(1-\chi_R)\rho(F_A)\|_2 \|\eta\|_{2,1}^2$$

can also be made arbitrarily small by taking $R$ large enough. The terms involving $\nabla_A \psi$ and $\bar{\partial}_A \psi$ can be bounded in the same manner. *Note that here we used the $L^2$-integrability of $F_A$ and $\nabla_A \psi$ for $n = 0$ solutions (since the only vortex solution on $\mathbb{C}$ with vortex number 0 consists of flat connection and constant Higgs field), which is not true for $n > 0$.* This concludes the proof of Fredholmness in the $l = 1$ case, which actually implies the Fredholmness for all $l$ as explained.

**5.2.4 The calculation of index** in this case ($\epsilon = 0$, $t > 0$) follows from the next lemma. (The lemma calculates the index for the $l = 1$ case; since the index does not depend on $l$, this gives the index for any $l$.)

**Lemma** *In the notation of Theorem 5.1.10, $\mathcal{D}'_c$ is a self-adjoint operator between $L_1^2(\mathcal{W}_1)$ and $L^2(\mathcal{W}_1)$. Thus $\mathrm{Ind}(\mathcal{D}'_c) = 0$.*

*Proof.* To prove that $\mathcal{D}'_c$ is actually self-adjoint, it suffices to show that $\mathrm{Dom}(\mathcal{D}'^t_c) \subset \mathrm{Dom}(\mathcal{D}'_c)$, where $\mathcal{D}'^t_c$ is the adjoint of $\mathcal{D}'_c$ in the sense of [16]. Let $\zeta \in \mathrm{Dom}(\mathcal{D}'^t_c)$, that is, $\zeta \in L^2(\mathcal{W}_1)$, with $\mathcal{D}'^t_c \zeta$ a distribution in $L^2(\mathcal{W}_1)$.

As $\mathcal{D}'_c$ is proved to be Fredholm, $\mathrm{Range}(\mathcal{D}'_c)$ is closed and we can approximate $\zeta$ by elements in $C_0^\infty$. But for these elements $\mathcal{D}'^t_c \zeta = \mathcal{D}'_c \zeta$, and we know from (5.19) that

$$\|\nabla_A \zeta\|_2 \leq \|\mathcal{D}'_c \zeta\|_2 + C\|\zeta\|_2 < \infty. \tag{5.22}$$

This shows that $\mathrm{Dom}(\mathcal{D}'^t_c) \subset L_1^2 = \mathrm{Dom}(\mathcal{D}'_c)$ and the lemma is proved. □

This proves assertion (3) of Theorem 5.1.10 in this case, and thus concludes the proof of Theorem 5.1.10 in the case $\epsilon = 0$ completely.

The case $t = 0$ only needs minor modifications of the above arguments.

**5.2.5 Theorem** *Let $c$ be a reference configuration for $t = 0$. Then $\mathcal{D}_c$ is uniformly bounded (in $c$) and is Fredholm between $V_l(\mathcal{W}_1)$ and $L_{l-1}^2(\mathcal{W}_1)$. Its index is zero.*

The index calculation in the above theorem follows from a standard excision argument (See also (5.38)–(5.41) and the proof of Proposition 6.3.1 assertion 2), and the elementary fact that elliptic operators on compact odd-dimensional manifolds have zero index. (See for example [21]). □

## 5.3 Boundedness and Fredholmness of $\mathcal{D}'_c$: the weighted case

This subsection contains the first half of the proof of Theorem 5.1.10 in the case $\epsilon > 0$: we prove assertion (0) and assertion (1) in the case $l = 1$. In fact, we shall only show the boundedness of $\mathcal{D}'_c$, since the proof for $\mathcal{D}'^\dagger_c$ is similar. The second half is in the next subsection.

We shall start with special cases, then generalize step by step: Case (i) deals with some standard configuration $c$ on $\mathbb{R}^3$ via a separation-of-variables argument; case (ii) generalizes this to certain special configurations $c$ over a general MEE by excision[4] ; Case (iii) is the completely general case obtained via perturbing case (ii).

CASE (I). *Let $f$ and $c = (v, 0)$ be as in 5.1.4, with $c$ satisfying the following additional condition: there exists some real number $d \geq \Re$, so that $\partial_f v = 0$ wherever $|f| > d$, and $\forall C \in \mathbb{R}, \|\partial_f v\big|_{f=C}\|_{1,\infty} < \varepsilon$ for a small enough number $\varepsilon$ independent of $c$. This condition will be useful in lemma 5.3.3.*

**5.3.1 To show that $\mathcal{D}'_c$ is bounded in case (i),** it suffices to show that $T' + N'$ is bounded. Given any $\zeta \in K_{2,c}$, decompose $\zeta$ as $\zeta = u + w$, where $u \in L^2_{2:\epsilon}(\mathcal{W}_1; \mathcal{K}_c)$; $w \in L^2_{2:\epsilon-1}(\mathcal{K}_c)$. Note that in this case, $\|\zeta\|_{K_{2,c}}$ is bounded above and below by multiples of $\|u\|_{2,2:\epsilon} + \|\nabla_A w\|_{2,1:\epsilon} + \|w\|_{\epsilon-1}$.

A direct computation shows that $\|(T' + N')u\|_{2,1:\epsilon} \leq C\|u\|_{2,2:\epsilon}$ and Since $(T' + N')w = T'w$, $\|(T' + N')w\|_{2,1:\epsilon} \leq \|\nabla_A w\|_{2,1:\epsilon}$. These together imply that $\|(T' + N')\zeta\|_{2,1:\epsilon} \leq C'\|\zeta\|_{K_{2,c}}$ for a $c$-independent constant $C'$. $\square$

### 5.3.2 Fredholmness of $\mathcal{D}'_c$ with $l = 1$ in case (i).

**Proposition** *Let $c$ be as above, and $\epsilon \in (1, 3/2)$. Then $\mathcal{D}'_c$ is Fredholm between $K_{1,c}$ and $L^2_{:\epsilon}$, of index $-2n$, and has trivial kernel.*

*Proof.* The first step is the following standard lemma which defines a partial right inverse for $\mathcal{D}'_c$.

**5.3.3 Lemma** *Let $c, \epsilon$ be as above. Then for any $h \in L^2_\epsilon(\mathcal{W}_1)$, there is a pair of $b(h) \in L^2_\epsilon(\mathcal{K}_c)$ and $u(h) \in L^2_{2:\epsilon}(\mathcal{W}_1; \mathcal{K}_c)$ such that $\mathcal{D}'_c \mathcal{D}'^\dagger_c u(h) - h = b(h)$. Furthermore, there are $c$-independent constants $C, C'$, such that*

$$\|u(h)\|_{2,2:\epsilon} \leq C\|h\|_{2:\epsilon}; \tag{5.23}$$

$$\|b(h)\|_{2:\epsilon} \leq C'\|h\|_{2:\epsilon}. \tag{5.24}$$

---
[4]however this method does not tell us about the index of the Fredholm operator. To obtain information on the index, further refinement using gluing theory is needed. See Proposition 6.3.1 for an example.

*Proof.* Try minimizing the functional $\mathfrak{f}: L^2_\epsilon(\mathcal{W}_1; \mathcal{K}_c) \to \mathbb{R}$,

$$\mathfrak{f}(\mu) := \frac{1}{2}\|\mathcal{D}'^\dagger \mu\|^2_{2:\epsilon} - \langle \mu, h \rangle_{2:\epsilon}. \tag{5.25}$$

It is easy to see that $\mathfrak{f}$ is bounded and convex, so to prove that there exists an unique minimum, it suffices to establish a coercive lower bound.

We first find a lower bound for

$$\|(N'^\dagger + T'^\dagger)\mu\|^2_{2:\epsilon} = \|N'^\dagger \mu\|^2_{2:\epsilon} + \|T'^\dagger \mu\|^2_{2:\epsilon} + \text{cross-terms}.$$

The cross terms above may be bounded as follows. (i) Notice that $T'^\dagger = -T' + m$, where $m = -i\varsigma^{-2\epsilon}|\nabla f|^{-1}\partial_f(|\nabla f|\varsigma^{2\epsilon})\gamma_3$ (cf. (3.1) for $\gamma_3$), and $N'^\dagger = -N'$. So $\sigma_{T'}\sigma_{N'^\dagger} + \sigma_{N'}\sigma_{T'^\dagger} = 0$, where $\sigma_T$ means the principal symbol of the operator $T$. Thus the cross terms vanish at the symbol level. (ii) Since $\|m\|_\infty$ is small by Definition 2.2.10, the other non-vanishing cross terms in $\|\mathcal{D}'^\dagger_c \mu\|^2_{2:\epsilon}$ can be estimated such as:

$$|\langle \mu, N'm\mu \rangle_{2:\epsilon}| \leq \varepsilon_0(\|N'\mu\|^2_{2:\epsilon} + \|u\|^2_{2:\epsilon}),$$

where $\varepsilon_0$ is a small positive number. (iii) The norms of terms involving $\partial_f c$ can be bounded by $\varepsilon_1 \|\mu\|^2_{2:\epsilon}$ since $|\partial_f c|$ is small by the definition of $c$.

Using the fact that $\|\bar{\partial}_A \alpha\|_{2:\epsilon}$ is small (since $(A,\alpha) \in \mathcal{A}^n$ is an approximate vortex solution) and the orthogonality of $\mu$ to $L^2_\epsilon(\mathcal{K}_c)$, we then are able to estimate

$$\|\mathcal{D}'^\dagger_c \mu\|^2_{2:\epsilon} \geq \frac{1}{2}\|(T' + N')^\dagger \mu\|^2_{2:\epsilon} - \varepsilon_0 \|\mu\|_{2:\epsilon'}$$
$$\geq C(\varepsilon\|(|\nabla f|\partial_f)^\dagger \mu\|^2_{2:\epsilon} + \|\nabla_{z,A}\mu\|^2_{2:\epsilon} + \|\alpha\mu\|^2_{2:\epsilon}) - \varepsilon_2 \|\mu\|^2_{2:\epsilon} \tag{5.26}$$
$$\geq C_1(\varepsilon\||\nabla f|\partial_f \mu\|^2_{2:\epsilon} + \|\nabla_{z,A}\mu\|^2_{2:\epsilon} + \|\mu\|^2_{2:\epsilon}) - C'\varepsilon\||\partial_f m|^{1/2}\mu\|^2_{2:\epsilon} \tag{5.27}$$
$$\geq C''(\|\nabla_A \mu\|^2_{2:\epsilon} + \|\mu\|^2_{2:\epsilon}), \tag{5.28}$$

where $\varepsilon$ is a small positive number depending only on $\varsigma$, $C''$ is $c$-independent, and from (5.26) to (5.27) we have used a generalization of sublemma 4.1.5. Finally, we have

$$|\langle \mu, h \rangle_{2:\epsilon}| \leq \|\mu\|_{2:\epsilon}\|h\|_{2:\epsilon} \leq \varepsilon'\|\mu\|_{2:\epsilon} + C_2\|h\|_{2:\epsilon} \tag{5.29}$$

for a small positive $\varepsilon'$. Putting (5.26)–(5.29) together, we thus obtained the desired coercive lower bound for the functional $\mathfrak{f}$ for $\mu \in L^2_\epsilon(\mathcal{W}_1, \mathcal{K}_c) \cap L^2_{1:\epsilon}(\mathcal{W}_1)$, and it has a unique minimum at, say, $u$, which satisfies

$$\langle \mathcal{D}'^\dagger_c \mu, \mathcal{D}'^\dagger_c u \rangle_{2:\epsilon} - \langle \mu, h \rangle_{2:\epsilon} = 0, \tag{5.30}$$

$\forall \mu \in L^2_\epsilon(\mathcal{W}_1, \mathcal{K}_c)$; hence there exists a $b \in L^2_\epsilon(\mathcal{K}_c)$ such that

$$\mathcal{D}'_c \mathcal{D}'^\dagger_c u - h = b. \tag{5.31}$$

We still need to verify the estimates (5.23), (5.24) for $u$ and $b$. Substituting $\mu = u$ into (5.30), we have the following inequality from (5.28), (5.29):

$$\|\nabla_A u\|_{2:\epsilon}^2 + \|u\|_{2:\epsilon}^2 \leq \zeta \|h\|_{2:\epsilon}^2 \tag{5.32}$$

for some $c$-independent constant $\zeta$. Now (5.30), (5.32) ensure that the projection of $\mathcal{D}'_c \mathcal{D}'^\dagger_c u$ onto $L^2_\epsilon(\mathcal{W}_1, \mathcal{K}_c)$ is finite. To show that $\mathcal{D}'_c \mathcal{D}'^\dagger_c u$ is $L^2_\epsilon(\mathcal{W}_1)$-integrable, we estimate $\langle \mathcal{D}'^\dagger_c \mu, \mathcal{D}'^\dagger_c u \rangle_{2:\epsilon}$ for all $\mu \in L^2_\epsilon(\mathcal{K}_c) \cap L^2_{1:\epsilon}(\mathcal{W}_1)$. First note that because $\Pi_c \mu = \mu$, $\Pi_c u = 0$, and $N'^\dagger \mu = 0$, we can write (cf. (5.11) for notation)

$$\langle \mathcal{D}'^\dagger_c \mu, \mathcal{D}'^\dagger_c u \rangle_{2:\epsilon}$$
$$= \langle T'^\dagger \Pi_c \mu, (T'^\dagger + N'^\dagger) u \rangle_{2:\epsilon} + \text{terms involving } R'$$
$$= \langle [T'^\dagger, \Pi_c] \mu, \mathcal{D}'^\dagger_c u \rangle_{2:\epsilon} + \langle T'^\dagger \mu, [\Pi_c, T'^\dagger] u \rangle_{2:\epsilon} + \text{terms involving } R',$$

the absolute value of which may be bounded using the definition of $c$ by

$$C \|\mu\|_{2:\epsilon} \|u\|_{2,1:\epsilon} \leq C' \|\mu\|_{2:\epsilon} \|h\|_{2:\epsilon}.$$

This implies that the $L^2_\epsilon(\mathcal{K}_c)$ component of, and hence the whole of $\|\mathcal{D}'_c \mathcal{D}'^\dagger_c u\|_{2:\epsilon}$ also, is bounded by $\|h\|_{2:\epsilon}$. By (5.31) this implies (5.24).

(5.23) is verified by (5.32) and the following estimate:

$$\|\nabla_A \nabla_A u\|_{2:\epsilon}^2 \leq \|\mathcal{D}'_c \mathcal{D}'^\dagger_c u\|_{2:\epsilon}^2 + C\|u\|_{2,1:\epsilon}^2 \leq C\|h\|_{2:\epsilon}^2.$$

$\square$

Define $P_c(h) := \mathcal{D}'^\dagger_c u(h)$; this is the desired partial right inverse. Note that when $c$ is a Seiberg-Witten solution on $\mathbb{R}^3$ given in §3.2, $P_c$ maps from $L^2_\epsilon(\mathcal{W}_1; \mathcal{K}_c)$ to $L^2_{1:\epsilon}(\mathcal{W}_1; \mathcal{K}_c)$, and $\mathcal{D}'_c P_c = 1 - \Pi_c$. This is not true for a general $c$ in case (i), since in general $[\Pi_c, \mathcal{D}'_c] \neq 0$.

The previous lemma then reduces the problem to (partially) inverting $\mathcal{D}'_c$ on $L^2_\epsilon(\mathcal{K}_c) \subset L^2_\epsilon(\mathcal{W}_1)$. More precisely, we decompose

$$h = u + b$$

with $u \in L^2_\epsilon(\mathcal{W}_1; \mathcal{K}_c); b \in L^2_\epsilon(\mathcal{K}_c)$. We look for a $q$ such that

$$\mathcal{D}'_c q = h \tag{5.33}$$

(at least partially). Write $q := v + w$, $v \in L^2_\epsilon(\mathcal{W}_1; \mathcal{K}_c)$, $w \in L^2_{\epsilon-1}(\mathcal{K}_c)$. (Note that by our definition of $c$, $K_{2,c} = L^2_\epsilon(\mathcal{W}_1; \mathcal{K}_c) \oplus L^2_{\epsilon-1}(\mathcal{K}_c)$). Projecting (5.33) onto $L^2_\epsilon(\mathcal{W}_1; \mathcal{K}_c)$ and $L^2_\epsilon(\mathcal{K}_c)$ respectively, we have:

$$(1 - \Pi_c) \mathcal{D}'_c v = i[\Pi_c, T']w + u + (1 - \Pi_c) R' w; \tag{5.34}$$
$$i[\Pi_c, T']v + i\Pi_c T' w + \Pi_c R' q = b. \tag{5.35}$$

From (5.34) we get an expression of $v$ in terms of $w$ using lemma 5.3.3 (i.e., $v = P_c(i[\Pi_c, T']w + u + (1 - \Pi_c)R'w))$; substituting into (5.35) we obtain an equation for $w$:

$$|\nabla f|\partial_f w - \gamma_3(i[\Pi_c, iT']w + i[\Pi_c, T']v - b + \Pi_c R'(v + w)) = 0. \qquad (5.36)$$

To solve (5.36), we identify elements in $L^2_\epsilon(\mathcal{K}_c)$ with $\mathbb{R}^{2n}$-valued functions over $\mathbb{R}$, and appeal to the obvious generalization of the following lemma to $\mathbb{R}^{2n}$-valued functions. From the asymptotic conditions on $f$ and the fact that $[\Pi_c, T']$ is supported on a bounded interval in $\mathbb{R}$, we see that the left hand side of the above equation is indeed an example of the operators in assertion 2 of the next lemma. Proposition 5.3.2 is then proved. $\square$

**5.3.4 Lemma** *Let $\varsigma$ be a function of 1-variable ($f$) defined by (2.9). Let $\mathcal{L}^2_{k:\epsilon}$ be the completion of the following weighted norms for $\mathbb{R}$-valued functions (Note that $\mathcal{L}^2_{0:\epsilon} = L_{2:\epsilon}$): Let $u(f) \in C_0^\infty(\mathbb{R})$,*

$$\|u\|_{\mathcal{L}^2_{k:\epsilon}} := \sum_{i=0}^{k} \|\varsigma^{i+\epsilon}\nabla^i f\|_2.$$

*Then for $\epsilon > 1$,*

1. *The operator $\frac{d}{df}$ is Fredholm between $\mathcal{L}^2_{k+1:\epsilon-1}$ and $\mathcal{L}^2_{k:\epsilon}$ of index $-1$. It has trivial kernel, and a 1-dimensional cokernel spanned by the function $\varsigma^{-2\epsilon}$.*

2. *Let $\nu(f)$ be a $C_k$ function such that $\sum_{i=0}^{k} |\nabla^{(i)}\nu|(f) \leq C\varsigma^3$. By perturbation, operators of the form $\frac{d}{df} + \nu(f)$ are also Fredholm between $\mathcal{L}^2_{k+1:\epsilon-1}$ and $\mathcal{L}^2_{k:\epsilon}$ with index $-1$.*

3. *The formal $L^2_\epsilon$-adjoint of $\frac{d}{df}$ is also Fredholm between $\mathcal{L}^2_{k+1:\epsilon-1}$ and $\mathcal{L}^2_{k:\epsilon}$, and has index $1$.*

4. *$\frac{d}{df}$ may be extended to be a Fredholm operator, $\hat{d}$, between $\hat{\mathcal{L}}^2_{k+1:\epsilon-1}$ and $\mathcal{L}^2_{k:\epsilon}$, where*

$$\hat{\mathcal{L}}^2_{k+1:\epsilon-1} := \{\lambda_d c_+ + (1 - \lambda_d)c_- + e : c_+, c_- \in \mathbb{R}; e \in \mathcal{L}^2_{k:\epsilon}\} \simeq \mathbb{R}^2 \times \mathcal{L}^2_{k+1:\epsilon-1},$$

*with the product Banach structure. This extended operator has index $1$: it has a 1-dimensional kernel consisting of the constant functions, and null cokernel. In fact, it has a bounded right inverse.*

*Proof.* 1. We will deal with the $k = 1$ case only; the $k > 1$ case follows from the $k = 1$ case via elliptic regularity (cf. §5.4).

We first change the variable from $f$ to $s := \int_{-\infty}^{f} \varsigma(u)du$. Let $\tilde{e}(s) := \varsigma(f(s))$. Note that by the definition of $\varsigma$ and $s$, $\tilde{e}(s) = Ce^{|s|/R}$ for all $s$ with large $|s|$. ($R$ here is as in (2.2.5).) Let $\tilde{\mathcal{L}}_{k:\epsilon}$ be the conventional exponential weighted Sobolev spaced defined by the norm:

$$\|f\|_{\tilde{\mathcal{L}}_{k:\epsilon}} := \sum_{i=0}^{k} \|\tilde{e}^{\epsilon} \nabla^i f\|_2.$$

We see that after the change of variables, $\mathcal{L}_{1:\epsilon-1}^2$ becomes $\tilde{\mathcal{L}}_{1:\epsilon-1/2}^2$; $\mathcal{L}_{0:\epsilon}^2$ becomes $\tilde{\mathcal{L}}_{0:\epsilon+1/2}$, and the operator $\frac{d}{df}$ becomes $\tilde{e}^{-1}\frac{d}{ds}$. In other words, we need only examine the operator $\frac{d}{ds}$ between $\tilde{\mathcal{L}}_{1:\epsilon-1/2}$ and $\tilde{\mathcal{L}}_{0:\epsilon-1/2}$, which is well-known to be Fredholm. For example, by identifying $\tilde{\mathcal{L}}_{k:\epsilon-1/2}$ with $L_k^2$ via multiplication by $\tilde{e}^{\epsilon-1/2}$, this can be seen by translating to the problem of $\frac{d}{ds} - (\epsilon - 1/2)\tilde{e}^{-1}\frac{d\tilde{e}}{ds}$ between $L_1^2$ and $L^2$, where $\tilde{e}^{-1}\frac{d\tilde{e}}{ds} = \pm R^{-1}$ for all $s$ with large enough $|s|$.

The statement about the kernel and cokernel is easy to check. Especially, it is easy to see that $\frac{d}{df}$ has a bounded inverse on the $\mathcal{L}_{0:\epsilon}^2$-orthogonal complement of the cokernel when $R$ is large enough. Because of the Fredholmness, it suffices to verify that $\|u\|_{\mathcal{L}_{1:\epsilon-1}^2} \leq C\|\frac{du}{df}\|_{\mathcal{L}_{0:\epsilon}^2}$ for all $u \in C_0^\infty$. Now write $u = \varsigma^{-\epsilon}g$; The above follows from the fact that in this case $\|g' - \epsilon\varsigma'\varsigma^{-1}g\|_2^2 \geq \|g'\|_2^2 + \epsilon^2\|\varsigma^{-1}g\|_2^2$ modulo small terms which may be absorbed to the first term.

*2.* is obvious. *3.* is similar to *1*.

*4.* The Fredholmness and the index calculation is obvious from the definition. For the assertion about the kernel/cokernel, note that $\hat{d}\lambda$ has nontrivial component in the cokernel of $\frac{d}{df}$. The boundedness of the right inverse follows from this and the last paragraph of the proof of *1.* above. □

Assertion 4 above together with Lemma 5.3.3 imply:

**5.3.5 Corollary** *(a) When $f = x_3$ and $c$ is an admissible Seiberg-Witten solution described in §3.2, then $\mathfrak{D}_c$ has a uniformly (in $c$) bounded right inverse, and its kernel is $\mathbb{C}^n$.*

*(b) More generally, if in 5.1.4, $R$ is large enough so that $\tilde{f}$ is close enough to $x_3$, and if $c(z_+, \tilde{f}) = (v(z_+), 0)$ with respect to (5.1), where $v$ is a vortex solution on $\mathbb{C}$, then the same statements hold.*

The uniform boundedness of the right inverse in part (a) above follows from the boundedness of $\hat{d}$ in the previous lemma, the $c$-independence of the constant $C$ in (5.23), and the uniform boundedness of $\mathfrak{D}_c^{'\dagger}$. The only difference of (b) from (a) is that here instead of $\hat{d}$, we have $\hat{d} + |\nabla f|^{-1}\Pi_c R'(1 + v')$. This is still bounded uniformly, because by assumption $\|R'\|_{C_{k-2}} \leq C\varepsilon$ is small. We remark that the uniform boundedness of the right inverse obtained in this Corollary will be important for the proof of the gluing theorem in §6.2.

CASE (II). *Let $c$ be a configuration on $M$ whose restriction to $M \backslash M_R$ agrees with a configuration $c_1$ on $\mathbb{R}^3$ of the form in case (i) above.*

Let $\chi_2 := \chi_R$; $\chi_1 := 1 - \chi_2$. Also let $\varphi_2 := \chi_{2R}$; $\varphi_1 := 1 - \chi_{R/2}$. Note that $\varphi_i$, $i = 1, 2$, are smooth cutoff functions of value 1 on the supports of $\chi_i$ respectively.

**5.3.6 The uniform boundedness of $\mathcal{D}'_c$ in case (ii)** is due to the following inequality:

$$\begin{aligned}
\|\mathcal{D}'_c \zeta\|_{2,1:\epsilon} &\leq \|\mathcal{D}'_c(\chi_1 \zeta)\|_{2,1:\epsilon} + \|\mathcal{D}'_c(\chi_2 \zeta)\|_{2,1:\epsilon} \\
&\leq C_1 \|\chi_1 \zeta\|_{K_{2,c_1}} + C_2 \|\chi_2 \zeta\|_{K_{2,c'}} \\
&\leq C \|\zeta\|_{K_{2,c}},
\end{aligned} \quad (5.37)$$

where $C$ is a constant depending on $c_1$, $c'$, but not on $\zeta$.

**5.3.7 To show Fredholmness of $\mathcal{D}'_c$ for $l = 1$ in case (ii),** we may construct a parametrix $Q_c$ of $\mathcal{D}'_c$ by defining

$$Q_c := \varphi_1 Q_1 \chi_1 + \varphi_2 Q_2 \chi_2, \quad (5.38)$$

where $Q_1$ is a parametrix for $\mathcal{D}'_{c_1}$ on $\mathbb{R}^3$ constructed in case (i) and $Q_2$ is a parametrix for the compact piece, which may for example be constructed as follows: We may always extend $M_{2R}$ to a compact closed manifold, say $X$. The spinor bundle $S\big|_{M_{2R}}$ extends trivially over $X$. We may also extend the configuration $c\big|_{M_{2R}}$ smoothly over $X$ (we call it $c_2$.) Because of the compactness of $X$ and the ellipticity of $\mathcal{D}'_{c_2}$, there is a parametrix to $\mathcal{D}'_{c_2}$ on $X$, which we call $Q_2$.

One may check that $Q_c$ is a parametrix for $\mathcal{D}'_c$ on $M$ by noting that the difference of $\mathcal{D}'_c Q_c$ or $Q_c \mathcal{D}'_c$ from the identity are sums of terms which are products of a compact operator with a bounded operator, which are still compact. (Note that the commutators of $\mathcal{D}'_c$ and the cutoff functions have compact supports, and therefore are compact.) For example, a direct computation shows that

$$\mathcal{D}'_c Q_c = \text{Id} + \mathcal{R}_c, \quad (5.39)$$

where

$$\begin{aligned}
\mathcal{R}_c :=& \sigma_{\mathcal{D}'}(d\phi_1) Q_1 \chi_1 + \sigma_{\mathcal{D}'}(d\phi_2) Q_2 \chi_2 + \phi_1 R_1 \chi_1 + \phi_2 R_2 \chi_2 \\
&+ \phi_1 (\mathcal{D}'_c - \mathcal{D}'_{c_1}) Q_1 \chi_1 + \phi_2 (\mathcal{D}'_c - \mathcal{D}'_{c_2}) Q_2 \chi_2,
\end{aligned} \quad (5.40)$$

where $R_1$, $R_2$ are compact operators defined by the formulae

$$\mathcal{D}'_{c_{1,2}} Q_{1,2} = 1 + R_{1,2}. \quad (5.41)$$

Note the term $\phi_1(\mathcal{D}'_c - \mathcal{D}'_{c_1}) Q_1 \chi_1$ is compact because $Q_1 \chi_1$ is bounded, and $\phi_1(\mathcal{D}'_c - \mathcal{D}'_{c_1})$ is supported on the region where $R \leq |x| \leq 2R$. This proves the Fredholmness of $\mathcal{D}'_c$ between $K_{1,c}$ and $L_{:\epsilon}(\mathcal{W}_1)$ in case (ii).

CASE (III): THE GENERAL CASE. The configuration in case (ii) is not a general configuration, but it is close because for any configuration $c' \in \mathcal{C}(M)$, there is a configuration $c$ in the form of case (ii), such that

$$c' = c + q, \text{ where } q \in Y.$$

In fact, we can choose $c$ such that $q = (1 - \chi_{R/2})q$. Also note that the two parameters $d$ and $R$ in the definition of $c$ are mutually independent. We shall later take $R \to \infty$, but leave $d$ fixed (depending only on $a_1, a_2$).

**5.3.8 The boundedness of $\mathcal{D}'_{c'}$ in the general case** follows from the inequality

$$\|\mathcal{D}'_{c'}\zeta - \mathcal{D}'_c\zeta\|_{2,1:\epsilon} \leq \|q\zeta\|_{2,1:\epsilon}$$
$$\leq C\|q\|_{3,1:\epsilon}\|\nabla(\varsigma\zeta)\|_{2,1} \leq C'\|q\|_Y \|\zeta\|_{K_{2,c}}. \tag{5.42}$$

**5.3.9 To show the Fredholmness of $\mathcal{D}'_c$ for $l = 1$ in the general case,** note that similarly to (5.42), we have

$$\|\mathcal{D}'_{c'}\zeta - \mathcal{D}'_c\zeta\|_{2:\epsilon} \leq C\|q\|_{Y_1}\|\zeta\|_{K_{1,c}}.$$

Here $\|q\|_{Y_1} = \|(1 - \chi_{R/2})q\|_{Y_1} \leq \varepsilon(R)$, where $\varepsilon(R)$ is an $R$-dependent positive number which may be made arbitrarily small by taking $R$ sufficiently large. (In particular, when $c'$ is a Seiberg-Witten solution or a reference configuration, $q \leq o^3$, and therefore in this case $\|(1 - \chi_{R/2})q\zeta\|_{2:\epsilon} \leq CR^{-3}\|\zeta\|_{K_{1,c}}$.)

We may then appeal to [16], IV.5.22, which implies that $\mathcal{D}'_{c'}$ is Fredholm.

The proof of Theorem 5.1.10 assertion (0) and assertion (1) in the $l = 1$ case is now complete.

## 5.4 Fredholmness for $\mathcal{D}'^\dagger_c$ and generalization to $l > 1$: the weighted case

This subsection contains the second half of the proof of Theorem 5.1.10. We shall first show that $\mathcal{D}'^\dagger_c$ is Fredholm for $l = 1$, then show that by elliptic regularity, this together with the Fredholmness of $\mathcal{D}'_c$ for $l = 1$ proved in the last subsection imply the Fredholmness of both $\mathcal{D}'_c$ and $\mathcal{D}'^\dagger_c$ for all $l$.

**5.4.1 The Fredholmness of $\mathcal{D}'^\dagger_c$ for $l = 1$.** We first review some terminology and facts.

Let $D$ be a linear differential operator mapping $\Gamma(V)$ to itself for some bundle $V$. Let $E$ and $F$ be two Hilbert spaces containing $C_0^\infty(V)$ (typically both taken to be $L^2(V)$), such that $E$ contains the domain of $D$ and $F$ contains the range of $D$. There are many possible choices for the domain of $D$; the following are two common choices: $\text{Dom}_{\min}(D)$, or equivalently the *minimal extension*, is the completion of $C_0^\infty(V)$ in $E$ with respect to the norm:

$$\|\zeta\|_D^2 := \|D\zeta\|_F^2 + \|\zeta\|_E^2.$$

The other choice, $\text{Dom}_{\max}(D)$, or equivalently the *maximal extension*, consists of all elements $\zeta \in E$ for which $D\zeta$ makes sense as a distribution in $F$. It follows almost immediately (cf. [16]) from the definitions that if we let $D : \text{Dom}_{\min}(D) \to F$, and $D^{\text{t}} : \text{Dom}_{\max}(D^\dagger) \to E$ (here $D^\dagger$ is the formal adjoint of $D$; $D^{\text{t}}$ is *the adjoint of $D$ in the sense of [16]*), then $D^{\text{t}}$ is Fredholm if $D$ is, and $\text{Ind}(D^{\text{t}}) = -\text{Ind}(D)$.

For our application here, taking $D = \mathcal{D}'_c$, $E = K_{0,c}$ and $F = L^2_\epsilon(\mathcal{W}_1)$, the Fredholmness of $\mathcal{D}_c^{'\dagger}$ for $l = 1$ thus follows from the following lemma.

**Lemma** *(a) $K_{1,c}$ coincides with the completion of $C_0^\infty$ with respect to the norm $\|\mathcal{D}'_c\zeta\|_{2:\epsilon} + \|\zeta\|_{K_{0,c}}$ (the minimal extension).*
*(b) The subspace in $L^2_\epsilon(\mathcal{W}_1)$ where $\mathcal{D}_c^{'\dagger}$ is defined as a distribution in $K_{0,c}$ (the maximal extension) coincides with $R_{1,c}$.*

*Proof.* The arguments are similar to the self-adjointness proof in §5.2. (cf. also [20] section 3, and [16] p.167.) Similar to §5.3 above, we may reduce the general case to case (i) on $\mathbb{R}^3$, since the excision and perturbation arguments in steps (ii)&(iii) can be easily modified.

For the rest of this proof, let $c$ be a configuration on $\mathbb{R}^3$ as in §5.3 case (i).

To show (a), we need to establish the inequality

$$\|\zeta\|_{K_{1,c}} \leq C\Big(\|\mathcal{D}'_c\zeta\|_{2:\epsilon} + \|\zeta\|_{K_{0,c}}\Big). \tag{5.43}$$

(The other direction is true by the continuity of $\mathcal{D}'_c$ between $K_{1,c}$ and $L^2_\epsilon$.) Similar to §5.3 case (i), decompose

$$\zeta = v + w,$$

where $v \in L^2_\epsilon(\mathcal{W}_1, \mathcal{K}_c)$ and $w \in L^2_{\epsilon-1}(\mathcal{K}_c)$. (5.43) then follows from a combination of the following two inequalities:

$$\|v\|_{2,1:\epsilon} \leq C_1(\|v\|_{2:\epsilon} + \|\mathcal{D}'_c v\|_{2:\epsilon}),$$

which is obtained by a straightforward computation similar to the arguments for (5.16), and

$$\|w\|_{K_{1,c}} \leq C_2\Big(\|T'w\|_{2:\epsilon} + \|w\|_{2:\epsilon-1}\Big)$$

(which follows easily from the definition of $K_{1,c}$), plus the fact that the cross term

$$|\langle \mathcal{D}'_c v, \mathcal{D}'_c w\rangle_{2:\epsilon}|$$
$$= \Big|\langle \mathcal{D}'_c v, [T', \Pi_c]w\rangle_{2:\epsilon} + \langle [T', \Pi_c]v, \mathcal{D}'_c w\rangle_{2:\epsilon}$$
$$+ \langle \mathcal{D}'_c v, R'w\rangle_{2:\epsilon} + \langle R'v, i(T' + N')w\rangle_{2:\epsilon}\Big|$$
$$\leq \varepsilon(\|\mathcal{D}'_c v\|^2_{2:\epsilon} + \|\mathcal{D}'_c w\|^2_{2:\epsilon}) + C(\|v\|^2_{2:\epsilon} + \|w\|^2_{2:\epsilon-1})$$

for some small number $\varepsilon$.

To show (b), let $\mathcal{D}'^t_c$ denote the adjoint of $\mathcal{D}'_c$ in the sense of [16]. Then $\mathcal{D}'^t_c$ is Fredholm since $\mathcal{D}'_c$ is, and $\text{Ind}(\mathcal{D}'^t_c) = -\text{Ind}(\mathcal{D}'_c)$. We want to show that $\text{Dom}(\mathcal{D}'^t_c) \subset R_{1,c}$. (The other direction follows from the continuity of $\mathcal{D}'^\dagger_c$). Let $\xi \in \text{Dom}(\mathcal{D}'^t_c) \subset L^2_\epsilon(\mathcal{W}_1)$. Approximate $\xi$ by $\xi_i \in C^\infty_0 \to \xi$, $i = 1, 2, \ldots$. By the Fredholmness of $\mathcal{D}'^t_c$, $\|\mathcal{D}'^t_c \xi_i\|_{K_{0,c}} \to \|\mathcal{D}'^t_c \xi\|_{K_{0,c}}$. The left hand side is easier to estimate, since for $\xi_i \in C^\infty_0$, $\mathcal{D}'^t_c \xi_i = \mathcal{D}'^\dagger_c \xi_i$, and integration by parts is applicable. So without loss of generality, we may assume that $\xi \in \mathbb{C}^\infty_0$. Again decompose

$$\xi = h + t,$$

where $h \in L^2_\epsilon(\mathcal{W}_1; \mathcal{K}_c)$ and $t \in L^2_\epsilon(\mathcal{K}_c)$. The relevant estimate,

$$\|\xi\|_{R_{1,c}} \leq C(\|\mathcal{D}'^\dagger \xi\|_{K_{0,c}} + \|\xi\|_{2:\epsilon}),$$

then follows from combining the inequalities (5.44), (5.45) below.

$$\|\nabla h\|_{2:\epsilon} \leq C\Big(\|\mathcal{D}'^\dagger_c h\|_{2:\epsilon} + \|h\|_{2:\epsilon}\Big) \tag{5.44}$$

is established in (5.28). Also, note that

$$T'^\dagger t = -(|\varsigma|^2|\nabla f|\partial_f t + (2\epsilon|\varsigma| + o^1)t)$$

where $|f| \gg R$; so

$$\|\partial_f t\|_{2:1+\epsilon} \leq C'\Big\|T'^\dagger t - (2\epsilon|\varsigma| + o^1)t\Big\|_{2:\epsilon-1}$$
$$\leq C''\Big(\|T'^\dagger t\|_{2:\epsilon-1} + \|t\|_{2:\epsilon}\Big). \tag{5.45}$$

We have now completed the proof that $\mathcal{D}'^\dagger_c$ is Fredholm between $R_{1,c}$ and $K_{0,c}$. □

**5.4.2 Generalizing to higher $l$.** Similar to §5.2, the next "elliptic regularity"-type of lemma shows that if $\mathcal{D}'_c$ or $\mathcal{D}'^\dagger_c$ are Fredholm for $l = k$, then they are also Fredholm for $l = k+1$. Together with the $l = 1$ case already proved, this proves the claims of Theorem 5.1.10 (1), (2) for all $l \in \mathbb{Z}^+$. Furthermore, the kernel, cokernel, and hence the index are independent of $l$. And since by 5.4.1 $\text{Ind}\,\mathcal{D}'^\dagger_c = -\text{Ind}\,\mathcal{D}'_c$ in the $l = 1$ case, this holds for all $l$.

**Lemma** *For any $l \in \mathbb{Z}^+$,*
*(a) if $\zeta \in K_{l,c}$, $\mathcal{D}'_c\zeta \in L^2_{l:\epsilon}(\mathcal{W}_1)$, then $\zeta \in K_{l+1,c}$;*
*(b) if $\xi \in R_{l,c}$; $\mathcal{D}'^\dagger_c\xi \in K_{l,c}$, then $\xi \in R_{l+1,c}$.*

*Proof.* The lemma follows from the following inequalities:

$$\|\nabla^{l+1}\zeta\|^2_{2:\epsilon} \leq C(\|\nabla^l \mathcal{D}'_c\zeta\|^2_{2:\epsilon} + \|\zeta\|^2_{K_{l,c}});$$
$$\|\nabla^{l+1}\xi\|^2_{2:\epsilon} + \|\partial^{l+1}_f(\Pi_c(1-\lambda)\xi)\|^2_{2:2+\epsilon} \leq C'(\|\nabla^l \mathcal{D}'^\dagger_c\xi\|^2_{2:\epsilon} + \|\xi\|^2_{R_{l,c}}).$$

These formulae can be verified by direct computation similar to the proof of Lemma 5.4.1. □

End of the proof of Theorem 5.1.10.

## 5.5 Boundedness and Fredholmness of $\mathcal{D}_c$

Theorem 5.1.10 has the following extension as a consequence:

**5.5.1 Theorem** *Let $t > 0$, $c \in \mathcal{C}$. Then $\mathcal{D}_c$ is bounded (uniformly in c) and Fredholm between $\hat{K}_1$ and $\hat{L}^2_{1:\epsilon}(\mathcal{W}_1)$, with index $4n + \mathrm{Ind}(\mathcal{D}'_c)$.*

*Proof.* Recall the decomposition of $\mathcal{D}_c$ from (5.13) (and adopt the notation therein).

To show the boundedness of $\mathcal{D}_c$, we only need to show that $\mathcal{X}_c$ is bounded. By inspection, the $L^2_{1:\epsilon}$ norm of the first and second terms in (5.14) are both bounded by $\frac{C}{R}(\|v_1\|_{L^2_1(\mathbb{C},T\mathbb{C}\oplus\mathbb{C})} + \|v_2\|_{L^2_1(\mathbb{C},T\mathbb{C}\oplus\mathbb{C})})$, and the $L^2_{1:\epsilon}$ norms of third and fourth terms are bounded by

$$C'(d^{1/2}\|a_1 - a_2\|_{L^2_{2,R}(\mathbb{C},T\mathbb{C}\oplus\mathbb{C})} + \|q\|_{2,2:\epsilon,R})\|v_i\|_{L^2_1(\mathbb{C},T\mathbb{C}\oplus\mathbb{C})}$$

for $i = 1, 2$ respectively. (The notation $\|\cdot\|_{\cdot,R}$ means the norm is taken outside of the disc or ball of radius $R$ at the center.)

As for the Fredholmness, since we know that $\Theta_{a_i}$ ($i = 1, 2$) and $\mathcal{D}'_c$ are Fredholm, from (5.13) we may conclude that $\mathcal{D}_c$ is Fredholm of index $4n + \mathrm{Ind}(\mathcal{D}'_c)$ (because $\mathrm{Ind}(\Theta_{a_i}) = 2n$) if $(0, 0, \mathcal{X}_c)$ can be regarded as a perturbation. This holds because from the form of (5.14), the operator norm of $\mathcal{X}_c$ may be made arbitrarily small if $R$ (the parameter in the cutoff function $\chi = \chi_R$) is taken to be large enough.

This finishes the proof of Theorem 5.5.1. □

# 6 The moduli space

Armed with the results obtained in previous sections, in §6.1 we apply the standard arguments to deduce the usual smoothness and invariance properties of the moduli spaces. In §6.2 we discuss the compactness of the moduli spaces. In the $t > 0$ case, the moduli space has non-compact components, and we also need an explicit description of the "ends" of the moduli space. This is provided in §6.3, 6.4. Using these results one may define some gauge-theoretic invariants of 3-manifolds. An outline of the definition, and a preliminary example of them was in [23] section 7.

Throughout this section we let $l = 2$. Higher $l$ versions of the results are straightforward (yet unnecessary, in view of lemma 5.1.2) generalizations.

We will be mainly dealing with the $t > 0$ case. The $t = 0$ case is similar and relatively simple.

### 6.1 Smoothness and invariance

Recall from §5.1 the definition of the moduli space $\mathcal{M}$ as $\mathcal{S}^{-1}(0)/\mathcal{G} \subset \mathcal{Q}$, where $\mathcal{S} : \mathcal{C} \to \mathcal{V}$ is given by (5.4). Note that for $\mathcal{S}$ to be zero, $\omega$ has to be closed.

It is well-known that (cf. [2] Proposition 4.2.23) a local model of $\mathcal{M}$ at $[c] \in \mathcal{M}$ is given by $f_c^{-1}(0)/\Gamma_c$, where $f_c$ is a map from $\text{Ker}\,\mathcal{D}_c = H_1$ to $\text{Coker}\,\mathcal{D}_c = H_2$ ($H_i$ are the $i$-th cohomology of the elliptic complex (5.6)) which vanishes at first order, and $\Gamma_c$ is the stabilizer of the gauge action at $c$. Thus $\mathcal{M}$ is smooth at $c$ when $\text{Coker}\,\mathcal{D}_c = \emptyset$ and $\Gamma_c = \{1\}$, and its dimension is $\text{Ind}(\mathcal{D}_c)$. Singularities occur in two occasions. The first kind of singularities has a cone-like neighborhood and comes from reducible configurations ($\Gamma_c \neq \{1\}$). It is in a sense more fundamental because it is already present at the level of $\mathcal{Q}$ and is usually difficult to get rid of by perturbations. The second kind of singularities appear at $c \in \mathcal{S}^{-1}(0)$, when $\mathcal{D}_c$ has a non-trivial cokernel (since $f_c$ vanishes at first order). As usual, the second kind of singularities can be killed by putting in a generic perturbation 2-form $w$. (Recall the definition of $w$ and $\omega$ from (2.4). We let $w \in L^2_{1:\epsilon}(\bigwedge^2 T^*M)$.

In our case, the first kind of singularities disappears also too when $w$ is introduced. More explicitly, for a reducible Seiberg-Witten solution $(A, 0)$, we have $F_A + i\omega = 0$. Integrating over a surface $H$ described in the end of 3.3.8,

$$\text{vortex number} = \int_H \frac{i}{2\pi} F_A = -\int_H \frac{\omega}{2\pi}$$

if $(A, 0)$ solves the Seiberg-Witten equation (2.2). Thus for a generic $\omega$ with $-\int_H \frac{\omega}{2\pi} \neq$ vortex number, reducible solutions do not occur, and we are happily in the situation where the moduli space is a smooth manifold. In particular, for an admissible $\theta$, $\theta$ is not even in $L^2_{\epsilon}(\bigwedge^2 T^*M)$; thus for the $t > 0$ case, reducible Seiberg-Witten solutions never occur. If $t = 0$ and $b_1(M) > 0$, we can choose a small $w = \omega$ such that $\int_Z \omega \neq 0$ for some compact cycle $Z$ to get rid of singularities.

**6.1.1 Remark** In the previous sections, we assumed $w = 0$. When $w$ is introduced, the corresponding Seiberg-Witten solutions might no longer satisfy the pointwise estimates in Proposition 3.3.3. However, they still lie in the configuration space $\mathcal{C}$ defined in section 4 if $w$ is small.

In this section, we fix $\theta$, and denote the moduli space corresponding to the admissible metric $g$ and perturbation $\omega = -t\frac{\theta}{2} + w$ as $\mathcal{M}_{g,w,t}$. We show below that these moduli spaces are cobordant.

Let Met be the *Banach manifold of admissible metrics* whose tangent space $\mathcal{H}$ is $L^2_9(\text{Sym}^2 T^*M_\Re)$ with respect to a smooth metric. (Though the $L^2_9$-*norm* depend on $g$, by lemma 2.2.12 the $L^2_9$-*space* does not.) The topology of Met is defined using the $L^2_9$-norm. *The space of (minor) perturbation 2-forms*, $\mathcal{B}$, is $L^2_{1:\epsilon}(\bigwedge^2 T^*M)$, $\epsilon \in (1, 3/2)$. (For higher $l$ cases, Met is modeled on $L^2_k(\text{Sym}^2 T^*M_\Re)$, $k \geq l+6$, and $\mathcal{B}$ is $L^2_{l-1:\epsilon}(\bigwedge^2 T^*M)$.)

Define the map $\tilde{\mathcal{S}} : \text{Met} \times \mathcal{B} \times \mathcal{C} \to \mathcal{V}$ by $\tilde{\mathcal{S}}(g,w,c) := \mathcal{S}_{g,w}(c)$, where $\mathcal{S}_{g,w}$ is given by (5.4) with the metric $g$, and with $\omega$ given by (2.4). (Note that by routine estimates, $\mathcal{C}, \mathcal{V}$ do not depend on the metric $g$.) The quotient $\tilde{\mathcal{M}}_t := \tilde{\mathcal{S}}^{-1}(0)/\mathcal{G}$ is called the *parameterized moduli space*. The next theorem shows that $\tilde{\mathcal{M}}_t$ is smooth and thus forms a cobordism between different moduli spaces.

**6.1.2 Theorem** *For any $t$, the map $\tilde{\mathcal{S}}$ is smooth and zero is a regular value.*

*Therefore for the $t > 0$ case, $\tilde{\mathcal{M}}_t = \coprod_{n=0}^{\infty} \tilde{\mathcal{M}}_t^n$ where each $\tilde{\mathcal{M}}_t^n$ is a Banach manifold, and the projection*

$$\varpi : \tilde{\mathcal{M}}_t^n \subset \mathcal{Q} \times \mathcal{H} \times \mathcal{B} \longrightarrow \mathcal{H} \times \mathcal{B}$$

*is Fredholm of index equal to* $\text{Ind}(\mathcal{D}_c)$. *Here $n$ is the vortex number of the elements in each component. As a consequence, for generic $(g, w) \in \mathcal{H} \times \mathcal{B}$, $\mathcal{M}_{g,w,t}^n = \varpi^{-1}(g, w)$ is a smooth manifold.*

*Similarly in the $t = 0$ case, for an MEE with $b_1 > 0$ and for generic $(g, w)$, the moduli space $\mathcal{M}_{g,w,0}$ is a 0-dimensional smooth manifold.*

*Proof.* First note as before that for the $t > 0$ case, we don't have to worry about the reducible solutions. For the $t = 0$ case, we need the extra condition $b_1 > 0$ on $M$ to guarantee the nonexistence of reducible solutions for generic $(g, w)$.

Below we prove the theorem for the $t > 0$ case. The $t = 0$ case is similar.

The smoothness of $\tilde{\mathcal{S}}$ follows from the boundedness of $\mathcal{D}_c$ (Theorem 5.1.10 (0)) (hence $D_c = d\mathcal{S}$ is bounded). For the regularity at 0, it suffices to prove that the restriction, $\mathcal{S}'$, of $\tilde{\mathcal{S}}$ to $\mathcal{B} \times \mathcal{C} \xrightarrow{\mathcal{S}'} \mathcal{V}$ has zero as a regular value[5]. Let $c$ be a solution to the Seiberg-Witten equations (2.2). Note that since $\mathcal{D}_c$ has closed range (being Fredholm), $\text{Im}(D_c) \subset \text{Ker}(d_c^*) \subset \mathcal{V}$ is closed. On the other hand, $d\mathcal{S}'(v, e) = D_c e + (-iv, 0) \in \Gamma(T^*M \oplus S)$; thus its image is the sum of two closed subspaces). Therefore it suffices to show that $\text{Coker}(d_c \mathcal{S}') = \emptyset$.

---
[5]Let $Y$ be replaced by $K_{2,c} \cap \Gamma(\mathcal{W}_0)$ (identifying $\Gamma(\mathcal{W}_0)$ with its natural embedding in $\Gamma(\mathcal{W}_1)$). cf. §5.1, footnote 3.

Since $d\mathcal{S}'$ surjects to the $L^2_1(\mathbb{C}, iN \oplus \mathbb{C} \otimes T\mathbb{C})^{\oplus 2}$ component of $\mathcal{V}$ (because the moduli space of vortex solutions on $\mathbb{C}$ are smooth), we need only to consider the $L^2_{1:\epsilon}(\mathcal{W}_0)$ component of $\mathcal{V}$. However, a by-now standard argument (cf. e.g. [29]) shows that there is no $(F, \chi) \in L^2_{1:\epsilon}(\mathcal{W}_0)$ in the cokernel of $d\mathcal{S}'$. Thus $d\mathcal{S}'$ (and hence $d\tilde{\mathcal{S}}$) is surjective at zero and therefore $\tilde{\mathcal{S}}^{-1}(0)$ is a Banach manifold. By Theorem 4.3.5, $\tilde{\mathcal{M}}_t = \tilde{\mathcal{S}}^{-1}(0)/\mathcal{G}$ is a Banach submanifold of $\mathcal{Q}$.

As a consequence, $\varpi : \tilde{\mathcal{M}}_t \to \mathcal{H} \times \mathcal{B}$ is Fredholm of index $\mathrm{Ind}(\mathcal{D}_c)$ by a standard argument (cf. e.g. [9], p.60). It follows then from the Sard-Smale theorem that $\mathcal{M}_{g,w,t} = \varpi^{-1}(g, w)$ is a smooth manifold for generic $(g, w)$. $\square$

As a simplest example of $\mathcal{M}$,

**6.1.3 Corollary** *The moduli space for Seiberg-Witten solutions on $\mathbb{R}^3$ with perturbation $t * dx_3$ ($t > 0$) is smooth; in fact, it is diffeomorphic to $\coprod_{n \geq 0} \mathrm{Sym}^n \mathbb{C}$.*

*Proof.* This is obvious from the work in §3.2. In fact, the embedding of the moduli in the configuration space $\mathcal{C}$ in the $*$-gauge is given simply by $(A, (\alpha, 0))$, where $(A, \alpha)$ is a vortex solution on $\mathbb{C}$ in the **(vor)** gauge. The diffeomorphism to $\coprod_{n \geq 0} \mathrm{Sym}^n \mathbb{C}$ is then given by taking the centers of the vortex solution $(A, \alpha)$. The smoothness follows from Corollary 5.3.5. $\square$

We say that $(g, w)$ is a *good pair* if the corresponding moduli space $\mathcal{M}_{g,w,t}$ is smooth.

The following corollary is useful for showing that the Seiberg-Witten invariants of MEE's are independent of metric and perturbations.

**6.1.4 Corollary** *(a) Fix the parameter $t > 0$. The moduli spaces $\mathcal{M}_{g_1,w_1,t}$, $\mathcal{M}_{g_2,w_2,t}$ corresponding to any two good pairs of $(g_1, w_1)$, $(g_2, w_2)$ are cobordant in the quotient space $\mathcal{Q}$.*

*When $t = 0$, the same holds with the additional requirement that $b_1(M) > 1$.*

*(b) Now fix a good pair $(g, w)$ and let $t$ vary. The moduli spaces $\mathcal{M}_{g,w,t}$ patch up to form a (non-compact) cobordism:*

$$\hat{\mathcal{M}}_{g,w} \xrightarrow{\pi} \mathbb{R}^+,$$

*with $\pi^{-1}(t) = \mathcal{M}_{g,w,t}$, $t \in (0, \infty)$. Furthermore, when the vortex number $n = 0$ and $\epsilon = 0$ (note that $\epsilon = 0$ is allowed only when $n = 0$), elements in $\mathcal{M}^0_{g,w,t}$ converge in $L^2_{1,loc}$ to points in $\mathcal{M}_{g,w,0}$ as $t$ goes to zero.*

**6.1.5 Remark** Note that $\mathcal{Q}_0$ and $\mathcal{Q}_t$, $t > 0$ have different topologies so the convergence $\mathcal{M}_{g,w,t} \to \mathcal{M}_{g,w,0}$ can only happen at the $L^2_{1,loc}$ level. Also, $\mathcal{M}_{g,w,t}$, $t > 0$ do not converge to the whole $\mathcal{M}_{g,w,0}$. For example, when $M$ is the connected sum of $\mathbb{R}^3$ and a Seifert fibered manifold with $b^1 = 0$, we expect (by the correspondence between Seiberg-Witten solutions and sets of gradient flow lines explained in §1.1 (A)) that the limit of $\mathcal{M}_{g,w,t}$, $t \to 0$ consists only of flat connections on $M$, while there might be other $t = 0$ Seiberg-Witten solutions in view of the computation in [32].

*Proof.* Part (a) is standard. The proof of part (b) follows the same strategy, but some preliminary work is required: since the moduli spaces $\mathcal{M}_{g,w,t}$ for different $t$'s lie in different spaces $\mathcal{Q}_t$, we need to construct a bigger ambient space in which $\hat{\mathcal{M}}_{g,w}$ lies. This is done by a "fibration argument" similar to that in §4.3. Recall the construction of the configuration space $\mathcal{C}_t$ as $\bigcup_n Y_t \times \mathcal{A}_t^n \times \mathcal{A}_t^n$. We showed in Lemma 4.3.3 that the space $Y_t$ does not depend on $t$. In addition, though $\mathcal{A}_t^n$ does depend on $t$, $\mathcal{A}_{t_1}^n$ and $\mathcal{A}_{t_2}^n$ are isomorphic for any $t_1, t_2 \in \mathbb{R}^+$ by the isomorphism $(a, \alpha)(z) \mapsto (\sqrt{t_2/t_1}\, a, \sqrt{t_2/t_1}\, \alpha)(\sqrt{t_1/t_2}\, z)$. Thus the union

$$\hat{\mathcal{C}} := \bigcup_{t \in \mathbb{R}^+} \mathcal{C}_t$$

has the product structure $\mathcal{C}_t \times \mathbb{R}^+$, and its quotient is $\hat{\mathcal{Q}} := \hat{\mathcal{C}}/\mathcal{G} = \mathcal{Q}_t \times \mathbb{R}^+$. The rest of the proof then follows the standard argument. $\square$

### 6.2 Compactness, and the ends of the moduli

We have the following version of compactness results.

**6.2.1 Lemma** *Let $(g, w)$ be a good pair. Then:*

*(a) When $t = 0$, $\mathcal{M}_{g,w,t}$ is compact.*

*(b) When $t > 0$, the projection map $\partial_+ \times \partial_- : \mathcal{M}_{g,w,t}^n \to \mathrm{Sym}^n(\mathbb{C}) \times \mathrm{Sym}^n(\mathbb{C})$ is proper.*

*Proof.* We shall only prove the $t > 0$ case, since the $t = 0$ case is entirely similar, using Proposition 3.3.11 instead of Proposition 3.3.3.

Let $t > 0$. We first assume that $w = 0$ in $\omega$. Suppose $\{c_1, c_2, \ldots\}$ is a sequence of elements in $\mathcal{M}^n$ such that $\partial_+ \times \partial_-(c_1), \partial_+ \times \partial_-(c_2), \ldots$ converges to $a \in \mathrm{Sym}^n(\mathbb{C}) \times \mathrm{Sym}^n(\mathbb{C})$. Given the $L^\infty$ bound on $\psi$ established in Lemma 3.1.5, by a standard elliptic bootstrapping argument (cf. [19, 31, 29]) one concludes that:

> *Over any compact subset $U$ of $M$, we may choose suitable representatives of $c_i$ in $\mathcal{C}$, denoted by $c_i^U$, so that a subsequence of $\{c_i^U\}$ converges over this subset.* (6.1)

Therefore: (1) one may choose representatives of $c_i$ in $\mathcal{C}$, denoted still by $c_i$, so that restricted to $M_{2R}$, $R \geq 2\mathfrak{R}$, a subsequence of $\{c_i\}$ converges in the $L_2^2$ norm. Furthermore, by Proposition 3.3.3 and (6.1) again, a diagonal argument shows that: (2) there is another subsequence of $\{c_i\}$ (without loss of generality, we assume it is $\{c_i\}$ itself) so that over $M\backslash M_R$ there are $e^{i\xi_i} \in \mathcal{G}$ such that $\{e^{i\xi_i} c_i\}$ converges in the $\hat{L}^2_{2;\epsilon}$-norm on $M\backslash M_R$. (Cf. Definition 5.1.8 for the notation).

By (1) and (2), there exists still another subsequence of $\{c_i\}$, so that $\{e^{i(1-\chi_R)\xi_i} c_i\}$ converges in the $\hat{L}^2_{2;\epsilon}$ norm over the whole $M$. The lemma is proved for $w = 0$.

Now when $w \neq 0$, we may choose a large enough $R$ so that $(g, w')$ is still good with $w' = \chi_R w$, and $\mathcal{M}_{g,w',t}$ is cobordant to $\mathcal{M}_{g,w,t}$. Now Proposition 3.3.3 still applies with $(g, w')$ because on $M\backslash M_{2R}$, $w' = 0$. The properness for $\mathcal{M}_{g,w',t}$ implies the properness for $\mathcal{M}_{g,w,t}$. $\square$

**6.2.2 Remark** (a) In fact, both $\partial_+$ and $\partial_-$ are proper in the $t > 0$ case. This follows from Corollary 3.3.9 (3), and an argument similar to the above proof. We shall henceforth call $\partial_+^{-1}(\mathrm{Sym}^n(\mathbb{C})\backslash K)$ "the end of $\mathcal{M}^n$" when

$$K = \{(z_1, z_2, \ldots, z_n) \in \mathrm{Sym}^n \mathbb{C} : z_i \in \mathbb{C}, \sum_{i=1}^n |z_i|^2 \leq R\}$$

for some large $R \gg \Re$. The properness of $\partial_+$ means that the complement of the end of the moduli is compact.

(b) In the $t = 0$ case, Lemma 6.2.1 together with Theorem 6.1.2 implies that one may define a version of Seiberg-Witten invariant of 3-manifolds (for $b_1 > 1$) as the $\mathbb{Z}_2$-count of the points in the moduli space. That this is independent of $g, w$ is guaranteed by Corollary 6.1.4. The story for the $t > 0$ case is more involved due to the noncompactness of the moduli. We need a more detailed description of the end of the moduli via gluing theory (cf. Theorem 6.3.3 below), to which we devote the rest of this section. In [23] section 7, the product structure of the end (Corollary 6.3.5 below) was used to construct a series of Seiberg-Witten invariants from the $t > 0$ moduli.

## 6.3 Local description of the moduli via gluing

The main goal of the rest of this section is to state and prove Theorem 6.3.3. This theorem describes parts of the moduli space for $t > 0$ Seiberg-Witten solutions by gluing solutions on $\mathbb{R}^3$ and some known solutions. This partially describes the higher-vortex-number strata in terms of the lower-vortex-number ones; in particular, by Lemma 6.2.1 this determines the ends of the moduli.

The result is obtained in four standard steps. The first three steps, (A), (B), (C) below, occupy the rest of this subsection. The last step, i.e. proving that the gluing map is a diffeomorphism, is in the next subsection.

For the rest of this section, let $t > 0$, and $(g, w)$ be a good pair.

(A) CONSTRUCTING THE APPROXIMATE SOLUTIONS.

This is done in the following steps.

1. Let $R$ be a real number; $R > 2\Re$. Let $C_1 > C_2 > C_3 > C_4$ be positive constants independent of $R$.

    Let $\gamma_R(x)$ be a smooth cutoff function which is $|z|$-dependent-only on $M\backslash M_\Re$, and is supported in $\{x \in M : |z| < C_2 R\}$ with value 1 on $\{x \in M : |z| < C_3 R\}$; $\gamma_R(x) := 1$ when $x \in M_\Re$.

    Let $y$ be a vortex solution of on $\mathbb{C}$ in the (vor)-gauge (cf. Appendix). Suppose its vortex number is $n$ ($n = 0$ is possible), whose centers fall in the region where $|z| > C_1 R$. Let $y^+(z^+, f) := (y(z^+), 0)$ in the decomposition (5.1) be a configuration on $U_+ \subset \mathbb{R}^3$.

    Let $\mathcal{N} \subset \mathcal{M}$ be a closed ball of radius $\rho$ about $[v_0] \in \mathcal{M}^m$, which may be represented as a subset in the gauge slice $\{v_0 + e : e \in T\mathcal{C}; d_{v_0}^* e = 0\} \subset \mathcal{C}$, where $v_0$ is in the $*$-gauge defined in the end of §4.3. Let $v$ be an element in the above subset. We require $\rho = \rho(R)$ to be small enough such that the centers of all such $v$ fall within the region where $|z| < C_4 R$, and that there exist a continuous family of right inverses $P_v$ of $\mathcal{D}_v$ with $\|P_v\| \leq C_5 \ \forall v$, for a constant $C_5$ independent of $v, R$. (Note that $P_v$ exists by the smoothness of $\mathcal{N}$, and it is uniformly bounded by the compactness of $\mathcal{N}$.)

2. Identify the spinor bundles over $U_+$ and $M$ on the support of $\gamma_R(1-\gamma_R)$ by requiring the $\alpha$-components of the spinor fields in $y^+$ and $v$ to be parallel. (Note that these components never vanish in this region by Proposition 3.3.3 and Appendix, point 3.)

    We define the configuration $c := y \# v := (\tilde{A}, \tilde{\psi}) \in \mathcal{C}$ as follows:
    (i) Let $c$ be $y^+$ on the region $M\backslash \text{support}(\gamma_R) \subset U_+$.
    (ii) On $M\backslash \text{support}(1 - \gamma_R)$, let $c = v$.
    (iii) On the support of $\gamma_R(1-\gamma_R)$, let

$$c := v + (1 - \gamma_R)(y^+ - v). \tag{6.2}$$

We will use $y$ to denote either a point in $\text{Sym}^n(\mathbb{C}\backslash D(R))$ or a corresponding vortex solution in the (vor)-gauge indiscriminantly. Similarly, we will use the same notation $v$ for an element in $\mathcal{N}$ or the corresponding Seiberg-Witten solution.

The above patching-up construction thus defines a smooth map

$$\# : \text{Sym}^n(\mathbb{C}\backslash D(R)) \times \mathcal{N} \to \mathcal{C}(M)$$

by $\#(y, v) := y\#v$. We denote by $\mathbf{Gl} := \text{Sym}^n(\mathbb{C}\backslash D(R)) \times \mathcal{N}$ "the space of gluing parameters".

(B) ERROR ESTIMATE AND UNIFORM INVERTIBILITY.

**6.3.1 Proposition** *Let $c := y \# v$, where $y$, $v$ are as in part (A) 1 above. Then there exists a $R_0$ such that for $R > R_0$, the map $\# : \mathbf{Gl} \to \mathcal{C}(M)$ is smooth and the following hold:*

1. $\|\mathcal{S}(y \# v)\|_{\mathcal{V}} \leq CR^{-3/2+\epsilon} + \|w\|_{2,1:\epsilon}$.

2. $\mathcal{D}_c$ *is right invertible for $c = y \# v$; $\mathrm{Ind}(\mathcal{D}_c) = 2n + \mathrm{Ind}(\mathcal{D}_v)$; and the right inverse $G_c$ is uniformly bounded with respect to $c$.*

*Proof.* The smoothness may be checked by routine estimates.

To verify assertion 1, write $v = (A', \psi')$, $y^+ = (A, \psi)$, and $(\tilde{A}, \tilde{\psi}) := y \# v$. We omit the subscript $R$ of the cutoff function and write it as $\gamma$ below.

Let's begin with the special case when $w = 0$.

It is easy to check that in this case

$$\bar{\partial}_{\tilde{A}} \tilde{\psi} = \rho(d\gamma)(\psi' - \psi) - \gamma(1-\gamma)\rho(A' - A)(\psi' - \psi) + (1-\gamma)o^3; \tag{6.3}$$
$$F_{\tilde{A}} - \rho^{-1} \circ \sigma(\tilde{\psi}, \tilde{\psi}) = d\gamma(A' - A) + (1-\gamma)o^3$$
$$+ \gamma(1-\gamma)[\sigma(\psi, \psi' - \psi) - \sigma(\psi', \psi' - \psi)].$$

Thus we see that the two $L_1^2(\mathbb{C}, iN \oplus \mathbb{C} \otimes T\mathbb{C})$ components of $\mathcal{S}(y \# v)$ (according to the decomposition of $\mathcal{V}$ in (5.3)) are bounded by

$$\|d\gamma(\partial_\pm v - \partial_\pm y^+)\|_{L_1^2(\mathbb{C}, T\mathbb{C} \oplus \mathbb{C})} + \|\gamma(1-\gamma)(\partial_\pm v - \partial_\pm y^+)^2\|_{L_1^2(\mathbb{C}, T\mathbb{C} \oplus \mathbb{C})},$$

which decays exponentially with $R$ for large $R$ by Appendix, point 3.

As for the $L_{1:\epsilon}^2(\mathcal{W}_0)$ component, a direct computation shows that it is bounded by

$$\|d\gamma(v - v_1)\|_{2,1:\epsilon} + 2\|\gamma(1-\gamma)(y^+ - v_1)(v - v_1)\|_{2,1:\epsilon} + \|\gamma(1-\gamma)(v - v_1)^2\|_{2,1:\epsilon}$$
$$+ 2\|\lambda_d \circ f(1 - \lambda_d \circ f)\gamma(1-\gamma)(\partial_+ v \circ z_+ - \partial_- v \circ z_-)^2\|_{2,1:\epsilon} + \|(1-\gamma)o^3\|_{2,1:\epsilon}$$

where $v_1(x) := \lambda_d \circ f \partial_+ v \circ z_+ + (1 - \lambda_d \circ f)\partial_- v \circ z_-$, with the usual interpretation. Recall from Proposition 3.3.3 that we have $|v - v_1| + |\nabla(v - v_1)| < C|x|^{-3}$ for large enough $|x|$. Together with Appendix pt. 6, this implies assertion 1.

When $w \neq 0$, $\|w\|_{2,1:\epsilon}$ is small, we may regard $w$ as a perturbation and use lemmas 6.3.2, 6.3.7 below to establish assertion 1.

To prove the right invertibility in assertion 2, notice that for large $R$, $f\big|_{M \setminus M_R}$ and $y^+$ can be extended to the whole $\mathbb{R}^3$ as 5.1.4 and Corollary 5.3.5. (We shall denote the extensions by the same notation.) From Corollary 5.3.5 (b), we know that $\mathcal{D}_{y^+}$ has a uniformly bounded right inverse. By assumption $\mathcal{D}_v$ also has a uniformly bounded right inverse. We denote the right inverses of $\mathcal{D}_{y^+}$ and $\mathcal{D}_v$ by $P_1, P_2$ respectively. Combine these $P_i$'s via partition of unity to construct an operator $Q$ similar to (5.38):

$$Q := \phi_1 P_1 \gamma_1 + \phi_2 P_2 \gamma_2, \tag{6.4}$$

where $\gamma_1 := \gamma$; $\gamma_2 := (1-\gamma)$, and $\phi_1, \phi_2$ are as usual, smooth cutoff functions with value 1 over the supports of $\gamma_1$, $\gamma_2$ respectively. Then just like the computations in section 5, we have $\mathcal{D}_c Q = 1 + \mathcal{R}$, where the operator $\mathcal{R}$ is given by the same formula (5.40) but with the cutoff function $\chi$ is replaced by $\gamma$ here. We see then that the operator norm of $\mathcal{R}$ goes to zero as $R \to \infty$, because $\|\sigma_\mathcal{D}(d\phi_i)\|_\infty \leq C_1/R$; $\|\gamma_1\gamma_2(y^+ - v)\|_\infty \leq C_2/R^3$, where $C_1$, $C_2$ are independent of $y$ and $v$. Therefore $1 + \mathcal{R}$ is invertible when $R$ is large enough, and we may take the desired right inverse of $\mathcal{D}_c$ to be $\mathcal{P}_c := Q(1 + \mathcal{R})^{-1}$, which is uniformly bounded because of the uniform boundedness of $P_1$ and $P_2$.

To prove the claim about the index of $\mathcal{D}_c$ in assertion 2, note that $\mathrm{Coker}(\mathcal{D}_c)$ is trivial since $\mathcal{D}_c$ is right-invertible. On the other hand, $\mathrm{Ker}(\mathcal{D}_c)$ is isomorphic to

$$\mathrm{Ker}(\mathcal{D}_y) \oplus \mathrm{Ker}(\mathcal{D}_v) \simeq \mathbb{R}^{2n} \oplus T\mathcal{N},$$

by mapping $(k_1, k_2) \in \mathrm{Ker}(\mathcal{D}_y) \oplus \mathrm{Ker}(\mathcal{D}_v)$ to

$$\gamma_1 k_1 + \gamma_2 k_2 - \mathcal{P}_c \Big( \sigma_\mathcal{D}(d\gamma)(k_1 - k_2) + \gamma_1\gamma_2(\mathcal{D}_y - \mathcal{D}_v)(k_1 - k_2) \Big).$$

$\square$

(C) PERTURBATION.

The key of this part is the next well-known lemma. It can be easily derived from the contraction mapping principle. Similar versions can be found in, for example, [2] Lemma 7.2.23 or [7].

**6.3.2 Lemma** *Let $\mathfrak{f}$ be a smooth map between two Banach spaces $E$ and $F$ which has an expansion as $\mathfrak{f}(q) = \eta + Lq + B(q)$ where $L$ is a linear map with a right inverse $G$, and $B$ is a quadratic form[6] satisfying the estimate:*

$$\|B(Gh_1) - B(Gh_2)\| \leq k(\|h_1\| + \|h_2\|)\|h_1 - h_2\| \tag{6.5}$$

*for any $h_1$, $h_2 \in F$ and some positive constant $k$. Then if $\|\eta\| \leq \frac{1}{10k}$, there exists a unique $h(\eta) \in F$, $\|h(\eta)\| \leq 1/(5k)$, such that $\mathfrak{f}(Gh(\eta)) = 0$. Moreover, $\|h(\eta)\| \leq 2\|\eta\|$, and $h(\eta)$ varies smoothly with $\eta$.*

**6.3.3 Theorem** *Let $t > 0$; $\mathcal{N}$ be as in Part (A). Then there exists $R'_0 \in \mathbb{R}^+$ such that for $R > R'_0$, and any $c \in \#(\mathbf{Gl})$, the following hold.*

*(1) There exists a unique $q(c) \in \hat{K}$ such that $\mathcal{S}(c + q(c)) = 0$, and $\|q(c)\|_{\hat{K}} \leq C\|\mathcal{S}(c)\|_\mathcal{V} \leq C' R^{-3/2+\epsilon}$. Thus this defines a map $\Upsilon : \mathrm{Sym}^n(\mathbb{C}\backslash D(R)) \times \mathcal{N} \to \mathcal{M}^{n+m}$,*

$$\Upsilon(y, v) := [y \# v + q(y \# v)] \in \mathcal{Q}.$$

---
[6] We will use the same notation for its associated bilinear form.

*(2)* $\Upsilon$ *is smooth and injective. Actually, it describes the local structure of the moduli space in the following sense. Let* $\mathcal{U} \subset \mathrm{Sym}^n(\mathbb{C}\backslash D(R)) \times \mathcal{N}$ *be an open subset, and let*

$$U(d) := \Big\{c' \in \mathcal{C} : \exists c \in \#(\mathcal{U}), \|c' - c\|_{\hat{K}} < d; \|\mathcal{S}(c')\|_{\mathcal{V}} < C_{\mathcal{D}}d\Big\},$$

*where* $C_{\mathcal{D}}$ *is a constant larger than the operator norm of* $\mathcal{D}$. *Then there exists a small positive number* $\nu$ *depending on* $R$ *and* $\mathcal{U}$, *such that* $U(\nu) \cap \mathcal{M}^n$ *is diffeomorphic to a neighborhood of* $\mathcal{U}$ *in* $\mathrm{Sym}^n(\mathbb{C}\backslash D_R) \times \mathcal{N}$.

**6.3.4 Remark** If in the pre-gluing construction in Part (A), $v$ is allowed to have centers in $\{x \in M : |z| > C_1 R\}$ in addition to those in $\{x \in M : |z| < C_4 R\}$, then part (1) of the above theorem still holds, and the same arguments define a smooth map $\Upsilon$, which however may no longer be injective.

$\Upsilon$ is usually called the *gluing map*.

**6.3.5 Corollary** *(1) The end of the vortex number 1 stratum of the moduli space for any MEE consists of finite copies of* $\mathbb{C}\backslash D(R)$.

*(2) Any point in the end of* $\mathcal{M}^n$ *(cf. Remark 6.2.2) has a neighborhood with a product structure* $\mathcal{N} \times \mathcal{R}$, *where* $\mathcal{N}$ *is a neighborhood in* $\mathcal{M}^p$, $p < n$, *and* $\mathcal{R}$ *is a neighborhood in* $\mathrm{Sym}^{n-p}(\mathbb{C})$. *In particular, any non-compact component of the moduli space has dimension at least 2.*

*Proof.* (1) Lemma 6.2.1 says that the vortex number 0 stratum consists of finite points. Taking this as $\mathcal{N}$ and taking $\mathcal{U}$ to be $(\mathbb{C}\backslash D(R_1)) \times \mathcal{N}$, $R_1 > R$, Theorem 6.3.3 shows that the end of $\mathcal{M}^1$ contains $(\mathbb{C}\backslash D(R_1)) \times \mathcal{N}$. On the other hand, this already describes the whole end for the following reasons. By Lemma 6.2.1, any Seiberg-Witten solution corresponding to an element in the end have a single center far away from $M_{R_1}$. Let $v_1$ be such a Seiberg-Witten solution on $M$ with vortex number 1, and let $y_0$ be the unique vortex number 0 solution on $\mathbb{R}^3$. By Remark 6.3.4, $y_0 \# v_1$ approximates a vortex number 0 solution on M, say $v_0$. Thus any such $v_1$ approximates a vortex number 0 solution $v_0$ inside the cylinder $|z| = R$. This together with the proof of Proposition 6.3.1 implies that $v_1 \in U(\nu)$ and by Theorem 6.3.3 (2), this means that $v_1 \in \Upsilon((\mathbb{C}\backslash D(R_1) \times \mathcal{N})$ for all such $v_1$.

(2) Let $c$ be a configuration representing this point in $\mathcal{M}^n$. let $r_1, r_2, \ldots, r_n$ be the distances of the $n$ zeroes of $\partial_+ c$ to the origin; we order them such that $r_1 \leq r_2 \leq \cdots \leq r_n$. Since $c$ lies on the end of $\mathcal{M}^n$, we may assume $r_1^2 + r_2^2 + \cdots + r_n^2 \geq (10R)^2 n^3$, where $R > \Re$ is the same as that in Theorem 6.3.3. If $r_1 > R$, then by applying Theorem 6.3.3 and Remark 6.3.4 in the manner of (1), the claim holds for $p = 0$. Otherwise $\max_{i=1}^n (r_{i+1} - r_i) > R$; say the maximum occurs at $i = p$. Then Theorem 6.3.3 again shows that $c = \Upsilon(x, y)$, where $x \in \mathcal{M}^p$ and $y \in \mathrm{Sym}^n(\mathbb{C}\backslash D(R))$, and a neighborhood of $c$ is $\mathcal{N} \times \mathcal{R}$, where $\mathcal{N}$ is a neighborhood of $x$, and $\mathcal{R}$ is a neighborhood of $y$. $\square$

*Proof of Theorem 6.3.3 (1).* For the existence of the map $\Upsilon$, apply Proposition 6.3.1 and Lemma 6.3.2.

It is not hard to see by the proof of Lemma 4.1.8 that for any $h$ in $\hat{K}$ or $\mathcal{V}$, it can be uniquely decomposed as

$$h = q + d_c\gamma,$$

for some $q \in \operatorname{Ker} d_c^*$. Let the projection operator $\pi'$ be

$$\pi'(h) := q.$$

Note that (cf. §5.1) when $c$ is a Seiberg-Witten solution, $\operatorname{Ker}(d_c^*) = \operatorname{Im}(D_c)$. Here $c$ is only an "approximate solution"; however we shall show that $\pi'D_c\big|_{\operatorname{Ker}(d_c^*)}$ is still right-invertible). We take in Lemma 6.3.2

$$\mathfrak{f}(q) := \pi'\mathcal{S}(c+q) = \pi'\mathcal{S}(c) + \pi'D_c q + B(q). \tag{6.6}$$

**6.3.6 Lemma** $L_c := \pi'D_c$ *has a uniformly bounded right inverse, $G_c$. In other words, $L_c G_c = \operatorname{Id}$; $\pi_0 := G_c L_c$ defines a projection operator that decomposes $\hat{K} = \operatorname{Im} G_c \oplus \operatorname{Ker} L_c$.*

*Proof.* By Proposition 6.3.1, $\mathcal{D}_c$ has a uniformly-bounded right inverse $\mathcal{P}_c$. Decomposing its domain and range as

$$i\Omega^0 \oplus i\Omega^1 \oplus \Gamma(S) = i\Omega^0 \oplus \operatorname{Ker}(d_c^*) \oplus \operatorname{Im}(d_c),$$

$\mathcal{D}_c$ has the following form with respect to the decomposition:

$$\begin{pmatrix} 0 & 0 & d_c^* \\ 0 & L_c & * \\ d_c & * & * \end{pmatrix} \tag{6.7}$$

where $*$ denote some operators with small operator norms (which goes to 0 as $R \to \infty$; especially, these $*$ are replaced by 0 if $c$ is a Seiberg-Witten solution). The blocks $d_c$ and $d_c^*$ in (6.7) are isomorphisms. Thus $L_c$ has a uniformly bounded right inverse since $\mathcal{D}_c$ does. $\square$

On the other hand, denoting $q := (a, \eta)$,

$$B(q) := \pi'\Big(-i\sigma(\eta,\eta), \frac{\rho(a)}{2}\eta\Big).$$

A routine computation confirms that $B(q)$ satisfies the estimate (6.5) (using the bound on the operator norm of $G_c$). To apply Lemma 6.3.2, we now only need to estimate $\eta := \pi'\mathcal{S}(c)$, which is done in Proposition 6.3.1. Lemma 6.3.2 together with the following lemma then show that $c + q$ satisfies $\mathcal{S}(c+q) = 0$.

**6.3.7 Lemma** *If $\pi' \mathcal{S}(c+q) = 0$ and $q$ is small enough then $(1-\pi')\mathcal{S}(c+q) = 0$.*

*Proof.* Write $(1-\pi')\mathcal{S}(c+q) =: d_c\xi$. Then

$$d_c^* d_c \xi = d_c^* \mathcal{S}(c+q) = i\operatorname{Im}\langle \eta, \bar{\partial}_{A+a}(\psi + \eta)\rangle. \tag{6.8}$$

Now if $\pi' \mathcal{S}(c+q) = 0$, then $\mathcal{S}(c+q) - d_c\xi = 0$, and thus $\bar{\partial}_{A+a}(\psi + \eta) = \xi\psi$. Substitute this back to (6.8), we have

$$d_c^* d_c \xi = i\operatorname{Im}\langle \eta, \xi\psi\rangle.$$

This implies, via a straightforward extension of Lemma 4.1.8,

$$\|\xi\|_{\hat{X}} \leq C\|d_c^* d_c \xi\|_{\hat{L}^2_{1:\epsilon}} \leq C'\|\xi\|_{\hat{L}^2_{2:\epsilon}} \|\eta\|_{\hat{L}^2_{2:\epsilon}}$$

by the $L^\infty$-bounds on $\psi$ and $\nabla_A\psi$. Thus $\xi = 0$ if $\|\eta\|_{\hat{L}^2_{2:\epsilon}} \leq C_1\|q\|_{\hat{K}}$ is very small. □

The proof of Theorem 6.3.3 (1) is now complete; assertion (2) will be proved in the next subsection.

## 6.4 The gluing map is a diffeomorphism

This subsection contains the last step of the proof of Theorem 6.3.3: showing that the gluing map $\Upsilon$ is a diffeomorphism (assertion (2)).

First, we note that the smoothness of $\Upsilon$ follows from the smoothness of: (a) #, (b) $G_c$ with respect to $c$, (c) the term $h(\eta)$ with respect to $\eta$ in Lemma 6.3.2, which all may be verified by direct computations.

On the other hand, the contraction mapping theorem argument in Lemma 6.3.2 shows that $c + q(c)$ are the only Seiberg-Witten solutions in $U(d)$ having the form $c + G_c(h)$. To show that all Seiberg-Witten solutions in $U(d)$ are of this form, we need the following gauge-fixing lemma.

**6.4.1 A gauge-fixing result.** Since $\mathcal{U} \subset \operatorname{Sym}^n(\mathbb{C}\backslash D(R)) \times \mathcal{N}$ is open, there is a positive number $\nu$ such that $(\partial \#(\operatorname{Sym}^n(\mathbb{C}\backslash D(R)) \times \mathcal{N})) \cap U(\nu) = \emptyset$.

**Proposition** *Let $\nu$ be the positive number described above. Then every element in $U(\nu)$ has the form $u(c + G_c(h))$ for some $c \in \operatorname{Im}(\#)$, $u \in \mathcal{G}$, and $h \in \operatorname{Ker} d_c^* \subset \mathcal{V}$.*

*Proof.* We follow [2] Chapter 7 and use the continuity method to show that given an element $c' \in U(\nu)$, we may find $h$, $c$, and $u$ as in the statement of the proposition so that

$$G_c h = u^{-1}c' - c. \tag{6.9}$$

Since $c' \in U(\nu)$, we can find a $c(0) \in \#(\mathcal{U})$ such that $\|c' - c(0)\|_{\hat{K}} < \nu$. Let's look at a path in $U(\nu)$ interpolating $c'$ and $c(0)$:

$$c'(t) = tc' + (1-t)c(0). \tag{6.10}$$

Let $J \subset [0,1]$ be the set of all $t$ for which there exists $u(t) =: e^{\xi(t)} \in \mathcal{G}$, $c(t) \in \operatorname{Im} \# \subset \mathcal{C}(M)$ and $h(t) \in \mathcal{V}$, such that

$$e^{-\xi(t)} c'(t) = c(t) - G_{c(t)}(h(t)), \text{ and} \tag{6.11}$$
$$\|h\|_{\mathcal{V}} < \varepsilon, \ \nu \ll \varepsilon \ll 1. \tag{6.12}$$

Since $0 \in J$, if we can prove that $J$ is both closed and open, then $J = [0,1]$ and the proposition is proved.

To prove that $J$ is closed, we establish some a priori bounds. Since we have $h = \mathfrak{f}(Gh) - \mathfrak{f}(0) - B(Gh)$ from (6.6), by the Definition of $U(\nu)$ and Proposition 6.3.1 (1),

$$\|h\|_{\mathcal{V}} \leq C_{\mathcal{D}} \nu + C_1 R^{-3/2+\epsilon} + C' \|h\|_{\mathcal{V}}^2.$$

For $R^{-3/2+\epsilon} \ll \nu \ll \varepsilon$, this togther with (6.12) implies that $\|h\|_{\mathcal{V}} \leq \frac{\varepsilon}{2}$, which is a closed condition though we started with the open condition (6.12). Given a sequence in $J$, say $\{t_i\}$, this implies that there is a subsequence (also denoted $\{t_i\}$) such that the corresponding $\{h(t_i)\}$ has a weak limit in $\mathcal{V}$. On the other hand, both $\{c'(t_i)\}$ and $\{c(t_i)\}$ converge too (after further taking subsequences): the former because of the compactness of the parameter space $[0,1]$; the latter, $\{c(t_i)\}$, converges because (6.11) implies that the centers of $c(t)$ are bounded, and because of the local compactness of the space of gluing parameters **Gl**. Lastly, $\{u(t_i)\}$ also has a weak limit in $\hat{X}$ by (6.11) and the convergence of $h(t_i)$, $c(t_i)$, $c'(t_i)$. Furthermore, the limits satisfy (6.11) again, which implies that $J$ is closed.

To show the openness, consider the linearization at $t$ of $c'(t) = u(t)(c(t) + G_{c(t)} h(t))$ (without loss of generality, set $u(t) = 1$):

$$\mathcal{T}_{(u,c,h)}(\xi, \tau, \zeta) = d_c \xi + l(\tau) + (\partial_\tau G_c) h + G_c \zeta + \xi G_c h, \tag{6.13}$$

where $\xi \in i\Omega^0(M)$, $\tau \in T\,\mathbf{Gl}$; $l$ denotes the linearization of the gluing map $\#$, and $\zeta \in \operatorname{Ker} d_c^* \subset \mathcal{V}$. Note that our choice of $\nu$ makes sure that for all $t \in J$, $c(t) \notin \partial \operatorname{Im}(\#)$.

We shall show that $\mathcal{T}$ is surjective, which implies that $J$ is open.

We begin with some definitions for $\xi, \zeta, \tau$. Define the norms

$$\|(\xi, \zeta)\|_{B_1} := \|\xi\|_{\hat{X}} + \|\zeta\|_{\mathcal{V}};$$
$$\|(\xi, \tau, \zeta)\|_{B_2} := \|(\xi, \zeta)\|_{B_1} + |\tau|.$$

$\mathcal{T}$ is then a bounded linear operator between $B_2$ and $\hat{K}$.

**6.4.2 Lemma** *In the notation above,*

$$\mathcal{T}_1(\xi, \zeta) := d_c \xi + G_c \zeta \tag{6.14}$$

*is a Fredholm operator between $B_1$ and $\hat{K}$ with index $-\operatorname{Ind} \mathcal{D}_c$, and it has null kernel.*

*Proof.* The closed range property and the triviality of the kernel is guaranteed by the following two estimates: (i) By acting on both sides of (6.14) by $\pi' D_c$, we get

$$\|\zeta\|_{\mathcal{V}} \leq \|\pi' D_c \mathcal{T}_1(\xi, \zeta)\|_{\mathcal{V}} + \|\pi' D_c d_c \xi\|_{\mathcal{V}}$$
$$\leq C(\|\mathcal{T}_1(\xi, \zeta)\|_{\hat{K}} + R^{-3}\|\xi\|_{\mathcal{V}}),$$

where we have used estimates for $\|D_c d_c\|_\infty$ obtained from (6.3).

(ii) Also, by an extension of lemma 4.2.5, we have

$$\|\xi\|_{\hat{X}} \leq C\|d_c \xi\|_{\hat{L}^2_{2;\epsilon}}$$
$$\leq C'\|d_c \xi\|_{\hat{K}} \leq C'\|\mathcal{T}_1(\xi, \zeta)\|_{\hat{K}} + C'\|G_c \zeta\|_{\hat{K}}$$
$$\leq C'\|\mathcal{T}_1(\xi, \zeta)\|_{\hat{K}} + C''\|\zeta\|_{\mathcal{V}},$$

Combining the above two inequalities, we see then when $R$ is large

$$\|(\xi, \zeta)\|_{B_1} \leq C\|\mathcal{T}_1(\xi, \zeta)\|_{\hat{K}}.$$

Next we turn to the Fredholmness and the index calculation. This is easy due to the observation that if the domain and range of $\mathcal{T}_1$ are decomposed as $i\Omega^0 \oplus \operatorname{Ker} d_c^*$ and $\operatorname{Im} d_c \oplus \operatorname{Ker} d_c^*$ respectively, $\mathcal{T}_1$ has the following matrix form:

$$\begin{pmatrix} d_c & 0 \\ 0 & G_c \end{pmatrix},$$

where the block $d_c$ is an isomorphism. Thus $\mathcal{T}_1$ is Fredholm as $G_c$ is, and has the same index, namely $-\operatorname{Ind} \mathcal{D}_c$. □

Now to continue the proof of the proposition, decompose $\mathcal{T} =: \mathcal{T}_1 + \mathcal{T}_2 + \mathcal{T}_3 + \mathcal{T}_4$, where

$$\mathcal{T}_2(\tau) := l(\tau);$$
$$\mathcal{T}_3(\tau) := (\partial_\tau G_c)h;$$
$$\mathcal{T}_4(\xi) := \xi G_c h.$$

Since $h$ is small, both $\mathcal{T}_3$ and $\mathcal{T}_4$ can be regarded as perturbation. For example, $\mathcal{T}_3$ is small by the following observation. By the relationship between $G_c$ and $\mathcal{P}_c$, we can find a constant $C$ such that $\|(\partial_\tau G_c)h\|_{\hat{K}} \leq C\|(\partial_\tau \mathcal{P}_c)h\|_{\hat{K}^1}$; so we can try to estimate $\|(\partial_\tau \mathcal{P}_c)h\|_{\hat{K}^1}$ instead (for notation cf. Definition 5.1.8). $\partial_\tau(\mathcal{P}_c)h$ can be directly computed from the gluing construction of $\mathcal{P}_c$ in Proposition 6.3.1 (see the paragraph following (6.4), cf. also [2]), which enables us to find $\|(\partial_\tau \mathcal{P}_c)h\|_{\hat{K}^1} \leq C|\tau|\|h\|_\mathcal{V}$.

Note that $\mathcal{T}_2$ has null kernel by the gluing construction: that is, there exists a positive constant $C$ such that $\|(1-\pi_0)\mathcal{T}_2(\tau)\|_{\hat{K}} \geq C|\tau|$. Decompose $\mathcal{T}_2$ as $\mathcal{T}_2 = \pi_0 \mathcal{T}_2 + (1-\pi_0)\mathcal{T}_2$. By Proposition 6.3.1,

$$\|\pi_0 \mathcal{T}_2(\tau)\|_{\hat{K}} = \|G_c L_c l(\tau)\|_{\hat{K}} \leq C'\|L_c l(\tau)\|_\mathcal{V} \leq \zeta R^{-3/2+\epsilon}|\tau|, \tag{6.15}$$

is very small; therefore $(1-\pi_0)\mathcal{T}_2$ is the dominant term.

With these facts understood, we argue that $\mathcal{T}$ has null kernel. Suppose the contrary. Then there exists $(\xi, \tau, \zeta)$, $\tau \neq 0$ (since $\mathcal{T}_1$ has null kernel and $\mathcal{T}_3, \mathcal{T}_4$ are small perturbations), such that

$$\mathcal{T}(\xi, \tau, \zeta) = d_c \xi + G_c \zeta + l(\tau) + \mathcal{T}_3(\tau) + \mathcal{T}_4(\xi) = 0. \tag{6.16}$$

Without loss of generality, we assume $|\tau| = 1$. Since $\pi' D_c \mathcal{T}(\xi, \tau, \zeta) = 0$, we have

$$\pi' D_c d_c \xi + \zeta + \pi' D_c l(\tau) + \pi' D_c(\mathcal{T}_3(\tau) + \mathcal{T}_4(\xi)) = 0; \text{ so}$$

$$\|\zeta\|_\mathcal{V} \leq \varepsilon_1(|\tau| + \|\xi\|_{\hat{X}}). \tag{6.17}$$

Similarly, letting $d_c^*$ act on (6.16), we obtain

$$\|\xi\|_{\hat{X}} \leq \varepsilon\|(\xi, \tau, \zeta)\|_{B_2}. \tag{6.18}$$

Substituting (6.17), (6.18) into (6.16), we get

$$|\tau| \leq C\|l(\tau)\|_{\hat{K}} \leq \varepsilon < 1,$$

which contradicts with the assumption on $\tau$. Hence $\operatorname{Ker} \mathcal{T} = \emptyset$.

On the other hand, by Proposition 6.3.1, $\operatorname{Ind}(\mathcal{D}) = \dim \mathbf{Gl}$, and therefore $\mathcal{T}$ has index zero. This implies that $\mathcal{T}$ is surjective, and therefore $J$ is open and can not be anything but $[0, 1]$. End of the proof of the proposition. □

**6.4.3 Finishing the proof of Theorem 6.3.3: injectivity of $\Upsilon$.** Proposition 6.4.1 and the proof of Theorem 6.3.3 (1) in §6.3 show that $\Upsilon$ is an immersion that surjects over $\mathcal{M} \cap U(\nu)$. To show that $\Upsilon$ is actually a diffeomorphism, it thus remains to show that it is injective.

Suppose the contrary, then there exist $c_1, c_2 \in \operatorname{Im} \#$, $h_1 \in \operatorname{Ker} d_{c_1}^*$, $h_2 \in \operatorname{Ker} d_{c_2}^*$, $u \in \mathcal{G}$ such that

$$c_1 + G_{c_1} h_1 = u(c_2 + G_{c_2} h_2). \tag{6.19}$$

Since the terms $G_{c_1}h_1$, $G_{c_2}h_2$ are small (with respect to the $\hat{K}$ norm), (6.19) implies $c_1 \sim uc_2$. By construction of the # map, this means $u \sim 1$; $c_1 \sim c_2$. Thus we can approximate (6.19) by its linearization.

$$-l(\tau) = d_{c_2}\xi + \xi G_{c_2}h_2 + \partial_\tau(G_{c_2}h) + o(\xi, \tau), \tag{6.20}$$

where $\tau := c_2 - c_1$. To simplify notation, we will henceforth denote $c := c_2$; $h := h_2$.

From the proof of Proposition 6.4.1, we saw that

$$\|(\xi, \tau, \zeta)\|_{B_2} \leq C\|d_c\xi + G_c\zeta + l(\tau) + \mathcal{T}_3(\tau) + \mathcal{T}_4(\xi)\|_{\hat{K}}.$$

Letting $h = \zeta = 0$ in the above inequality and combined with (6.20), we have

$$\|\xi\|_{\hat{X}} + |\tau| \leq C\|d_c\xi + l(\tau)\|_{\hat{K}} \leq C'\|\partial_\tau(G_ch) + \xi G_ch\|_{\hat{K}} + o(\xi, \tau)$$
$$\leq C'\|\partial_\tau(G_ch)\|_{\hat{K}} + \varepsilon\|\xi\|_{\hat{X}} + o(\xi, \tau). \tag{6.21}$$

where $\varepsilon$ is a small number, and we have therefore

$$\|\xi\|_{\hat{X}} + |\tau| \leq C_2\|\partial_\tau(G_ch)\|_{\hat{K}}. \tag{6.22}$$

To estimate the right-hand side of the above, note that

$$\|\partial_\tau(G_ch)\|_{\hat{K}} \leq \|(\partial_\tau G_c)h\|_{\hat{K}} + \|G_c(\partial_\tau h)\|_{\hat{K}}.$$

As we have already estimated: $\|(\partial_\tau G_c)h\|_{\hat{K}} < \varepsilon'|\tau|$ in the proof of Proposition 6.4.1, we may concentrate on the term $\|G_c(\partial_\tau h)\|_{\hat{K}} \leq C\|\partial_\tau h\|_{\mathcal{V}}$. $\|\partial_\tau h\|_{\mathcal{V}}$ may be estimated by differentiating $\pi'\mathcal{S}(c) = h + \pi'B(G_ch, G_ch)$:

$$\partial_\tau \pi'\mathcal{S}(c) = \partial_\tau h + 2\pi'B\Big(G_c(\partial_\tau h) + (\partial_\tau G_c)h, G_ch\Big), \tag{6.23}$$

which implies $\|\partial_\tau h\|_{\mathcal{V}} \leq \varepsilon_2|\tau|$. Summing up, $\|\partial_\tau(G_ch)\|_{\hat{K}} \leq \varepsilon|\tau|$; this together with (6.22) gives $\|\xi\|_{\hat{X}} + |\tau| \leq \varepsilon_1|\tau|$, with $\varepsilon_1 \ll 1$. This implies $|\tau| = \|\xi\|_{\hat{X}} = 0$. That is, $\Upsilon$ is injective. This finishes the proof of Theorem 6.3.3. $\square$

## Appendix: Review of Vortex Solutions on $\mathbb{C}$

We recapitulate some useful properties of the vortex solutions on $\mathbb{C}$, which we apply extensively in this paper. More details may be found in [40, 15], and [43] section 2, [42] section 4.

Let $z$ be a complex coordinate on $\mathbb{C}$.

1. The vortex solutions on $\mathbb{C}$ are classified by the integer $n = \int_\mathbb{C} \frac{i}{2\pi} F_{A^E}$, usually called the "vortex number". This integer coincides with the number of points (counted with multiplicity) where the "Higgs field" $\alpha$ vanishes, and it is finite iff $F_A$, $\bar{\partial}_A\alpha$ are both $L^2$-integrable. We always assume this is true.

2. For a vortex solution $(A, \alpha)$, $\alpha$ has the following $L^\infty$-bound: $|\alpha| \leq 1$, where the equality happens only when $n = 0$, and in that case $|\alpha| \equiv 1$.

3. For a vortex solution $(A, \alpha)$, $|F_A|$, $(1 - |\alpha|^2)$ and $|\nabla_{A^E} \alpha|$ decay exponentially far away from the zeros of $\alpha$.

4. The gauge group $\mathrm{Map}(\mathbb{C}, U(1))$ acts on the space of vortex solutions as follows:
$$e^{i\xi}(A, \alpha) = (A - id\xi, e^{i\xi}\alpha).$$
The moduli space of vortex solutions under this gauge action is isomorphic to $\coprod_{n \in \mathbb{N}} \mathrm{Sym}^n \mathbb{C}$, where $n$ is the vortex number. The isomorphism is obtained by assigning a vortex solution the zeros of $\alpha$.

5. The moduli space of vortex solutions may be embedded as a gauge slice in the configuration space $\Omega^1(\mathbb{C}) \times \mathbb{C}^\infty(\mathbb{C}, \mathbb{C})$ such that the Higgs field $\alpha$ satisfies
$$\alpha = fz^n$$
when $z$ is large enough. Here $f$ is a positive function. In this paper we call this gauge (**vor**), and the image of the moduli space under this embedding $\mathcal{M}_{vortex}$.

6. The following operator on $i\Omega^1(\mathbb{C}) \times \mathbb{C}^\infty(\mathbb{C}, \mathbb{C})$ is the composition of the linearization of the vortex equation with a gauge fixing condition: for $c = (A', \alpha')$,
$$\Theta_c(b, \lambda) = (-4\partial b + \bar{\alpha}'\lambda, 2\bar{\partial}_{A'}\lambda + b\alpha'). \tag{A.1}$$
Here $(A', \alpha')$ is a solution to the vortex equation, and the gauge condition named is
$$\delta_c^1(b, \lambda) := 2d^* b + i \,\mathrm{Im}(\bar{\alpha}'\lambda) = 0. \tag{A.2}$$
$\Theta_c$ is a bounded Fredholm operator between $L_1^2$ and $L^2$, with null cokernel and $\mathrm{Ker}\,\Theta_c = \mathbb{C}^n$, $n$ being the vortex number of $c$. Furthermore, if $t \in \mathrm{Ker}\,\Theta_c$, then both $|t|$ and $|\nabla t|$ decay exponentially with $|z|$, for large enough $|z|$. $\mathrm{Ker}\,\Theta_c$ may be identified with the tangent space to the moduli space of vortex solutions at $c$, under the previously mentioned embedding. (Point 5 above.)

[41] C. Taubes, *The Seiberg-Witten and the Gromov invariants*, Math. Res. Lett., **2**, 221 (1995).

[42] C. Taubes, *From Seiberg-Witten to Gromov*, Journal of AMS **9**, 845 (1996).

[43] C. Taubes, *From Gromov to Seiberg-Witten*, J. Differential Geom. 51 (1999), no. 2, 203–334.

[44] C. Taubes, *Counting Connections and Curves*, J. Differential Geom. 52 (1999), no. 3, 453–609.

[45] C. Taubes, *The Seiberg-Witten invariants and symplectic forms*, Math. Res. Lett. **1**, 809 (1994).

[46] C. Taubes, *Counting pseudo-holomorphic submanifolds of dimension 4*, J. Diff. Geom., **44**, 818 (1996).

[47] C. Taubes, *Seiberg-Witten invariants and pseudo-holomorphic subvarieties for self-dual, harmonic 2-forms*, Geom. Topol. 3 (1999), 167–210

[48] C. Taubes, *Moduli spaces and Fredholm theory for pseudoholomorphic subvarieties associated to self-dual, harmonic* 2*-forms.* Asian J. Math. 3 (1999), no. 1, 275–324.

[49] C. Taubes, *The structure of pseudo-holomorphic subvarieties for a degenerate almost complex structure and symplectic form on* $S^1 \times B^3$. Geom. Topol. 2 (1998), 221–332

[50] C. Taubes, *The Seiberg-Witten invariants and 4-manifolds with essential tori*, Geom. Topol. 5 (2001) 441-519.

[51] V. Turaev, *Torsion invariants of* $\text{Spin}^c$*-structures on 3-danifolds*, Math. Res. Lett., **4**, 679 (1997).

[52] V. Turaev, *A combinatorial formulation for Seiberg-Witten invariants of 3-manifolds*, preprint, 1998.

[53] E. Witten, *Monopoles and four-manifolds*, Math. Res. Lett. **1**, 769 (1994).